\documentclass[aos,MSNbibl,nameyear,noautosecdot,dvips]{arximspdf}
\usepackage{mathbh}
\usepackage{dcolumn}
\usepackage{graphicx}

\doi{10.1214/13-AOS1175} 
\volume{42}
\issue{2}
\pubyear{2014}
\firstpage{413}
\lastpage{468}

\makeatletter
\newcommand{\widebar}{\overline}
\newcommand{\AR}{\mathit{AR}}
\newcolumntype{d}[1]{D{.}{.}{#1}}

\newcommand{\rright}{\right}
\newcommand{\lleft}{\left}
\newcommand{\eqref}[1]{(\ref{#1})}
\newtheorem{theorem}{Theorem}
\newtheorem{lemma}{Lemma}
\newproclaim{rem}{Remark}
\newcommand{\argmin}{\mathop{\operatorname{argmin}}}
\def\Var{\operatorname{Var}}
\def\Cov{\operatorname{Cov}}
\def\df{\operatorname{df}}
\def\rank{\operatorname{rank}}
\def\sign{\operatorname{sign}}
\def\supp{\operatorname{supp}}
\def\hy{\hat{y}}
\def\ty{\tilde{y}}
\def\hbeta{\hat{\beta}}
\def\tbeta{\tilde{\beta}}
\def\Exp{\operatorname{Exp}}
\def\RSS{\mathrm{RSS}}
\makeatother

\begin{document}
\begin{frontmatter}

\title{A significance test for the lasso\thanksref{T11}}
\runtitle{A significance test for the lasso}

\relateddois{T11}{Discussed in \relateddoi{d}{10.1214/13-AOS1175A},
\relateddoi{d}{10.1214/13-AOS1175B},
\relateddoi{d}{10.1214/13-AOS1175C},
\relateddoi{d}{10.1214/13-AOS1175D},
\relateddoi{d}{10.1214/13-AOS1175E} and
\relateddoi{d}{10.1214/14-AOS1175F};
rejoinder at \relateddoi{r}{10.1214/14-AOS1175REJ}.}

\begin{aug}
\author[A]{\fnms{Richard} \snm{Lockhart}\thanksref{t1}\ead[label=e1]{lockhart@sfu.ca}},
\author[B]{\fnms{Jonathan} \snm{Taylor}\thanksref{t2}\ead[label=e2]{jonathan.taylor@stanford.edu}},\\
\author[C]{\fnms{Ryan J.} \snm{Tibshirani}\thanksref{t3}\ead[label=e3]{ryantibs@cmu.edu}}
\and
\author[B]{\fnms{Robert} \snm{Tibshirani}\corref{}\ead[label=e4]{tibs@stanford.edu}\thanksref{t4}}
\runauthor{Lockhart, Taylor, Tibshirani and Tibshirani}
\affiliation{Simon Fraser University, Stanford University, Carnegie
Mellon University\\ and Stanford University}
\address[A]{R. Lockhart\\
Department of Statistics\\
\quad and Actuarial Science\\
Simon Fraser University\\
Burnaby, British Columbia V5A 1S6\\
Canada\\
\printead{e1}}
\address[B]{J. Taylor\\
Department of Statistics\\
Stanford University\\
Stanford, California 94305\\
USA\\
\printead{e2}\hspace*{36pt}}
\address[C]{R. J. Tibshirani\\
Departments of Statistics\\
\quad and Machine Learning\\
Carnegie Mellon University\\
229B Baker Hall\\
Pittsburgh, Pennsylvania 15213\\
USA\\
\printead{e3}}
\address[D]{R. Tibshirani\\
Department of Health, Research \& Policy\\
Department of Statistics\\
Stanford University\\
Stanford, California 94305\\
USA\\
\printead{e4}}
\end{aug}
\thankstext{t1}{Supported by the Natural Sciences and Engineering Research Council of Canada.}
\thankstext{t2}{Supported by NSF Grant DMS-12-08857 and AFOSR Grant 113039.}
\thankstext{t3}{Supported by NSF Grant DMS-13-09174.}
\thankstext{t4}{Supported by NSF Grant DMS-99-71405 and NIH Grant EB-001988.}
\pdftitle{A significance test for the lasso}

\received{\smonth{3} \syear{2013}}
\revised{\smonth{9} \syear{2013}}

%
\begin{abstract}
In the sparse linear regression setting, we consider testing the
significance of the predictor variable that enters the current
lasso model, in the sequence of models visited along the
lasso solution path.
We propose a simple test statistic based on lasso
fitted values, called the \textit{covariance test statistic},
and show that when the true model is linear, this
statistic has an $\operatorname{Exp}(1)$ asymptotic distribution under
the null
hypothesis (the null being that all truly active variables are
contained in the current lasso model). Our proof of this result for
the special case of the first predictor to enter the model (i.e.,
testing for a single significant predictor variable against the global
null) requires only weak assumptions on the predictor matrix $X$.
On the other hand, our proof for a general step in the lasso path
places further technical assumptions on $X$ and the generative
model, but still allows for the important high-dimensional case $p>n$,
and does not necessarily require that the current lasso model achieves
perfect recovery of the truly active variables.

Of course, for testing the significance of an additional variable
between two nested
linear models, one typically uses the chi-squared test,
comparing the drop in residual sum of squares (RSS) to a $\chi^2_1$
distribution. But when this additional variable is not fixed, and has
been chosen adaptively or greedily, this test is no
longer appropriate: adaptivity makes the drop in RSS stochastically
much larger than
$\chi^2_1$ under the null hypothesis. Our analysis explicitly accounts
for adaptivity,
as it must, since the lasso builds an adaptive \mbox{sequence} of linear
models as the tuning
parameter $\lambda$ decreases. In this analysis, shrinkage plays a key
role: though
additional variables are chosen adaptively, the \mbox{coefficients} of lasso
active variables
are shrunken due to the $\ell_1$ penalty. Therefore, the test
statistic (which is based on
lasso fitted values) is in a sense balanced by these two opposing
properties---adaptivity
and shrinkage---and its null distribution is tractable and
asymptotically $\operatorname{Exp}(1)$.
\end{abstract}

%
\end{frontmatter}

\setcounter{footnote}{5}

\section{\texorpdfstring{Introduction.}{Introduction}}\label{secintro}

We consider the usual linear regression setup, for an outcome vector
$y\in\mathbb{R}^n$ and matrix of predictor variables $X \in\mathbb
{R}^{n\times p}$:
%
%
\begin{equation}
\label{eqmodel} y = X\beta^* + \varepsilon,\qquad \varepsilon\sim N \bigl(0,
\sigma^2 I \bigr),
\end{equation}
where $\beta^* \in\mathbb{R}^p$ are unknown coefficients to be
estimated.
[If an intercept term is desired, then we can still assume a model of
the form \eqref{eqmodel} after centering $y$ and the columns of $X$;
see Section~\ref{secprostate} for more details.]
We focus on the lasso estimator [\citet{lasso}, \citet{bp}],
defined as
%
%
\begin{equation}
\label{eqlasso} \hbeta= \argmin_{\beta\in\mathbb{R}^p} \frac{1}{2}\| y-X\beta
\|_2^2 + \lambda\|\beta\|_1,
\end{equation}
where $\lambda\geq0$ is a tuning parameter, controlling the
level of sparsity in $\hbeta$.
Here, we assume that the columns
of $X$ are in general position in order to ensure uniqueness of the
lasso solution [this is quite a weak condition, to be discussed
again shortly; see also \citet{lassounique}].

There has been a considerable amount of recent work dedicated to the
lasso problem, both in terms of computation and theory. A
comprehensive summary of the
literature in either category would be too long for our purposes here, so
we instead give a short summary: for computational work, some
relevant contributions are
\citet{pco}, \citet{fista}, \citet{glmnet},
\citet{nesta}, \citet{admm}, \citet{tfocs};
and for theoretical work see, for example,
\citet{persistence}, \citet{fuchs}, \citet{cs},
\citet{nearopt}, \citet{irrepcond}, \citet{sharp},
\citet{nearideal}.
Generally speaking, theory for the lasso\vspace*{1pt} is focused on bounding the estimation
error $\|X\hbeta-X\beta^*\|_2^2$ or $\|\hbeta-\beta^*\|_2^2$, or ensuring
exact recovery of the underlying model, $\supp(\hbeta) = \supp(\beta^*)$
[with $\supp(\cdot)$ denoting the support function];
favorable results in both respects can be shown under
the right assumptions on the generative model \eqref{eqmodel} and
the predictor matrix $X$. Strong theoretical backing, as well
as fast algorithms, have made the lasso a highly popular tool.

Yet, there are still major gaps in our understanding of the lasso as an
estimation
procedure. In many real applications of the lasso, a practitioner will
undoubtedly
seek some sort of inferential guarantees for his or her computed lasso
model---but, generically, the usual constructs like $p$-values, confidence
intervals, etc., do not exist for lasso estimates. There is a small but growing
literature dedicated to inference for the lasso, and important
progress has certainly been made, with many methods being based on resampling or data splitting;
we review this work in Section~\ref{secrelated}. The current
paper focuses on a significance test for lasso models that does not
employ resampling or data splitting, but instead uses the full data
set as given, and proposes a test statistic that has a
simple and exact asymptotic null distribution.

Section~\ref{seccovtest} defines the problem that we are trying to solve,
and gives the details of our proposal---the covariance test statistic.
Section~\ref{secorthx} considers an orthogonal predictor
matrix $X$, in which case the statistic greatly simplifies. Here,
we derive its $\Exp(1)$ asymptotic distribution using
relatively simple arguments from extreme value theory. Section~\ref
{secgenx} treats a
general (nonorthogonal) $X$, and under some regularity conditions,
derives an $\Exp(1)$ limiting distribution for the covariance test
statistic, but through a different method of proof that relies on
discrete-time Gaussian
processes. Section~\ref{secsim} empirically verifies convergence of
the null distribution to $\Exp(1)$ over a variety of problem setups.
Up until this point, we have assumed that the error variance $\sigma^2$
is known; in Section~\ref{secsigma}, we discuss the case of
unknown~$\sigma^2$. Section~\ref{secrealdata} gives some real data
examples. Section~\ref{secextensions} covers extensions to
the elastic net, generalized linear models, and the Cox model
for survival data.
We conclude with a discussion in Section~\ref{secdiscussion}.

\section{\texorpdfstring{Significance testing in linear modeling.}{Significance testing in linear modeling}}\label{seccovtest}

Classic theory for significance testing in linear regression operates on
two fixed nested models. For example, if $M$ and $M\cup\{j\}$
are fixed subsets of $\{1,\ldots, p\}$, then to test the significance
of the $j$th predictor in the model (with variables in) $M\cup\{j\}$,
one naturally uses the chi-squared test, which computes the drop in
residual sum of squares (RSS) from regression on $M\cup\{j\}$ and $M$,
%
%
\begin{equation}
\label{eqchisq} R_j = (\RSS_M - \RSS_{M\cup\{j\}})/
\sigma^2
\end{equation}
and compares this to a $\chi_1^2$
distribution. (Here, $\sigma^2$ is assumed to be known; when $\sigma^2$
is unknown, we use the sample variance in its place, which results in
the $F$-test, equivalent to the $t$-test, for testing the significance of
variable~$j$.)

\begin{figure}

\includegraphics{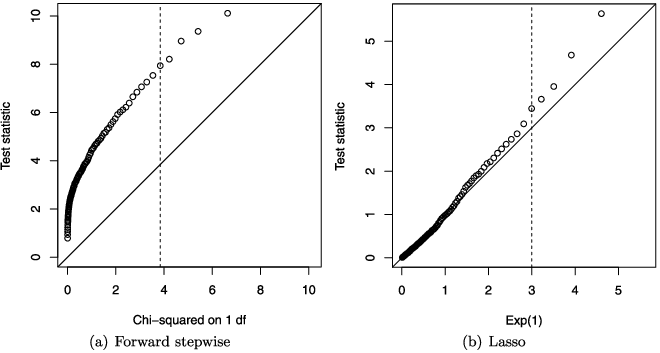}

\caption{A simple example with $n=100$ observations and $p=10$ orthogonal
predictors. All true regression coefficients are zero,
$\beta^*=0$. On the left is a quantile--quantile plot, constructed over
1000 simulations, of the standard
chi-squared statistic $R_1$ in \protect\eqref{eqchisq}, measuring the
drop in residual sum of squares for the first predictor to enter in
forward stepwise regression, versus the $\chi_1^2$ distribution.
The dashed vertical line marks the 95\% quantile of the $\chi_1^2$
distribution. The right panel shows
a quantile--quantile plot of the covariance test statistic $T_1$ in
\protect\eqref{eqcovtest} for the first predictor to enter in the
lasso path,
versus its asymptotic null distribution $\Exp(1)$. The covariance test
explicitly accounts for the adaptive nature of lasso modeling, whereas
the usual chi-squared test is not appropriate for adaptively selected
models, for example, those produced by forward stepwise regression.}
\label{fig1}
\end{figure}

Often, however, one would like to run the same test for $M$ and
$M\cup\{j\}$ that are not fixed, but the outputs of an adaptive or
greedy procedure. Unfortunately,
adaptivity invalidates the use of a $\chi_1^2$ null distribution for
the statistic \eqref{eqchisq}. As a simple example, consider forward
stepwise regression: starting with an empty model $M=\varnothing$, we
enter predictors one at a time, at each
step choosing the predictor $j$ that gives the largest
drop in
residual sum of squares. In other words, forward stepwise regression
chooses $j$ at each step in order to maximize $R_j$ in
\eqref{eqchisq}, over all $j\notin M$. Since
$R_j$ follows a $\chi_1^2$ distribution under the null
hypothesis for each fixed~$j$, the maximum possible $R_j$ will clearly
be stochastically larger than $\chi_1^2$ under the null. Therefore,
using a chi-squared test to evaluate the significance of a predictor
entered by forward stepwise regression would be far too liberal
(having type I error much larger than the nominal level). Figure~\ref
{fig1}(a) demonstrates this point by displaying the quantiles
of $R_1$ in forward stepwise regression (the chi-squared
statistic for the first predictor to enter) versus those
of a $\chi_1^2$ variate, in the fully null case (when $\beta^*=0$).
A test at the $5\%$ level, for example, using the $\chi^2_1$
cutoff of $3.84$, would have an actual type I error of about $39\%$.

%
%

The failure of standard testing methodology when applied to forward
stepwise regression is not an anomaly---in general, there seems to be
no direct way\vspace*{1pt} to carry out the significance tests designed for fixed
linear models in an adaptive setting.\footnote{It is important to
mention that a simple application of sample splitting can yield proper
$p$-values for an adaptive procedure like forward stepwise:
for example, run forward stepwise regression on one-half of the observations
to construct a sequence of models, and use the other half to evaluate
significance via the usual chi-squared test. Some of the related work
mentioned in Section~\ref{secrelated} does essentially this, but with
more sophisticated
splitting schemes. Our proposal uses the entire data set as given, and
we do not consider sample splitting or resampling techniques. Aside
from adding a layer of complexity, the use of sample splitting can
result in a loss of power in significance testing.}
Our aim is hence to provide a (new)
significance test for the predictor variables chosen adaptively by
the lasso, which we describe next.

\subsection{\texorpdfstring{The covariance test statistic.}{The covariance test statistic}}\label{seccovtestdef}

The test statistic that we propose here is constructed from the lasso
solution path, that is, the solution $\hbeta(\lambda)$ in \eqref{eqlasso}
a function of the tuning parameter $\lambda\in[0,\infty)$. The lasso
path can be computed by the well-known LARS algorithm of
\citet{lars} [see also \citeauthor{homotopy1} (\citeyear{homotopy1,homotopy2})],
which traces out the solution as $\lambda$ decreases from $\infty$ to~$0$. Note that when $\rank(X)<p$ there are possibly many lasso solutions
at each $\lambda$ and, therefore, possibly many solution paths; we assume
that the columns of $X$ are in general position,\footnote{Points
$X_1,\ldots, X_p\in\mathbb{R}^n$ are said to be in \textit{general position}
provided that no $k$-dimensional affine subspace $L \subseteq\mathbb{R}^n$,
$k < \min\{n,p\}$, contains more than $k+1$ elements of
$\{\pm X_1, \ldots, \pm X_p\}$, excluding
antipodal pairs. Equivalently: the affine span of any $k+1$ points
$s_1X_{i_1}, \ldots, s_{k+1}X_{i_{k+1}}$, for any signs
$s_1,\ldots, s_{k+1} \in\{-1,1\}$, does not contain any element of the set
$\{\pm X_i\dvtx i\neq i_1,\ldots, i_{k+1}\}$.} implying that there is
a unique
lasso solution at each $\lambda>0$, and hence a unique path. The assumption
that $X$ has columns in general position is a very weak one [much weaker,
e.g., than assuming that $\rank(X)=p$]. For example, if the entries of $X$
are drawn from a continuous probability distribution on $\mathbb
{R}^{np}$, then the
columns of $X$ are almost surely in general position, and this is true
regardless of the sizes of $n$ and $p$; see~\mbox{\citet{lassounique}}.

Before defining our statistic, we briefly review some properties of the
lasso path.
\begin{itemize}
\item The path $\hbeta(\lambda)$ is a continuous and piecewise linear
function of $\lambda$, with knots (changes in slope)
at values $\lambda_1 \geq\lambda_2 \geq\cdots\geq
\lambda_r \geq0$ (these knots depend on $y,X$).
\item At $\lambda=\infty$, the solution $\hbeta(\infty)$ has no active
variables (i.e., all variables have zero coefficients); for decreasing
$\lambda$, each knot $\lambda_k$ marks the entry or removal of some
variable from the current active set (i.e., its coefficient becomes
nonzero or zero, resp.). Therefore, the active set, and also the
signs of active coefficients, remain constant in between knots.
\item At any point $\lambda$ in the path, the corresponding active
set $A=\supp(\hbeta(\lambda))$ of the lasso solution indexes a
linearly independent set of predictor variables, that is, $\rank(X_A)=|A|$,
where we use $X_A$ to denote the columns of $X$ in~$A$.
\item For a general $X$, the number of knots in the lasso path is
bounded by~$3^p$ (but in practice this bound is usually very loose).
This bound comes from the following realization: if at some knot
$\lambda_k$, the active set is $A=\supp(\hbeta(\lambda_k))$ and the
signs of active coefficients are $s_A=\sign(\hbeta_A(\lambda_k))$,
then the active set and signs cannot again be $A$ and $s_{A}$ at some
other knot $\lambda_\ell\neq\lambda_k$. This in particular means
that once a variable enters the active set, it cannot
immediately leave the active set at the next step.
\item For a matrix $X$ satisfying the positive cone condition
(a~restrictive condition that covers, e.g., orthogonal matrices), there
are no variables removed from the active set as $\lambda$ decreases
and, therefore, the number of knots is~$p$.
\end{itemize}

We can now precisely define the problem that we are trying to
solve: at a given step in the lasso path (i.e., at a given knot), we
consider testing the significance of the variable that enters the
active set. To this end, we propose a test
statistic defined at the $k$th step of the path.

First, we define some needed quantities.
Let $A$ be the active set just
before $\lambda_k$, and suppose that predictor $j$ enters at $\lambda_k$.
Denote by $\hbeta(\lambda_{k+1})$ the solution at the next knot in the
path $\lambda_{k+1}$, using predictors $A\cup\{ j\}$. Finally, let
$\tbeta_A(\lambda_{k+1})$ be the solution of the lasso problem
using only the active predictors $X_A$, at $\lambda=\lambda_{k+1}$.
To be perfectly explicit,
%
%
\begin{equation}
\label{eqtbeta} \tbeta_A(\lambda_{k+1}) =
\argmin_{\beta_A \in\mathbb{R}^{|A|}} \frac{1}{2}\|y-X_A\beta_A
\|_2^2 + \lambda_{k+1} \|\beta_A
\|_1.
\end{equation}
We propose the \textit{covariance test statistic} defined by
%
%
\begin{equation}
\label{eqcovtest} T_k = \bigl( \bigl\langle y, X\hbeta(
\lambda_{k+1}) \bigr\rangle- \bigl\langle y,X_A
\tbeta_A(\lambda_{k+1}) \bigr\rangle \bigr) /
\sigma^2.
\end{equation}
%
Intuitively, the covariance statistic in \eqref{eqcovtest}
is a function of the difference between $X\hbeta$ and $X_A\tbeta_A$,
the fitted values given by incorporating the $j$th predictor into the current
active set, and leaving it out, respectively. These fitted values are
parameterized
by $\lambda$, and so one may ask: at which value of $\lambda$ should this
difference be evaluated? Well, note first that $\tbeta_A(\lambda
_k)=\hbeta_A(\lambda_k)$,
that is, the solution of the reduced problem at $\lambda_k$ is simply
that of the full
problem, restricted to the active set $A$ (as verified by
the KKT conditions). Clearly then, this means that we cannot evaluate
the difference
at $\lambda=\lambda_k$, as the $j$th variable has a zero coefficient
upon entry
at $\lambda_k$, and hence
\[
X\hbeta(\lambda_k)=X_A\hbeta_A(
\lambda_k)=X_A\tbeta_A(\lambda_k).
\]
Indeed, the natural choice for the tuning parameter in \eqref
{eqcovtest} is
$\lambda=\lambda_{k+1}$: this allows the $j$th coefficient to have
its fullest
effect on the fit $X\hbeta$ before the entry of the next variable at
$\lambda_{k+1}$ (or possibly, the deletion of a variable from $A$ at
$\lambda_{k+1}$).

Secondly, one may also ask about the particular choice of function of
$X\hbeta(\lambda_{k+1})-X_A\tbeta_A(\lambda_{k+1})$. The covariance
statistic in \eqref{eqcovtest} uses an inner product of this difference
with $y$, which can be roughly thought of as an (uncentered)
covariance, hence
explaining its name.\footnote{From its definition in \eqref
{eqcovtest}, we
get $T_k = \langle y-\mu, X\hbeta(\lambda_{k+1}) \rangle-
\langle y-\mu, X_A \tbeta_A(\lambda_{k+1}) \rangle+
\langle\mu, X\hbeta(\lambda_{k+1})-X_A \tbeta_A(\lambda
_{k+1})\rangle$ by
expanding $y=y-\mu+\mu$, with $\mu=X\beta^*$ denoting the true mean.
The first two terms are now really empirical covariances, and the last
term is
typically small. In fact, when $X$ is orthogonal, it is not hard to see
that this
last term is exactly zero under the null hypothesis.}
At a high level, the larger the covariance
of $y$ with $X\hbeta$ compared to that with $X_A\tbeta_A$, the more important
the role of variable $j$ in the proposed model $A\cup\{j\}$. There
certainly may be other functions that would seem appropriate here, but the
covariance form in \eqref{eqcovtest} has a distinctive advantage:
this statistic admits a simple and exact asymptotic null distribution.
In Sections~\ref{secorthx} and~\ref{secgenx}, we show that under the
null hypothesis that the current lasso model contains all truly
active variables, $A \supseteq\supp(\beta^*)$,
\[
T_k \stackrel{d} {\rightarrow}\Exp(1),
\]
that is, $T_k$ is asymptotically distributed as a standard exponential random
variable, given reasonable assumptions on $X$ and the magnitudes of the nonzero
true coefficients. [In some cases, e.g., when we have a strict
inclusion $A \supsetneq\supp(\beta^*)$, the use of an $\Exp(1)$
null distribution is actually conservative, because the limiting
distribution of $T_k$ is stochastically smaller than $\Exp(1)$.]
In the above limit, we are considering both $n,p \rightarrow\infty$;
in Section~\ref{secgenx}, we allow for the possibility $p>n$, the
high-dimensional case.


See Figure~\ref{fig1}(b) for a quantile--quantile plot of $T_1$ versus
an $\Exp(1)$ variate for the same fully null example ($\beta^*=0$)
used in Figure~\ref{fig1}(a); this shows that the weak convergence to
$\Exp(1)$
can be quite fast, as the quantiles are decently matched
even for $p=10$. Before proving this limiting distribution in
Sections~\ref{secorthx} (for an orthogonal $X$) and~\ref{secgenx}
(for a general $X$),
we give an example of its application to real data, and discuss issues related
to practical usage. We also derive useful alternative expressions
for the statistic, present a connection to degrees of freedom, review
related work, and finally, discuss the null hypothesis in more detail.

\subsection{\texorpdfstring{Prostate cancer data example and practical issues.}{Prostate cancer data example and practical issues}}\label{secprostate}

We consider a training set of 67 observations and 8 predictors, the goal
being to predict log of the PSA level of men who had surgery for prostate
cancer. For more details, see \citet{esl} and the
references therein.
Table~\ref{tabprostate} shows the results of forward stepwise regression
and the lasso. Both methods entered the same predictors in the same order.
The forward stepwise $p$-values are smaller than the lasso
$p$-values, and would enter four predictors at level $0.05$.
The latter would enter only one or maybe two predictors.
However, we know that the forward stepwise $p$-values are inaccurate,
as they are based on a null distribution that does not account
for the adaptive choice of predictors. We now make several remarks.

%
%
\begin{table}
\tabcolsep=0pt
\tablewidth=270pt
\caption{Forward stepwise and lasso applied to the prostate cancer
data example. The error variance is estimated by $\hat{\sigma}^2$,
the MSE of
the full model. Forward stepwise regression $p$-values are based on comparing
the drop in residual sum of squares (divided by $\hat{\sigma}^2$) to
an $F(1,n-p)$
distribution (using $\chi^2_1$ instead produced slightly smaller $p$-values).
The lasso $p$-values use a simple modification of the covariance test
\protect\eqref{eqcovtest} for unknown variance, given in
Section~\protect\ref{secsigma}. All $p$-values are rounded to 3
decimal places}\label{tabprostate}
\begin{tabular*}{\tablewidth}{@{\extracolsep{\fill}}@{}lccc@{}}
\hline
\textbf{Step} & \textbf{Predictor entered} & \textbf{Forward stepwise} & \textbf{Lasso}
\\
\hline
1 & lcavol & 0.000 & 0.000 \\
2 & lweight & 0.000 & 0.052 \\
3 & svi & 0.041 & 0.174 \\
4 & lbph & 0.045 & 0.929 \\
5 & pgg45 & 0.226 & 0.353 \\
6 & age & 0.191 & 0.650 \\
7 & lcp & 0.065 & 0.051 \\
8 & gleason & 0.883 & 0.978\\
\hline
\end{tabular*} \vspace*{-3pt}
\end{table}

%
\begin{rem}\label{rem1} The above example implicitly assumed that one might
stop entering variables into the model when the computed $p$-value rose
above some threshold. More generally, our proposed test statistic and
associated $p$-values could be used as the basis for multiple testing
and false discovery rate control methods for this problem; we leave
this to future work.
\end{rem}

%
\begin{rem}\label{rem2} In the example, the lasso entered a predictor
into the active set at each step. For a general $X$, however, a given
predictor variable may enter the active set more than once along the
lasso path, since it may leave the active set at some point. In this
case, we treat each entry as a separate problem. Our test is
specific to a step in the path, and not to a predictor variable at large.
\end{rem}

%
\begin{rem}\label{rem3} For the prostate cancer data set, it is
important to
include an intercept in the model. To accommodate this, we ran
the lasso on centered $y$ and column-centered $X$ (which is
equivalent to including an unpenalized intercept term in the lasso
criterion), and then applied the covariance test (with the centered
data).
In general, centering $y$ and the columns of $X$ allows us to
account for the effect of an intercept term, and still use a model
of the form \eqref{eqmodel}. From a theoretical perspective,
this centering step creates a weak dependence between the components
of the error vector $\varepsilon\in\mathbb{R}^n$. If originally we assumed
i.i.d. errors, $\varepsilon_i \sim N(0,\sigma^2)$, then after centering
$y$ and the columns of $X$, our new errors are of the form
$\tilde{\varepsilon}_i = \varepsilon_i - \bar{\varepsilon}$,
where $\bar{\varepsilon}=\sum_{j=1}^n \varepsilon_j/n$. It is easy see that
these new errors are correlated:
\[
\Cov(\tilde{\varepsilon}_i,\tilde{\varepsilon}_j) = -
\sigma^2/n\qquad\mbox{for } i \neq j.\vadjust{\goodbreak}
\]
One might imagine that such correlation would cause problems for our
theory in Sections~\ref{secorthx} and~\ref{secgenx}, which assumes
i.i.d. normal errors in the model \eqref{eqmodel}. However, a~careful look at the arguments in these sections reveals that the only
dependence on $y$ is through $X^T y$, the inner products of $y$ with
the columns of $X$. Furthermore,
\[
\Cov \bigl(X_i^T \tilde{\varepsilon},
X_j^T \tilde{\varepsilon} \bigr) = \sigma^2
X_i^T \biggl(I- \frac{1}{n}\mathbh{1}
\mathbh{1}^T \biggr) X_j = \sigma^2
X_i^T X_j\qquad\mbox{for all } i,j,
\]
which is the same as it would have been without centering (here
$\mathbh{1}\mathbh{1}^T$ is the matrix of all $1$s, and we used
that the columns of $X$ are centered). Therefore, our arguments in
Sections~\ref{secorthx} and~\ref{secgenx} apply equally well to
centered data, and centering has no effect on the asymptotic
distribution of $T_k$.
\end{rem}

%
\begin{rem}\label{rem4} By design, the covariance test is applied in a
sequential manner, estimating $p$-values for each predictor variable as
it enters the model along the lasso path. A more difficult problem is
to test the significance of any of the active predictors in a model
fit by the lasso, at some arbitrary value of the tuning parameter
$\lambda$. We discuss this problem briefly in Section~\ref{secdiscussion}.
\end{rem}

\subsection{\texorpdfstring{Alternate expressions for the covariance statistic.}{Alternate expressions for the covariance statistic}}\label{secalternate}

Here, we derive two alternate forms for the covariance statistic in
\eqref{eqcovtest}. The first lends some insight into the role
of shrinkage, and the second is helpful for the convergence results
that we
establish in Sections~\ref{secorthx} and~\ref{secgenx}. We rely on
some basic properties
of lasso solutions; see, for example, \citet{lassodf2},
\citet{lassounique}.
To remind the reader, we are assuming that $X$ has columns in general position.

For any fixed $\lambda$, if the lasso solution has active set
$A=\supp(\hbeta(\lambda))$ and signs $s_A=\sign(\hbeta_A(\lambda))$,
then it can be written explicitly (over active variables) as
\[
\hbeta_A(\lambda) = \bigl(X_A^T
X_A \bigr)^{-1} X_A^T y - \lambda
\bigl(X_A^T X_A \bigr)^{-1}
s_A.
\]
In the above expression, the first term $(X_A^T X_A)^{-1} X_A^T y$ simply
gives the regression coefficients of $y$ on the active variables $X_A$, and
the second term $-\lambda(X_A^T X_A)^{-1} s_A$ can be thought of as a shrinkage
term, shrinking the values of these
coefficients toward zero. Further, the lasso fitted value at $\lambda$ is
%
%
\begin{equation}
\label{eqlassofit} X\hbeta(\lambda) = P_A y - \lambda
\bigl(X_A^T \bigr)^+ s_A,
\end{equation}
where $P_A = X_A(X_A^T X_A)^{-1} X_A^T$ denotes the projection onto the
column space
of~$X_A$, and $(X_A^T)^+=X_A(X_A^T X_A)^{-1}$ is the (Moore--Penrose)
pseudoinverse
of $X_A^T$.

Using the representation \eqref{eqlassofit} for the fitted values, we can
derive our first alternate expression for the covariance statistic in
\eqref{eqcovtest}.
If $A$ and $s_A$ are the active set and signs just before the knot
$\lambda_k$, and $j$ is the variable\vadjust{\goodbreak} added to the active set at
$\lambda_k$, with sign~$s$ upon entry, then by
\eqref{eqlassofit},
\[
X\hbeta(\lambda_{k+1}) = P_{A\cup\{j\}} y - \lambda_{k+1}
\bigl(X_{A\cup\{j\}}^T \bigr)^+ s_{A\cup\{j\}},
\]
where $s_{A\cup\{j\}}=\sign(\hbeta_{A\cup\{j\}}(\lambda_{k+1}))$.
We can equivalently write $s_{A\cup\{j\}}=\break (s_A,s)$, the concatenation
of $s_A$ and the sign $s$ of the $j$th coefficient when it entered
(as no sign changes could have occurred inside of the
interval $[\lambda_k,\lambda_{k+1}]$, by definition of the knots). Let
us assume for the moment that the solution of reduced lasso problem
\eqref{eqtbeta} at $\lambda_{k+1}$ has all variables active and
$s_A = \sign(\tbeta_A(\lambda_{k+1}))$---remember,
this holds for the reduced problem at $\lambda_k$, and we will return
to this assumption shortly. Then, again by \eqref{eqlassofit},
\[
X_A\tbeta_A(\lambda_{k+1}) = P_A
y - \lambda_{k+1} \bigl(X_A^T \bigr)^+
s_A
\]
and plugging the above two expressions into \eqref{eqcovtest},
%
%
\begin{eqnarray}\label{eqcovtestproj}
T_k &=& y^T (P_{A\cup\{j\}}-P_A)
y / \sigma^2
\nonumber\\[-8pt]\\[-8pt]
&&{}  - \lambda_{k+1} \cdot y^T \bigl(
\bigl(X_{A\cup\{j\}}^T \bigr)^+ s_{A\cup\{j\}} -
\bigl(X_A^T \bigr)^+ s_A \bigr) /
\sigma^2.\nonumber
\end{eqnarray}
Note that the first term above is
$y^T (P_{A\cup\{j\}}-P_A) y / \sigma^2 =
(\|y-P_A y\|_2^2 - \|y-P_{A\cup\{j\}}y\|_2^2)/\sigma^2$,
which is exactly the chi-squared statistic for testing the significance
of variable $j$, as in \eqref{eqchisq}. Hence, if $A,j$ were fixed, then
without the second term, $T_k$ would have a $\chi_1^2$ distribution under
the null. But of course $A,j$ are not fixed, and so much like we saw
previously with forward stepwise regression, the first term in \eqref
{eqcovtestproj}
will be generically larger than $\chi_1^2$, because $j$ is chosen
adaptively based on
its inner product with the current lasso residual vector.
Interestingly, the second
term in \eqref{eqcovtestproj} adjusts for this adaptivity: with this
term, which is composed
of the shrinkage factors in the solutions of the two relevant lasso
problems (on $X$~and~$X_A$), we prove in the coming sections that $T_k$ has an asymptotic
$\Exp(1)$ null distribution. Therefore, the presence of the second term
restores the (asymptotic) mean of $T_k$ to $1$, which is what it would
have been
if $A,j$ were fixed and the second term were missing. In short,
adaptivity and
shrinkage balance each other out.

This insight aside, the form \eqref{eqcovtestproj} of the covariance statistic
leads to a second representation that will be useful for the
theoretical work
in Sections~\ref{secorthx} and~\ref{secgenx}. We call this the \textit
{knot form}
of the covariance statistic, described in the next lemma.

%
%
\begin{lemma}
\label{lemcovtestknot}
Let $A$ be the active set just before the $k$th step in the lasso
path, that is, $A=\supp(\hbeta(\lambda_k))$, with $\lambda_k$ being the
$k$th knot. Also, let $s_A$ denote\vspace*{1pt} the signs of the active coefficients,
$s_A=\sign(\hbeta_A(\lambda_k))$, $j$ be the predictor that enters the
active set at $\lambda_k$, and $s$ be its sign upon entry. Then,
assuming that
%
%
\begin{equation}
\label{eqnosignchange} s_A = \sign \bigl(\tbeta_A(
\lambda_{k+1}) \bigr)
\end{equation}
or in other words, all coefficients are active in the reduced lasso problem
\eqref{eqtbeta} at $\lambda_{k+1}$ and have signs $s_A$, we have
%
%
\begin{equation}
\label{eqcovtestknot} T_k = C(A,s_A,j,s) \cdot
\lambda_k (\lambda_k - \lambda_{k+1}) /
\sigma^2,
\end{equation}
where
\[
C(A,s_A,j,s) = \bigl\| \bigl(X_{A\cup\{ j\}}^T \bigr)^+
s_{A \cup\{j\}} - \bigl(X_A^T \bigr)^+ s_A
\bigr\|_2^2
\]
and $s_{A\cup\{j\}}$ is the concatenation of $s_A$ and $s$.
\end{lemma}

The proof starts with expression \eqref{eqcovtestproj}, and arrives at
\eqref{eqcovtestknot} through simple algebraic manipulations. We
defer it until
Appendix~\ref{appcovtestknot}.

When does the condition \eqref{eqnosignchange} hold? This was a key assumption
behind both of the forms \eqref{eqcovtestproj} and \eqref
{eqcovtestknot} for the
statistic. We first note that
the solution $\tbeta_A$ of the reduced lasso problem has signs $s_A$
at $\lambda_k$,
so it will have the same signs $s_A$ at $\lambda_{k+1}$ provided that
no variables are
deleted from the active set in the solution path $\tbeta_A(\lambda)$
for $\lambda\in[\lambda_{k+1},\lambda_k]$. Therefore, assumption
\eqref{eqnosignchange} holds:
\begin{longlist}[{[3]}]
\item[{[1]}] When $X$ satisfies the positive cone condition (which includes
$X$ orthogonal),
because no variables ever leave the active set in this case. In fact,
for $X$ orthogonal,
it is straightforward to check that $C(A,s_A,j,s)=1$, so
$T_k=\lambda_k(\lambda_k-\lambda_{k+1})/\sigma^2$.
\item[{[2]}] When $k=1$ (we are testing the first variable to enter), as a variable
cannot leave the active set right after it has entered. If $k=1$ and $X$
has unit normed columns, $\|X_i\|_2=1$ for $i=1,\ldots, p$, then we again
have $C(A,s_A,j,s)=1$
(note that $A=\varnothing$), so $T_1=\lambda_1(\lambda_1-\lambda
_2)/\sigma^2$.
\item[{[3]}] When $s_A = \sign((X_A)^+ y)$, that is, $s_A$ contains the signs
of the least squares
coefficients on $X_A$, because the same active set and signs cannot
appear at two different knots in the lasso path (applied here to the
reduced lasso
problem on~$X_A$).
\end{longlist}
The first and second scenarios are considered in Sections~\ref
{secorthx} and
\ref{secgenxfirst}, respectively. The third scenario is actually
somewhat general
and occurs, for example, when $s_A = \sign((X_A)^+ y) = \sign(\beta
_A^*)$; in this case, both the lasso
and least squares on $X_A$ recover the signs of the true coefficients.
Section~\ref{secgenxgen} studies the general $X$ and $k \geq1$ case,
wherein this
third scenario is important.

\subsection{\texorpdfstring{Connection to degrees of freedom.}{Connection to degrees of freedom}}\label{secdf}

There is an interesting connection between the covariance
statistic in \eqref{eqcovtest} and the degrees of freedom of a fitting
procedure. In the regression setting \eqref{eqmodel}, for an estimate
$\hy$
[which we think of as a fitting procedure $\hy=\hy(y)$],
its degrees of freedom is typically defined [\citet{bradbiased}] as
%
%
\begin{equation}
\label{eqdf} \df(\hy)=\frac{1}{\sigma^2}\sum_{i=1}^n
\Cov(y_i,\hy_i).
\end{equation}
In words, $\df(\hy)$ sums the covariances of each observation $y_i$
with its fitted value~$\hy_i$. Hence, the more adaptive a fitting
procedure, the higher this covariance, and the greater its degrees of
freedom. The covariance test evaluates the significance
of adding the $j$th predictor via something loosely like a sample
version of degrees of freedom, across two models: that fit on
$A\cup\{j\}$, and that on $A$. This was more or~less the inspiration
for the current work.

Using the definition \eqref{eqdf}, one can reason [and confirm by
simulation, just as in Figure~\ref{fig1}(a)] that with $k$
predictors entered into the model, forward stepwise regression had
used substantially more than $k$ degrees of freedom.
But something quite remarkable happens when we consider the lasso:
for a model containing $k$ nonzero coefficients, the degrees of
freedom of the lasso fit is equal to
$k$ (either exactly or in expectation, depending on the assumptions)
[\citet{lars}, \citet{lassodf}, \citet{lassodf2}].
Why does this happen? Roughly speaking, it is the same adaptivity versus
shrinkage phenomenon at play. [Recall our discussion in the last section
following the expression \eqref{eqcovtestproj} for the covariance statistic.]
The lasso adaptively chooses the active predictors, which costs extra
degrees of freedom; but it also shrinks the nonzero coefficients
(relative to the usual least squares estimates), which decreases the
degrees of freedom just the right amount, so that the total is simply $k$.
\subsection{\texorpdfstring{Related work.}{Related work}}\label{secrelated}

There is quite a lot of recent work related to the
proposal of this paper.
\citet{screenclean} propose a procedure for variable selection
and $p$-value estimation in high-dimensional linear models based on
sample splitting, and this idea was extended by \citet{mein2009}.
\citet{stabselect} propose a generic method using resampling
called ``stability selection,'' which controls the expected
number of false positive variable selections. \citet{minnier2011}
use perturbation resampling-based procedures to approximate the
distribution of a general class of penalized parameter
estimates. One big difference with the work here: we propose a
statistic that utilizes the data as given and does not employ any
resampling or sample splitting.

\citet{zhangconf} derive confidence intervals for
contrasts of high-dimensional regression coefficients,
by replacing the usual score vector with the residual
from a relaxed projection (i.e., the residual from sparse linear
regression). \citet{buhlsignif}
constructs $p$-values for coefficients in high-dimensional regression
models, starting with ridge estimation and then employing a bias
correction term that uses the lasso.
Even more recently,
\citet{vdgsignif}, \citeauthor{montahypo2} (\citeyear{montahypo2,montahypo1})
all present approaches for debiasing the lasso
estimate based on estimates of the inverse covariance matrix of the
predictors. (The latter work focuses on the
special case of a predictor matrix $X$ with i.i.d. Gaussian rows; the
first two consider a general matrix $X$.) These debiased lasso
estimates are asymptotically normal, which allows one to compute
$p$-values both marginally for an individual coefficient, and
simultaneously for a group of coefficients. All of the work mentioned
in the present paragraph provides a way to make inferential
statements about preconceived
predictor variables of interest (or preconceived groups of interest);
this is in contrast to our work, which instead deals directly with
variables that have been adaptively selected by the lasso procedure.
We discuss this next.

\subsection{\texorpdfstring{What precisely is the null hypothesis?}{What precisely is the null hypothesis}}\label{secwhat?}

The referees of a preliminary version of this manuscript expressed
some confusion with regard to the null distribution considered by
the covariance test. Given a fixed number of steps $k \geq1$
along the lasso path, the covariance test examines the
set of variables $A$ selected by the lasso before the $k$th step
(i.e., $A$ is the current active set not including the variable
to be added at the $k$th step). In particular, the null
distribution being tested is
%
%
\begin{equation}
\label{eqnullhypo1} H_0\dvtx A \supseteq\supp \bigl(\beta^* \bigr),
\end{equation}
where $\beta^*$ is the true underlying coefficient vector in the model
\eqref{eqmodel}. For $k=1$, we have $A=\varnothing$ (no
variables are selected before the first step), so this reduces to a
test of the global null hypothesis: $\beta^*=0$. For $k > 1$, the
set $A$ is random (it depends on $y$), and hence the null hypothesis
in \eqref{eqnullhypo1} is itself a random event. This makes the
covariance test a \textit{conditional hypothesis test} beyond the first
step in the path, as the null hypothesis that it considers is
indeed a function of the observed data. Statements about its null
distribution must therefore be made conditional on the event that $A
\supseteq\supp(\beta^*)$, which is
precisely what is done in Sections~\ref{secorthxgen} and
\ref{secgenxgen}.

Compare the null hypothesis in \eqref{eqnullhypo1} to a null
hypothesis of the form
%
%
\begin{equation}
\label{eqnullhypo2} H_0\dvtx S \cap\supp \bigl(\beta^* \bigr) =
\varnothing,
\end{equation}
where $S \subseteq\{1,\ldots, p\}$ is a fixed subset. The latter
hypothesis, in \eqref{eqnullhypo2}, describes the setup considered
by \citet{zhangconf}, \citet{buhlsignif},
\citet{vdgsignif}, \citeauthor{montahypo2} (\citeyear{montahypo2,montahypo1}). At face value, the hypotheses
\eqref{eqnullhypo1} and \eqref{eqnullhypo2} may appear similar
[the test in \eqref{eqnullhypo1} looks just like that in
\eqref{eqnullhypo2} with $S=\{1,\ldots, p\}\setminus A$], but they are
fundamentally very different. The difference is that the null
hypothesis in \eqref{eqnullhypo1} is random, whereas that in
\eqref{eqnullhypo2} is fixed; this makes the covariance test
a conditional hypothesis test, while the tests constructed in all of
the aforementioned work are traditional (unconditional) hypothesis
tests. It should be made clear that the goal of our work and
these works also differ.
Our test examines an adaptive subset of variables $A$
deemed interesting by the lasso procedure; for such a goal, it
seems necessary to consider a random null hypothesis, as theory
designed for tests of fixed hypotheses would not be valid
here.\footnote{In principle, fixed hypothesis tests can be used along
with the appropriate correction for multiple comparisons in order to test
a random null hypotheses. Aside from being conservative, it is
unclear how to efficiently carry out such a procedure when the random
null hypothesis consists of a group of coefficients (as opposed to a
single one).}
The main goal of \citet{zhangconf}, \citet{buhlsignif},
\citet{vdgsignif}, \citeauthor{montahypo2} (\citeyear{montahypo2,montahypo1}), it appears, is to construct a new set
of variables, say $\widetilde{A}$, based on testing the hypotheses in
\eqref{eqnullhypo2} with $S=\{j\}$ for $j=1,\ldots, p$. Though the
construction of this new set $\widetilde{A}$ may have started from a lasso
estimate, it need not be true that $\widetilde{A}$ matches the lasso
active set $A$, and ultimately it is this new set $\widetilde{A}$ (and
inferential statements concerning $\widetilde{A}$) that these authors
consider the point of interest.


\section{\texorpdfstring{An orthogonal predictor matrix $X$.}{An orthogonal predictor matrix $X$}}\label{secorthx}

We examine the special case of an orthogonal predictor matrix
$X$, that is, one that satisfies $X^T X = I$. Even though the
results here can be seen as special cases of those for
a general $X$ in Section~\ref{secgenx}, the arguments in the current
orthogonal $X$ case rely on relatively straightforward extreme value theory
and are hence much simpler
than their general $X$ counterparts (which analyze the knots in the lasso
path via Gaussian process theory). Furthermore, the $\Exp(1)$ limiting
distribution for the covariance statistic translates in the orthogonal
case to a few interesting and
previously unknown (as far as we can tell) results on the order
statistics of
independent standard $\chi_1$ variates. For these reasons, we discuss the
orthogonal $X$ case in detail.

As noted in the discussion following Lemma~\ref{lemcovtestknot} (see the
first point), for an orthogonal $X$, we know that the covariance
statistic for testing the entry of the variable at step $k$ in the
lasso path is
\[
T_k = \lambda_k(\lambda_k-
\lambda_{k+1})/\sigma^2.
\]
Again using orthogonality, we rewrite $\|y-X\beta\|_2^2 =
\|X^T y - \beta\|_2^2 + C$ for a constant $C$ (not depending
on $\beta$) in the criterion in \eqref{eqlasso}, and then we can
see that the lasso solution at any given value of $\lambda$ has
the closed-form:
\[
\hbeta_j(\lambda) = S_\lambda \bigl(X_j^T
y \bigr),\qquad j=1,\ldots, p,
\]
where $X_1, \ldots, X_p$ are columns of $X$, and $S_\lambda\dvtx
\mathbb{R}\rightarrow\mathbb{R}$ is the soft-thresholding function,
\[
S_\lambda(x) = \cases{ x-\lambda, &\quad if $x > \lambda$,
\cr
0, &\quad
if $-\lambda\leq x \leq\lambda$,
\cr
x+\lambda, &\quad if $x < \lambda$.}
\]
Letting $U_j=X_j^Ty$, $j=1,\ldots, p$, the knots in the lasso
path are simply the values of $\lambda$ at which the coefficients become
nonzero (i.e., cease to be thresholded),
\[
\lambda_1=|U_{(1)}|,\qquad \lambda_2=|U_{(2)}|,\qquad \ldots,\qquad \lambda_p=|U_{(p)}|,
\]
where $|U_{(1)}| \geq|U_{(2)}| \geq\cdots\geq|U_{(p)}|$
are the order statistics of $|U_1|, \ldots, |U_p|$ (somewhat of an abuse
of notation). Therefore,
\[
T_k = |U_{(k)}|\bigl(|U_{(k)}|-|U_{(k+1)}|\bigr)/
\sigma^2.
\]
Next, we study the special case $k=1$, the test for the first predictor
to enter the active set along the lasso path.
We then examine the case $k\geq1$, the test at a general step
in the lasso path.

\subsection{\texorpdfstring{The first step, $k=1$.}{The first step, $k=1$}}\label{secorthxfirst}

Consider the covariance test statistic for the first
predictor to enter the active set, that is, for $k=1$,
\[
T_1 = |U_{(1)}|\bigl(|U_{(1)}|-|U_{(2)}|\bigr)/
\sigma^2.
\]
We are interested in the distribution of $T_1$ under the null hypothesis;
since we are testing the first predictor to enter, this is
\[
H_0\dvtx y \sim N \bigl(0,\sigma^2 I \bigr).
\]
Under the null, $U_1,\ldots, U_p$ are i.i.d., $U_j \sim N(0,\sigma^2)$,
and so $|U_1|/\sigma, \ldots, |U_p|/\sigma$ follow a $\chi_1$ distribution
(absolute value of a standard Gaussian). That $T_1$ has an asymptotic
$\Exp(1)$ null distribution is now given by the next result.

%
%
\begin{lemma}
\label{lemtoptwo}
Let $V_1 \geq V_2 \geq\cdots\geq V_p$ be the order statistics
of an independent sample of $\chi_1$ variates (i.e.,
they are the sorted absolute values of an independent sample of
standard Gaussian variates). Then
\[
V_1(V_1-V_2) \stackrel{d} {\rightarrow}
\Exp(1)\qquad\mbox{as } p\rightarrow\infty.
\]
\end{lemma}

This lemma reveals a remarkably simple limiting distribution
for the largest of independent $\chi_1$ random variables times
the gap between the largest two; we skip its proof, as it is a special
case of the following generalization.

%
%
\begin{lemma}
\label{lemanytwo}
If $V_1 \geq V_2 \geq\cdots\geq V_p$ are the order statistics
of an independent sample of $\chi_1$ variates, then for any fixed $k
\geq1$,
\begin{eqnarray*}
&& \bigl(V_1(V_1-V_2),V_2(V_2-V_3),
\ldots, V_k(V_k-V_{k+1}) \bigr)
\\
&&\qquad \stackrel{d} {\rightarrow} \bigl(\Exp(1),\Exp(1/2), \ldots,\Exp(1/k) \bigr) \qquad \mbox{as } p
\rightarrow\infty,
\end{eqnarray*}
where the limiting distribution (on the right-hand side above) has independent
components. To be perfectly clear, here and throughout we use
$\Exp(\alpha)$ to denote the exponential distribution with scale
parameter $\alpha$ (not rate parameter $\alpha$), so that if $Z\sim
\Exp(\alpha)$, then $\mathbb{E}[Z]=\alpha$.
\end{lemma}

\begin{pf}
The $\chi_1$ distribution has CDF
\[
F(x) = \bigl(2\Phi(x)-1 \bigr)1\{x \geq0\},
\]
where $\Phi$ is the standard normal CDF. We first compute
\[
\lim_{t\rightarrow\infty} \frac{F''(t)(1-F(t))}{(F'(t))^2} = \lim_{t\rightarrow\infty} -
\frac{t(1-\Phi(t))}{\phi(t)} = -1,
\]
the last equality using Mills' ratio.
Theorem 2.2.1 in
\citet{evt} then implies that, for constants $a_p=F^{-1}(1-1/p)$
and $b_p=pF'(a_p)$,
\[
b_p(V_1-a_p) \stackrel{d} {\rightarrow}-
\log{E_0},
\]
where $E_0$ is a standard exponential variate, so $-\log{E_0}$
has the standard (or type~I) extreme value distribution. Hence,
according to Theorem 3 in \citet{weissman}, for any fixed $k
\geq1$, the random variables $W_0=b_p(V_{k+1}-a_p)$ and
$W_i=b_p(V_i-V_{i+1})$, $i=1,\ldots,k$, converge jointly:
\[
(W_0,W_1,W_2,\ldots,W_k)
\stackrel{d} {\rightarrow}(-\log{G_0}, E_1/1,E_2/2,
\ldots,E_k/k),
\]
where $G_0,E_1,\ldots,E_k$ are independent, $G_0$ is Gamma
distributed with scale parameter 1 and shape parameter $k$, and
$E_1,\ldots,E_k$ are standard exponentials. Now note that
\begin{eqnarray*}
V_i(V_i-V_{i+1}) & =& \Biggl(a_p+
\frac{W_{0}}{b_p}+\sum_{j=i}^k
\frac{W_j}{b_p} \Biggr)\frac{W_i}{b_p}
\\
& =& \frac{a_p}{b_p}W_i + \frac{1}{b_p^2} \Biggl(W_0+
\sum_{j=i}^kW_j
\Biggr)W_i.
\end{eqnarray*}
We claim that $a_p/b_p \rightarrow1$; this would give the desired
result as the second term converges to zero, using
$b_p\rightarrow\infty$. Writing $a_p,b_p$ more explicitly, we see
that $1-1/p = 2\Phi(a_p)-1$, that is, $1-\Phi(a_p) = 1/(2p)$, and
$b_p = 2p\phi(a_p)$.
Using Mills' inequalities,
\[
\frac{\phi(a_p)}{a_p} \frac{1}{1+1/a_p^2} \leq1-\Phi(a_p) \leq
\frac{\phi(a_p)}{a_p}
\]
and multiplying by $2p$,
\[
\frac{b_p}{a_p} \frac{1}{1+1/a_p^2} \leq1 \leq\frac{b_p}{a_p}.
\]
Since $a_p\rightarrow\infty$, this means that $b_p/a_p \rightarrow1$,
completing the proof.
\end{pf}


Practically, Lemma~\ref{lemanytwo} tells us that under the global
null hypothesis
$y\sim N(0,\sigma^2)$, comparing the covariance statistic $T_k$ at the
$k$th step
of the lasso path to an $\Exp(1)$ distribution is increasingly
conservative [at the
first step, $T_1$ is asymptotically $\Exp(1)$, at the second step,
$T_2$ is
asymptotically $\Exp(1/2)$, at the third step, $T_3$ is asymptotically
$\Exp(1/3)$,
and so forth]. This progressive conservatism is favorable, if we place
importance on
parsimony in the fitted model: we are less and less likely to incur a
false rejection of
the null hypothesis as the size of the model grows. Moreover, we know that
the test statistics $T_1,T_2,\ldots$ at successive steps are
independent, and hence so
are the corresponding $p$-values; from the point of view of multiple
testing corrections,
this is nearly an ideal scenario.

Of real interest is the distribution of $T_k$, $k\geq1$, not under
the global null hypothesis, but rather, under the weaker null
hypothesis that all
variables excluded from the current lasso model are truly inactive
(i.e., they have zero coefficients in the true model).
We study this in next section.

\subsection{\texorpdfstring{A general step, $k\geq1$.}{A general step, $k>=1$}}\label{secorthxgen}

We suppose that exactly $k_0$ components of the true coefficient vector
$\beta^*$ are nonzero, and consider testing the entry of the predictor at
step $k=k_0+1$. Let $A^*=\supp(\beta^*)$ denote the true active set
(so $k_0=|A^*|$), and let $B$ denote the event that all truly active
variables are added at steps $1,\ldots, k_0$,
%
%
\begin{equation}
\label{eqorthxb} B = \Bigl\{ \min_{j\in A^*} |U_j| >
\max_{j\notin A^*} |U_j| \Bigr\}.
\end{equation}
We show that under the null hypothesis (i.e., conditional on $B$),
the test statistic $T_{k_0+1}$ is asymptotically
$\Exp(1)$, and further, the test statistic $T_{k_0+d}$ at a future step
$k=k_0+d$ is asymptotically $\Exp(1/d)$.

The basic idea behind our argument is as follows: if we assume that the nonzero
components of $\beta^*$ are large enough in magnitude, then it is not
hard to
show (relying on orthogonality, here) that the truly active predictors
are added to the model along the first $k_0$ steps of the lasso path,
with probability tending to one. The test statistic at the $(k_0+1)$st
step and beyond would therefore
depend on the order statistics of $|U_i|$ for truly inactive variables
$i$, subject to
the constraint that the largest of these values is smaller than the smallest
$|U_j|$ for truly active variables $j$. But with our strong signal assumption,
that is, that the nonzero entries of $\beta^*$ are large in absolute
value, this
constraint has essentially no effect, and we are back to studying the
order statistics
from a $\chi_1$ distribution, as in the last section. This is made
precise below.

%
%
\begin{theorem}
\label{thmorthxgen}
Assume that $X\in\mathbb{R}^{n\times p}$ is orthogonal, and $y\in
\mathbb{R}^n$ is
drawn from the normal regression model \eqref{eqmodel}, where the
true coefficient vector $\beta^*$ has $k_0$ nonzero components.
Let $A^*=\supp(\beta^*)$ be the true active set, and assume that the
smallest nonzero true coefficient is large compared to
$\sigma\sqrt{2\log{p}}$,
\[
\min_{j\in A^*} \bigl|\beta^*_j\bigr| - \sigma\sqrt{2\log{p}}
\rightarrow\infty\qquad\mbox{as } p\rightarrow\infty.
\]
Let $B$ denote the event in \eqref{eqorthxb}, namely, that the first
$k_0$ variables entering the model along the lasso path are those in
$A^*$. Then $\mathbb{P}(B) \rightarrow1$ as
$p\rightarrow\infty$, and for each fixed $d \geq0$, we have
\[
(T_{k_0+1},T_{k_0+2},\ldots, T_{k_0+d}) \stackrel{d} {
\rightarrow} \bigl(\Exp(1),\Exp(1/2),\ldots,\Exp(1/d) \bigr) \qquad\mbox {as } p
\rightarrow\infty.
\]
The same convergence in distribution holds conditionally on $B$.
\end{theorem}

\begin{pf}
We first study $\mathbb{P}(B)$. Let $\theta_p=\min_{i\in A^*} |\beta
^*_i|$, and
choose $c_p$ such that
\[
c_p-\sigma\sqrt{2\log{p}} \rightarrow\infty\quad\mbox{and}\quad
\theta_p - c_p \rightarrow\infty.
\]
Note that $U_j \sim N(\beta^*_j,\sigma^2)$, independently
for $j=1,\ldots, p$. For $j\in A^*$,
\[
\mathbb{P}\bigl(|U_j| \leq c_p\bigr) = \Phi \biggl(
\frac{c_p-\beta^*_i}{\sigma} \biggr) - \Phi \biggl(\frac{-c_p-\beta
^*_i}{\sigma} \biggr) \leq\Phi
\biggl(\frac{c_p-\theta_p}{\sigma} \biggr) \rightarrow0,
\]
so
\[
\mathbb{P} \Bigl(\min_{j\in A^*} |U_j| >
c_p \Bigr) = \prod_{j\in A^*}
\mathbb{P}\bigl(|U_j| > c_p\bigr) \rightarrow1.
\]
At the same time,
\[
\mathbb{P} \Bigl(\max_{j\notin A^*} |U_j| \leq
c_p \Bigr) = \bigl(\Phi(c_p/\sigma)-\Phi(-c_p/
\sigma) \bigr)^{p-k_0} \rightarrow1.
\]
Therefore, $\mathbb{P}(B) \rightarrow1$.
This in fact means that $\mathbb{P}(E|B)-\mathbb{P}(E) \rightarrow0$ for
any sequence of events $E$, so only the weak convergence
of $(T_{k_0+1},\ldots, T_{k_0+d})$ remains to be proved. For this, we let
$m=p-k_0$, and $V_1 \geq V_2 \geq\cdots\geq V_m$ denote the
order statistics of the sample $|U_j|$, $j\notin A^*$ of independent
$\chi_1$ variates. Then, on the event~$B$, we have
\[
T_{k_0+i} = V_i(V_i-V_{i+1}) \qquad
\mbox{for } i=1,\ldots, d.
\]
As $\mathbb{P}(B) \rightarrow1$, we have in general
\[
T_{k_0+i} = V_i(V_i-V_{i+1}) +
o_\mathbb{P}(1) \qquad\mbox{for } i=1,\ldots, d.
\]
Hence, we are essentially back in the setting of the last section,
and the desired convergence result follows from the same arguments
as those for Lemma~\ref{lemanytwo}.
\end{pf}

%

\section{\texorpdfstring{A general predictor matrix $X$.}{A general predictor matrix $X$}}\label{secgenx}

In this section, we consider a general predictor matrix $X$, with
columns in general position.
Recall that our proposed covariance test statistic \eqref{eqcovtest}
is closely intertwined with the knots $\lambda_1 \geq\cdots\geq
\lambda_r$ in the lasso path, as it was defined in terms of difference
between fitted values at successive knots. Moreover, Lemma
\ref{lemcovtestknot}
showed that (provided there are no sign changes in the reduced lasso
problem over $[\lambda_{k+1},\lambda_k]$) this test statistic can be
expressed even more explicitly in terms of the values of these knots.
As was the case in the last section, this knot form is quite important
for our analysis here.
Therefore, it is helpful to recall [\citet{lars,lassounique}] the
precise formulae for the knots in the lasso path.
If $A$ denotes the active set and $s_A$ denotes the
signs of active coefficients at a knot $\lambda_k$,
\[
A = \supp \bigl(\hbeta(\lambda) \bigr),\qquad s_A = \sign \bigl(
\hbeta_A(\lambda_k) \bigr),
\]
then the next knot $\lambda_{k+1}$ is given by
%
%
\begin{equation}
\label{eqnextknot} \lambda_{k+1} = \max \bigl\{\lambda_{k+1}^\mathrm{join},
\lambda_{k+1}^\mathrm{leave} \bigr\},
\end{equation}
where $\lambda_{k+1}^\mathrm{join}$ and
$\lambda_{k+1}^\mathrm{leave}$ are the values of $\lambda$
at which, if we were to decrease the tuning parameter from $\lambda_k$
and continue along the current (linear) trajectory for the lasso
coefficients, a variable would join and leave the active set $A$,
respectively. These values are\footnote{In expressing the joining
and leaving times in the forms \eqref{eqnextjoin} and
\eqref{eqnextleave}, we are implicitly assuming that $\lambda_{k+1}<
\lambda_k$, with strict inequality. Since $X$ has columns in general
position, this is true for (Lebesgue) almost every $y$, or in other
words, with probability one taken over the normally distributed errors
in \eqref{eqmodel}.}
%
%
\begin{equation}
\label{eqnextjoin}\quad  \lambda_{k+1}^\mathrm{join} = \max
_{j\notin A, s\in\{-1,1\}} \frac{X_j^T (I-P_A) y} {
s - X_j^T (X_A^T)^+ s_A} \cdot1 \biggl\{\frac{X_j^T (I-P_A) y} {
s - X_j^T (X_A^T)^+ s_A} <
\lambda_k \biggr\},
\end{equation}
where recall $P_A=X_A(X_A^T X_A)^{-1} X_A^T$, and
$(X_A^T)^+ = X_A (X_A^T X_A)^{-1}$; and
%
%
\begin{equation}
\label{eqnextleave} \lambda_{k+1}^\mathrm{leave} = \max
_{j \in A} \frac{ [(X_A)^+y]_j} {
[(X_A^T X_A)^{-1} s_A]_j} \cdot1 \biggl\{ \frac{ [(X_A)^+y]_j} {
[(X_A^T X_A)^{-1} s_A]_j } <
\lambda_k \biggr\}.
\end{equation}
%

As we did in Section~\ref{secorthx} with the orthogonal $X$ case, we
begin by studying the asymptotic distribution of the covariance
statistic in the special case $k=1$ (i.e., the first model along the path),
wherein the expressions for the next knot \eqref{eqnextknot},
\eqref{eqnextjoin}, \eqref{eqnextleave} greatly simplify. Following
this, we study the more difficult case $k\geq1$. For the sake of
readability, we defer the proofs and most technical details until the
\hyperref[app11]{Appendix}.

\subsection{\texorpdfstring{The first step, $k=1$.}{The first step, $k=1$}}\label{secgenxfirst}

We assume here that $X$ has unit normed columns:
$\|X_i\|_2=1$, for $i=1,\ldots, p$; we do this mostly for simplicity of
presentation, and the generalization to
a matrix $X$ whose columns are not unit normed is given in the next
section (though the exponential limit is now a conservative upper
bound).
As per our discussion following Lemma~\ref{lemcovtestknot} (see the
second point), we know that
the first predictor to enter the active set along the lasso path
cannot leave at the next step, so the constant sign condition
\eqref{eqnosignchange} holds, and by Lemma~\ref{lemcovtestknot}
the covariance statistic for testing the entry of the first
variable can be written as
\[
T_1 = \lambda_1(\lambda_1-
\lambda_2)/\sigma^2
\]
(the leading factor $C$ being
equal to one since we assumed that $X$ has unit normed columns).
Now let $U_j=X_j^T y$, $j=1,\ldots, p$, and $R=X^T X$. With $\lambda
_0=\infty$,
we have $A=\varnothing$, and trivially, no variables can leave the active
set. The first knot is hence given by \eqref{eqnextjoin}, which
can be expressed as
%
%
\begin{equation}
\label{eqlambda1} \lambda_1 = \max_{j=1,\ldots, p, s\in\{-1,1\}} s
U_j.
\end{equation}
Letting $j_1,s_1$ be the first variable to enter and its sign (i.e.,
they achieve
the maximum in the above expression), and recalling that $j_1$ cannot
leave the active set immediately after it has entered, the second knot
is again given by \eqref{eqnextjoin}, written as
\[
\lambda_2 = \max_{j\neq j_1, s\in\{-1,1\}} \frac{s U_j - s R_{j,j_1}
U_{j_1}} {
1 - s s_1 R_{j,j_1}} \cdot1
\biggl\{ \frac{s U_j - s R_{j,j_1} U_{j_1}} {
1 - s s_1 R_{j,j_1}} < s_1U_{j_1} \biggr\}.
\]
The general position assumption on $X$ implies that $|R_{j,j_1}|<1$,
and so
$1-s s_1 R_{j,j_1}>0$, all $j\neq j_1$, $s\in\{-1,1\}$. It is easy to
show then
that the indicator inside the maximum above can be dropped, and hence
%
%
\begin{equation}
\label{eqlambda2} \lambda_2 = \max_{j\neq j_1, s\in\{-1,1\}}
\frac{s U_j - s R_{j,j_1} U_{j_1}} {
1 - s s_1 R_{j,j_1}}.
\end{equation}
Our goal now is to calculate the asymptotic distribution of
$T_1=\lambda_1(\lambda_1-\lambda_2)/\sigma^2$, with $\lambda_1$
and $\lambda_2$ as
above, under the null hypothesis; to be clear, since we are testing the
significance of
the first variable to enter along the lasso path, the null hypothesis is
%
%
\begin{equation}
\label{eqnull} H_0\dvtx y \sim N \bigl(0,\sigma^2 I
\bigr).
\end{equation}
The strategy that we use here for the general $X$ case---which differs from
our extreme value theory approach for the orthogonal $X$ case---is to treat
the quantities inside the maxima in expressions
\eqref{eqlambda1}, \eqref{eqlambda2} for $\lambda_1,\lambda_2$ as
discrete-time
Gaussian processes. First, we consider the zero mean Gaussian process
%
%
\begin{equation}
\label{eqg} g(j,s) = s U_j\qquad\mbox{for } j=1,\ldots, p, s \in
\{-1,1\}.
\end{equation}
We can easily compute the covariance function of this process:
\[
\mathbb{E} \bigl[g(j,s)g \bigl(j',s' \bigr) \bigr] =
ss' R_{j,j'} \sigma^2,
\]
where the expectation is taken over the null distribution in \eqref{eqnull}.
From \eqref{eqlambda1}, we know that the first knot is simply
\[
\lambda_1 = \max_{j, s} g(j,s).
\]
%
In addition to \eqref{eqg}, we consider the process
%
%
\begin{equation}
\label{eqh} \qquad h^{(j_1,s_1)}(j,s) = \frac{g(j,s) - ss_1 R_{j,j_1}
g(j_1,s_1)}{1-ss_1 R_{j,j_1}} \qquad\mbox{for } j
\neq j_1, s\in\{-1,1\}. 
\end{equation}
%
An important property: for fixed $j_1,s_1$, the entire process
$h^{(j_1,s_1)}(j,s)$ is independent of $g(j_1,s_1)$. This
can be seen by verifying that
\[
\mathbb{E} \bigl[ g(j_1,s_1) h^{(j_1,s_1)}(j,s)
\bigr] = 0
\]
and noting that $g(j_1,s_1)$ and $h^{(j_1,s_1)}(j,s)$, all $j \neq j_1$,
$s\in\{-1,1\}$, are jointly normal.
Now define
%
%
\begin{equation}
\label{eqm} M(j_1,s_1) = \max_{j\neq j_1, s}
h^{(j_1,s_1)}(j,s)
\end{equation}
and from the above we know that for fixed $j_1,s_1$, $M(j_1,s_1)$ is independent
of $g(j_1,s_1)$. If $j_1,s_1$ are instead treated as random variables
that maximize
$g(j,s)$ (the argument maximizers being almost surely unique), then from
\eqref{eqlambda2} we see that the second knot is $\lambda_2=M(j_1,s_1)$.
Therefore, to study the distribution of
$T_1=\lambda_1(\lambda_1-\lambda_2)/\sigma^2$,
we are interested in the random variable
\[
g(j_1,s_1) \bigl( g(j_1,s_1)
- M(j_1,s_1) \bigr) / \sigma^2
\]
on the event
\[
\bigl\{ g(j_1,s_1) > g(j,s)\mbox{ for all } (j,s)
\neq(j_1,s_1) \bigr\}.
\]
It turns out that this event, which concerns the argument maximizers
of $g$, can be rewritten as an event concerning only the relative
values of
$g$ and $M$ [see \citet{taylorvalid} for the analogous result for
continuous-time processes].
%
%
\begin{lemma}
\label{lemeventdual}
With $g,M$ as defined in \eqref{eqg}, \eqref{eqh}, \eqref{eqm},
we have
\[
\bigl\{ g(j_1,s_1) > g(j,s) \mbox{ for all }
(j,s)\neq(j_1,s_1) \bigr\} = \bigl\{
g(j_1,s_1) > M(j_1,s_1) \bigr
\}.
\]
\end{lemma}
This is an important realization because the dual representation
$\{g(j_1,s_1) > M(j_1,s_1)\}$ is more tractable, once
we partition the space over the possible argument minimizers $j_1,s_1$, and
use the fact that $M(j_1,s_1)$ is independent of $g(j_1,s_1)$ for fixed
$j_1,s_1$. In this vein, we express the distribution of
$T_1=\lambda_1(\lambda_1-\lambda_2)/\sigma^2$ in terms of the sum
\begin{eqnarray*}
&& \mathbb{P}(T_1>t)
\\
&&\qquad = \sum_{j_1,s_1}
\mathbb{P} \bigl(g(j_1,s_1) \bigl( g(j_1,s_1)
- M(j_1,s_1) \bigr) / \sigma^2 > t,
g(j_1,s_1) > M(j_1,s_1)
\bigr).
\end{eqnarray*}
The terms in the above sum can be simplified: dropping for
notational convenience the
dependence on $j_1,s_1$, we have
\[
g(g-M)/\sigma^2 > t,\qquad g > M\quad\Longleftrightarrow\quad g/\sigma> u(t,M/\sigma),
\]
where $u(a,b)= (b+\sqrt{b^2+4a})/2$, which follows by simply solving
for $g$ in the quadratic equation $g(g-M)/\sigma^2 = t$. Therefore,
%
%
\begin{eqnarray}\label{eqt1int}
\mathbb{P}(T_1>t) &=& \sum_{j_1,s_1}
\mathbb{P} \bigl(g(j_1,s_1)/\sigma> u
\bigl(t,M(j_1,s_1)/\sigma \bigr) \bigr)
\nonumber\\[-8pt]\\[-8pt]
&=& \sum_{j_1,s_1} \int
_0^\infty\widebar{\Phi} \bigl(u(t,m/\sigma) \bigr)
F_{M(j_1,s_1)}(dm),\nonumber
\end{eqnarray}
where $\widebar{\Phi}$ is the standard normal survival function
(i.e., $\widebar{\Phi}=1-\Phi$, for $\Phi$ the standard normal CDF),
$F_{M(j_1,s_1)}$ is the distribution
of $M(j_1,s_1)$, and we have used the fact that $g(j_1,s_1)$ and $M(j_1,s_1)$
are independent for fixed $j_1,s_1$, as well as $M(j_1,s_1) \geq0$.
Continuing from \eqref{eqt1int}, we can write the difference between
\mbox{$\mathbb{P}(T_1>t)$}
and the standard exponential tail, $\mathbb{P}(\Exp(1)>t)=e^{-t}$, as
%
%
\begin{eqnarray}\label{eqt1bd}
&& \bigl|\mathbb{P}(T_1>t) - e^{-t} \bigr|
\nonumber\\[-8pt]\\[-8pt]
&&\qquad  = \Biggl| \sum
_{j_1,s_1} \int_0^\infty
\biggl(\frac{\widebar{\Phi} (u(t,m/\sigma) )} {
\widebar{\Phi}(m/\sigma)} - e^{-t} \biggr) \widebar{\Phi}(m/\sigma)
F_{M(j_1,s_1)}(dm) \Biggr|,\nonumber
\end{eqnarray}
where we used the fact that
\[
\sum_{j_1,s_1} \int_0^\infty
\widebar{\Phi} (m/\sigma) F_{M(j_1,s_1)}(dm) = \sum
_{j_1,s_1} \mathbb{P} \bigl(g(j_1,s_1) >
M(j_1,s_1) \bigr) = 1.
\]
We now examine the term inside the braces in
\eqref{eqt1bd}, the difference between a ratio of normal survival
functions and $e^{-t}$; our next lemma shows that this term vanishes as
$m\rightarrow\infty$.

%
%
\begin{lemma}
\label{lemexpratio}
For any $t\geq0$,
\[
\frac{\widebar{\Phi} (u(t,m) )}{\widebar{\Phi}(m)} \rightarrow e^{-t} \qquad \mbox{as } m\rightarrow
\infty.
\]
\end{lemma}

Hence, loosely speaking, if each $M(j_1,s_1) \rightarrow\infty$ fast
enough as
$p\rightarrow\infty$, then the right-hand side in \eqref{eqt1bd}
converges to
zero, and $T_1$ converges weakly to $\Exp(1)$. This is made precise below.

%
%
\begin{lemma}
\label{lemexpconv}
Consider $M(j_1,s_1)$ defined in \eqref{eqh}, \eqref{eqm} over
$j_1=1,\ldots, p$ and $s_1\in\{-1,1\}$. If for any fixed $m_0>0$
%
%
\begin{equation}
\label{eqmtoinf} \sum_{j_1,s_1} \mathbb{P}
\bigl(M(j_1,s_1) \leq m_0 \bigr) \rightarrow
0 \qquad\mbox{as } p\rightarrow\infty,
\end{equation}
then the right-hand side in \eqref{eqt1bd} converges to zero as
$p\rightarrow\infty$,
and so $\mathbb{P}(T_1 > t) \rightarrow e^{-t}$ for all $t \geq0$.
\end{lemma}

The assumption in \eqref{eqmtoinf} is written in terms of random
variables whose distributions are induced by the steps along the lasso
path; to make
our assumptions more transparent, we show that \eqref{eqmtoinf} is
implied by a
conditional variance bound involving the predictor matrix $X$ alone,
and arrive at
the main result of this section.

%
%
\begin{theorem}
\label{thmgenxfirst}
Assume that $X\in\mathbb{R}^{n\times p}$ has unit normed columns in general
position, and let $R=X^T X$.
Assume also that there is some $\delta>0$ such that for each
$j=1,\ldots, p$,
there exists a subset of indices $S \subseteq\{1,\ldots, p\}\setminus\{
j\}$ with
%
%
\begin{equation}
\label{eqcv} 1 - R_{i,S\setminus\{i\}} (R_{S\setminus\{i\},S\setminus\{
i\}})^{-1}
R_{S\setminus\{i\},i} \geq\delta^2 \qquad\mbox{for all } i\in S,
\end{equation}
and the size of $S$ growing faster than $\log{p}$,
%
%
\begin{equation}
\label{eqsgrow} |S| \geq d_p\qquad\mbox{where } \frac{d_p}{\log{p}}
\rightarrow\infty\mbox{ as } p\rightarrow\infty.
\end{equation}
The under the null distribution in \eqref{eqnull}
[i.e., $y$ is drawn from the regression model~\eqref{eqmodel}
with $\beta^*=0$], we have
$\mathbb{P}(T_1 > t) \rightarrow e^{-t}$ as $p\rightarrow\infty$ for all
$t\geq0$.
\end{theorem}

\begin{rem*}
Conditions \eqref{eqcv} and \eqref{eqsgrow} are
sufficient to ensure \eqref{eqmtoinf}, or in other words, that each
$M(j_1,s_1)$ grows as in
$\mathbb{P}(M(j_1,s_1) \leq m_0) = o(1/p)$, for any fixed~$m_0$. While
it is
true that $\mathbb{E}[M(j_1,s_1)]$ will typically grow as $p$ grows, some
assumption is required so that $M(j_1,s_1)$ concentrates
around its mean faster than standard Gaussian concentration
results (such as the Borell-TIS inequality) imply.

Generally speaking, the assumptions \eqref{eqcv} and \eqref{eqsgrow}
are not very strong. Stated differently, \eqref{eqcv} is a lower
bound on the variance of $U_i=X_i^T y$, conditional on
$U_\ell=X_\ell^T y$ for all $\ell\in S \setminus\{i\}$.
Hence, for any $j$, we require
the existence of a subset $S$ not containing $j$ such that
the variables $U_i$, $i\in S$, are not too correlated, in
the sense that the conditional variance of any one given all the others
is bounded below. This subset $S$ has to be larger in size than
$\log{p}$, as made clear in \eqref{eqsgrow}. Note that, in fact, it
suffices to find a total of two disjoint subsets $S_1, S_2$ with the
properties \eqref{eqcv} and \eqref{eqsgrow}, because then for any
$j$, either one or the other will not contain $j$.

An example of a matrix $X$ that does not satisfy
\eqref{eqcv} and \eqref{eqsgrow} is one with fixed rank
as $p$ grows. (This, of course, would also not
satisfy the general position assumption.) In this case, we would not
be able to find a subset of the variables $U_i=X_i^T y$, $i=1,\ldots,
p$, that is both linearly independent and has size larger than
$r=\rank(X)$, which violates the conditions. We note that in general,
since $|S| \leq\rank(X) \leq n$, and $|S|/\log{p} \rightarrow
\infty$, conditions \eqref{eqcv} and \eqref{eqsgrow} require
that $n/\log{p} \rightarrow\infty$.
\end{rem*}

\subsection{\texorpdfstring{A general step, $k\geq1$.}{A general step, $k>=1$}}\label{secgenxgen}

In this section, we no longer assume that $X$ has unit normed columns
(in any case, this provides no simplification in deriving the
null distribution of the test statistic at a general step in the lasso
path). Our arguments here have more or less the same form as they did
in the last section, but overall the calculations are more complicated.

Fix an integer $k_0 \geq0$, subset $A_0 \subseteq
\{1,\ldots, p\}$ containing the true active set
$A_0 \supseteq A^* = \supp(\beta^*)$, and sign vector $s_{A_0} \in
\{-1,1\}^{|A_0|}$. Consider the event
%
%
\begin{eqnarray}
\label{eqgenxb} B &=& \biggl\{\mbox{the solution at step $k_0$ in the
lasso path has active set $A=A_0$,}\nonumber\hspace*{-18pt}
\\
&&\hspace*{6pt} \mbox{signs $s_A=\sign \bigl((X_{A_0})^+y
\bigr)=s_{A_0}$, and the next two knots are given by}\hspace*{-18pt}
\\
&&\hspace*{53pt} \lambda_{k_0+1} = \max_{j\notin A\cup\{j_{k_0}\}, s\in\{-1,1\}} \frac
{X_j^T (I-P_A) y} {
s - X_j^T (X_A^T)^+ s_A},
\lambda_{k_0+2}=\lambda_{k_0+2}^\mathrm{join} \biggr\}.\nonumber\hspace*{-18pt}
\end{eqnarray}
We assume that $\mathbb{P}(B) \rightarrow1$ as $p \rightarrow\infty$.
In words, this is assuming that with probability approaching
one:
the lasso estimate at step $k_0$ in the path has support $A_0$ and
signs $s_{A_0}$; the least squares estimate on
$A_0$ has the same signs as this lasso estimate; the knots at steps
$k_0+1$ and $k_0+2$ correspond to joining events; and in particular,
the maximization defining the joining event
at step $k_0+1$ can be taken to be unrestricted, that is, without the
indicators constraining the individual arguments to be $<
\lambda_{k_0}$.
Our goal is to characterize the asymptotic distribution of the
covariance statistic $T_k$ at the step $k=k_0+1$, under the null
hypothesis (i.e., conditional on the event~$B$). We will comment on
the stringency of the assumption that $\mathbb{P}(B) \rightarrow1$ following
our main result in Theorem~\ref{thmgenxgen}.

First note that on $B$, we have $s_A=\sign((X_A)^+y)$, and
as discussed in the third point
following Lemma~\ref{lemcovtestknot}, this implies
that the solution of the reduced
problem \eqref{eqtbeta} on $X_A$ cannot incur any sign changes over
the interval $[\lambda_k,\lambda_{k+1}]$. Hence, we
can apply Lemma~\ref{lemcovtestknot} to write the covariance
statistic on $B$ as
\[
T_k = C(A,s_A,j_k,s_k) \cdot
\lambda_k(\lambda_k-\lambda_{k+1}) /
\sigma^2,
\]
where\vspace*{2pt} $C(A,s_A,j_k,s_k)=\|(X_{A\cup\{ j_k\}}^T)^+ s_{A \cup\{j_k\}} -
(X_A^T)^+ s_A \|_2^2$, $A$ and $s_A$ are the active set and signs at
step $k-1$, and $j_k$ is the variable added to the active
set at step~$k$, with sign $s_k$. Now,
analogous to our definition in the last section, we define the
discrete-time Gaussian process
%
%
\begin{equation}
\label{eqg2} g^{(A,s_A)}(j,s) = \frac{X_j^T (I-P_A) y} {
s - X_j^T (X_A^T)^+ s_A} \qquad\mbox{for } j
\notin A, s \in\{-1,1\}.
\end{equation}
For any fixed $A,s_A$, the above process has mean
zero provided that $A \supseteq A^*$. Additionally, for any such fixed
$A,s_A$, we can compute its covariance function
%
%
\begin{equation}
\label{eqg2cov} \qquad\mathbb{E} \bigl[g^{(A,s_A)}(j,s) g^{(A,s_A)}
\bigl(j',s' \bigr) \bigr] = \frac{X_j^T (I-P_A) X_{j'} \sigma^2}{[s -
X_j^T (X_A^T)^+ s_A]
[s' - X_{j'}^T (X_A^T)^+ s_A]}.
\end{equation}
Note that on the event $B$, the $k$th knot in the lasso path is
\[
\lambda_k = \max_{j\notin A, s\in\{-1,1\}} g^{(A,s_A)}(j,s).
\]
For fixed $j_k,s_k$, we also consider the process
%
%
\begin{eqnarray}\label{eqh2}
g^{(A\cup\{j_k\},s_{A\cup\{j_k\}})} (j,s) = \frac{X_j^T
(I-P_{A\cup
\{j_k\}}) y} {
s - X_j^T (X_{A\cup\{j_k\}}^T)^+ s_{A\cup\{j_k\}}}
\nonumber\\[-8pt]\\[-8pt]
\eqntext{\displaystyle\mbox{for } j \notin A \cup\{j_k\}, s \in\{-1,1\}}
\end{eqnarray}
(above, $s_{A\cup\{j_k\}}$ is the concatenation of $s_A$ and $s_k$)
and its achieved maximum value, subject to being less than
the maximum of $g^{(A,s_A)}$,
%
%
\begin{eqnarray}\label{eqm2}
\qquad M^{(A,s_A)}(j_k,s_k)
&=&  \max_{j\notin A \cup\{j_k\} s \in\{-1,1\}} g^{(A\cup\{j_k\},s_{A\cup\{j_k\}
})}(j,s)
\nonumber\\[-8pt]\\[-8pt]
&&{}\times 1 \Bigl\{ g^{(A\cup\{j_k\},s_{A\cup\{j_k\}})} (j,s) < \max_{j\notin A,
s\in\{-1,1\}}
g^{(A,s_A)}(j,s) \Bigr\}.\nonumber
\end{eqnarray}
If $j_k,s_k$ indeed maximize $g^{(A,s_A)}$, that is, they correspond to
the variable added to the active set at $\lambda_k$ and its sign
(note that these are almost surely unique),
then on~$B$, we have $\lambda_{k+1} = M^{(A,s_A)}(j_k,s_k)$.
To study the distribution of $T_k$ on $B$,
we are therefore interested in the random variable
\[
C(A,s_A,j_k,s_k) \cdot g^{(A,s_A)}(j_k,s_k)
\bigl( g^{(A,s_A)}(j_k,s_k)-M^{(A,s_A)}(j_k,s_k)
\bigr)/\sigma^2
\]
on the event
%
%
\begin{equation}
\label{eqev2}\quad  E(j_k,s_k) = \bigl\{ g^{(A,s_A)}(j_k,s_k)
> g^{(A,s_A)}(j,s) \mbox{ for all } (j,s)\neq(j_k,s_k)
\bigr\}.
\end{equation}
Equivalently, we may write
\begin{eqnarray*}
&& \mathbb{P} \bigl(\{T_k > t\} \cap B \bigr)
\\
&&\qquad = \sum
_{j_k,s_k} \mathbb{P} \bigl( \bigl\{ C(A,s_A,j_k,s_k)
\cdot g^{(A,s_A)}(j_k,s_k)
\\
&&\hspace*{65pt}{}\times \bigl(g^{(A,s_A)}(j_k,s_k)-M^{(A,s_A)}(j_k,s_k)
\bigr)/\sigma^2 >t \bigr\} \cap E(j_k,s_k)
\bigr).
\end{eqnarray*}
Since $\mathbb{P}(B)\rightarrow1$, we have in general
%
%
\begin{eqnarray}\label{eqtksum}
&& \mathbb{P}(T_k > t)\nonumber
\\
&&\qquad = \sum
_{j_k,s_k} \mathbb{P} \bigl( \bigl\{C(A_0,s_{A_0},j_k,s_k)
\cdot g^{(A_0,s_{A_0})}(j_k,s_k) \hspace*{-18pt}
\nonumber\\[-8pt]\\[-8pt]
&&\hspace*{64pt}{}\times \bigl(g^{(A_0,s_{A_0})}(j_k,s_k)-M^{(A_0,s_{A_0})}(j_k,s_k)
\bigr)/\sigma^2 >t \bigr\} \cap E(j_k,s_k)
\bigr)\nonumber\hspace*{-18pt}
\\
&&\quad\qquad{} + o(1),\nonumber\hspace*{-18pt}
\end{eqnarray}
where we have replaced all instances of $A$ and $s_A$ on the
right-hand side above with the fixed subset $A_0$ and sign
vector $s_{A_0}$. This is a helpful
simplification, because in what follows we may now take $A=A_0$
and $s_A=s_{A_0}$ as fixed, and consider the distribution of the
random processes $g^{(A_0,s_{A_0})}$ and $M^{(A_0,s_{A_0})}$. With
$A=A_0$ and $s_A=s_{A_0}$ fixed, we drop the notational dependence
on them and write these processes as $g$ and $M$. We also write
the scaling factor $C(A_0,s_{A_0},j_k,s_k)$ as $C(j_k,s_k)$.

The setup in \eqref{eqtksum} looks very much like the one in the
last section [and to draw an even sharper parallel, the scaling factor
$C(j_k,s_k)$ is actually equal to one over the variance of
$g(j_k,s_k)$, meaning that $\sqrt{C(j_k,s_k)} \cdot g(j_k,s_k)$ is
standard normal for fixed
$j_k,s_k$, a fact that we will use later in the proof of Lemma
\ref{lemexpconv2}]. However, a major complication is that
$g(j_k,s_k)$ and $M(j_k,s_k)$ are no longer independent for fixed
$j_k,s_k$. Next, we derive a dual representation for
the event \eqref{eqev2} (analogous to Lemma~\ref{lemeventdual} in the last section), introducing a triplet of
random variables $M^+,M^-,M^0$---it turns out that $g$ is independent
of this triplet, for fixed~$j_k,s_k$.

%
%
\begin{lemma}
\label{lemeventdual2}
Let $g$ be as defined in \eqref{eqg2}
(with $A,s_A$ fixed at $A_0,s_{A_0}$). Let $\Sigma_{j,j'}$ denote
the covariance function of $g$ [short form for the expression
in \eqref{eqg2cov}].\footnote{To be perfectly clear, here
$\Sigma_{j,j'}$ actually depends on $s,s'$, but our notation
suppresses this dependence for brevity.}
Define
%
%
\begin{eqnarray}
S^+(j,s) &=& \biggl\{ \bigl(j',s'
\bigr) \dvtx j'\notin A \cup\{j\}, \frac{\Sigma_{j,j'}}{\Sigma_{jj}} < 1 \biggr\},
\nonumber\\[-8pt]\label{eqmp}  \\[-8pt]
M^+(j,s) &=& \max_{(j',s') \in S^+(j,s)} \frac{ g(j',s') - (\Sigma_{j,j'}/\Sigma
_{jj}) g(j,s)} {
1 - \Sigma_{j,j'}/\Sigma_{jj}},\nonumber
\\
S^-(j,s) &=& \biggl\{ \bigl(j',s'
\bigr) \dvtx j'\notin A \cup\{j\}, \frac{\Sigma_{j,j'}}{\Sigma_{jj}} > 1 \biggr\},
\nonumber\\[-8pt]\label{eqmm} \\[-8pt]
M^-(j,s) &=& \min_{(j',s') \in S^-(j,s)} \frac{ g(j',s') - (\Sigma_{j,j'}/\Sigma
_{jj}) g(j,s)} {
1 - \Sigma_{j,j'}/\Sigma_{jj}},\nonumber
\\
S^0(j,s) &=& \biggl\{ \bigl(j',s'
\bigr)\dvtx j'\notin A \cup\{j\}, \frac{\Sigma_{j,j'}}{\Sigma_{jj}} = 1 \biggr\},
\nonumber\\[-8pt]\label{eqm0} \\[-8pt]
M^0(j,s) &=& \max_{(j',s') \in S^0(j,s)} g \bigl(j',s'
\bigr) - (\Sigma_{j,j'}/\Sigma_{jj}) g(j,s).\nonumber
\end{eqnarray}
Then\vspace*{1pt} the event $E(j_k,s_k)$ in \eqref{eqev2}, that
$j_k,s_k$ maximize $g$, can be written as an intersection of events
involving $M^+,M^-,M^0$:
%
%
\begin{eqnarray}
&& \bigl\{ g(j_k,s_k) > g(j,s) \mbox{ for all } (j,s)
\neq(j_k,s_k) \bigr\}\nonumber
\\
&&\qquad  = \bigl\{ g(j_k,s_k)
> 0 \bigr\} \cap \bigl\{ g(j_k,s_k) >
M^+(j_k,s_k) \bigr\}
\\
&&\quad\qquad{} \cap \bigl\{ g(j_k,s_k) < M^-(j_k,s_k)
\bigr\} \cap \bigl\{0 > M^0(j_k,s_k) \bigr\}.\nonumber
\end{eqnarray}
\end{lemma}

As a result of Lemma~\ref{lemeventdual2}, continuing from
\eqref{eqtksum}, we can decompose the tail probability of $T_k$
as
%
%
\begin{eqnarray}\label{eqtksum2}
&& \mathbb{P}(T_k > t) \nonumber
\\
&&\qquad = \sum
_{j_k,s_k} \mathbb{P} \bigl(C(j_k,s_k)
\cdot g(j_k,s_k) \bigl(g(j_k,s_k)-M(j_k,s_k)
\bigr)/\sigma^2 >t, g(j_k,s_k) > 0,
\nonumber\\[-9pt]\\[-9pt]
&&\hspace*{71pt}  g(j_k,s_k) > M^+(j_k,s_k),
g(j_k,s_k) < M^-(j_k,s_k), 0 >
M^0(j_k,s_k) \bigr)\nonumber
\\
&&\quad\qquad{} + o(1).\nonumber
\end{eqnarray}
A key point here is that, for fixed $j_k,s_k$, the triplet
$M^+(j_k,s_k)$, $M^-(j_k,s_k)$, $M^0(j_k,s_k)$ is independent of
$g(j_k,s_k)$, which is true because
\[
\mathbb{E} \bigl[ g(j_k,s_k) \bigl(g(j,s) - (
\Sigma_{j_k,j}/\Sigma_{j_k,j_k}) g(j_k,s_k)
\bigr) \bigr] = 0
\]
and $g(j_k,s_k)$, along with
$g(j,s) - (\Sigma_{j_k,j}/\Sigma_{j_k,j_k}) g(j_k,s_k)$,
for all $j,s$, form a jointly Gaussian collection of random
variables. If we were to now replace $M$ by $M^+$ in the first line of
\eqref{eqtksum2}, and define a modified statistic $\widetilde{T}_k$ via
its tail probability,
%
%
\begin{eqnarray}\label{eqttksum}
&& \mathbb{P}(\widetilde{T}_k > t) \nonumber
\\
&&\qquad = \sum
_{j_k,s_k} \mathbb{P} \bigl(C(j_k,s_k)
\cdot g(j_k,s_k) \bigl(g(j_k,s_k)-M^+(j_k,s_k)
\bigr)/\sigma^2 >t,
\nonumber\\[-9pt]\\[-9pt]
&&\hspace*{64pt} g(j_k,s_k) > 0, g(j_k,s_k) > M^+(j_k,s_k),\nonumber
\\
&&\hspace*{122pt}
g(j_k,s_k) < M^-(j_k,s_k), 0 >
M^0(j_k,s_k) \bigr),\nonumber
\end{eqnarray}
then arguments similar to those in the second half of
Section~\ref{secgenxfirst} give a (conservative) exponential limit
for $\mathbb{P}(\widetilde{T}_k>t)$.

%
%
\begin{lemma}
\label{lemexpconv2}
Consider $g$ as defined in \eqref{eqg2}
(with $A,s_A$ fixed at $A_0,s_{A_0}$), and $M^+,M^-,M^0$ as defined in
\eqref{eqmp}, \eqref{eqmm}, \eqref{eqm0}.
Assume that for any fixed $m_0$,
%
%
\begin{equation}
\label{eqmptoinf} \sum_{j_k,s_k} \mathbb{P}
\bigl(M^+(j_k,s_k) \leq m_0/
\sqrt{C(j_k,s_k)} \bigr) \rightarrow0 \qquad{as } p
\rightarrow\infty.
\end{equation}
Then the modified statistic $\widetilde{T}_k$ in \eqref{eqttksum}
satisfies $\lim_{p\rightarrow\infty} \mathbb{P}(\widetilde{T}_k > t)
\leq e^{-t}$,
for all $t\geq0$.
\end{lemma}

Of course, deriving the limiting distribution of $\widetilde{T}_k$ was not
the goal, and it remains to relate $\mathbb{P}(\widetilde{T}_k>t)$ to
$\mathbb{P}(T_k>t)$. A fortuitous calculation shows that the two seemingly
different quantities $M^+$ and $M$---the former of which is
defined as the maximum of particular functionals of $g$,
and the latter concerned with the joining event at step
$k+1$---admit a very simple relationship: $M^+(j_k,s_k) \leq
M(j_k,s_k)$ for the maximizing $j_k,s_k$. We use this to bound the
tail of $T_k$.

%
%
\begin{lemma}
\label{lemexpconv3}
Consider $g,M$ as defined in \eqref{eqg2}, \eqref{eqh2},
\eqref{eqm2} (with $A,s_A$ fixed at $A_0,s_{A_0}$), and consider
$M^+$ as defined in \eqref{eqmm}. Then for any fixed $j_k,s_k$, on
the event $E(j_k,s_k)$ in \eqref{eqev2}, we have
\[
M^+(j_k,s_k) \leq M(j_k,s_k).
\]
Hence, if we assume as in Lemma~\ref{lemexpconv2}
the condition \eqref{eqmptoinf},
then $\lim_{p\rightarrow\infty} \mathbb{P}(T_k > t) \leq
e^{-t}$ for all $t \geq0$.
\end{lemma}

Though Lemma~\ref{lemexpconv3} establishes a (conservative)
exponential limit for the covariance statistic $T_k$, it
does so by enforcing assumption \eqref{eqmptoinf}, which is phrased
in terms of the tail distribution of a random process defined at the
$k$th step in the lasso path. We translate this into an explicit
condition on the covariance structure in \eqref{eqg2cov}, to make the
stated assumptions for exponential convergence more concrete.

%
%
\begin{theorem}
\label{thmgenxgen}
Assume that $X\in\mathbb{R}^{n\times p}$ has columns in general
position, and $y\in\mathbb{R}^n$ is drawn from the normal regression model
\eqref{eqmodel}. Assume that for a fixed integer $k_0 \geq0$,
subset $A_0 \subseteq\{1,\ldots, p\}$ with
$A_0\supseteq A^*=\supp(\beta^*)$, and sign vector $s_{A_0} \in
\{-1,1\}^{|A_0|}$, the event $B$ in \eqref{eqgenxb}
satisfies $\mathbb{P}(B) \rightarrow1$ as $p \rightarrow\infty$.
Assume that there exists a constant $0<\eta\leq1$ such that
%
%
\begin{equation}
\label{eqirrep} \bigl\|(X_{A_0})^+ X_j\bigr\|_1 \leq1-
\eta\qquad\mbox{for all } j \notin A_0.
\end{equation}
Define the matrix $R$ by
\[
R_{ij} = X_i^T (I-P_{A_0})
X_j\qquad\mbox{for } i,j \notin A_0.
\]
Assume that the diagonal elements in $R$ are all of the same order,
that is, $R_{ii}/R_{jj} \leq C$ for all $i,j$ and some constant $C>0$.
Finally
assume that, for each fixed $j \notin A_0$, there is a set
$S \subseteq\{1,\ldots, p\} \setminus(A_0 \cup\{j\})$ such that\vadjust{\goodbreak}
for~all~$i\in S$,
%
%
\begin{eqnarray}
\label{eqcv2}  \bigl[ R_{ii} - R_{i,S\setminus\{i\}}
(R_{S\setminus\{i\},S\setminus\{i\}})^{-1} R_{S\setminus\{i\},i} \bigr] / R_{ii}
&\geq&\delta^2,
\\
\label{eqcovvar}  |R_{ij}|/R_{jj} &<& \eta/ (2-\eta),
\\
\label{eqirrep2} \bigl\|(X_{A_0\cup\{j\}})^+ X_i \bigr\|_1 &<& 1,
\end{eqnarray}
where $\delta>0$ is a constant (not depending on $j$),
and the size of $S$ grows faster than $\log{p}$,
%
%
\begin{equation}
\label{eqsgrow2} |S| \geq d_p\qquad\mbox{where } \frac{d_p}{\log{p}}
\rightarrow\infty\mbox{ as } p\rightarrow\infty.
\end{equation}
Then at step $k=k_0+1$, we have $\lim_{p\rightarrow\infty} \mathbb
{P}(T_k >
t) \leq e^{-t}$ for all $t \geq0$. The same result holds for the tail
of $T_k$ conditional on $B$.
\end{theorem}

%
\begin{rem}\label{rem11} If $X$ has unit normed columns, then by taking
$k_0=0$ (and accordingly, $A_0=\varnothing$, $s_{A_0}=\varnothing$) in
Theorem~\ref{thmgenxgen}, we essentially recover the result of Theorem
\ref{thmgenxfirst}. To see this, note that with $k_0=0$ (and
$A_0,s_{A_0}=\varnothing$), we have $\mathbb{P}(B)=1$\vadjust{\goodbreak} for all finite
$p$ (recall
the arguments given at the beginning of Section~\ref{secgenxfirst}).
Also, condition \eqref{eqirrep} trivially holds
with $\eta=1$ because $A_0=\varnothing$. Next, the matrix $R$ defined in
the theorem reduces to $R=X^T X$, again because $A_0=\varnothing$; note
that $R$ has all diagonal elements equal to one, because $X$ has unit
normed columns. Hence, \eqref{eqcv2} is the same as condition
\eqref{eqcv} in Theorem~\ref{thmgenxfirst}. Finally, conditions
\eqref{eqcovvar} and \eqref{eqirrep2} both reduce to $|R_{ij}| < 1$,
which always holds as $X$ has columns in general position.
Therefore, when $k_0=0$, Theorem~\ref{thmgenxgen}
imposes the same conditions as Theorem~\ref{thmgenxfirst}, and
gives essentially the same result---we say
``essentially'' here is because the former gives a conservative
exponential limit for $T_1$, while the latter gives an exact
exponential limit.
\end{rem}

%
\begin{rem}\label{rem12} If $X$ is orthogonal, then for any $A_0$, conditions
\eqref{eqirrep} and \mbox{\eqref{eqcv2}--\eqref{eqsgrow2}} are trivially
satisfied [for the latter set of conditions, we can take, e.g.,
$S=\{1,\ldots, p\} \setminus(A_0 \cup\{j\})$].
With an additional condition on the strength of the true nonzero
coefficients, we can assure that $\mathbb{P}(B) \rightarrow1$ as
$p\rightarrow\infty$ with
$A_0=A^*$, $s_{A_0} = \sign(\beta_{A_0}^*)$, and $k_0 = |A_0|$, and
hence prove a conservative exponential limit for $T_k$; note that this
is precisely what is done in Theorem~\ref{thmorthxgen} (except that
in this case, the exponential limit is proven to be exact).
\end{rem}

%
\begin{rem}\label{rem13} Defining $U_i = X_i^T(I-P_{A_0})y$ for $i
\notin A_0$,
the condition \eqref{eqcv2} is a lower bound on the ratio of the
conditional variance of $U_i$ on $U_\ell$, $\ell\notin S$, to the
unconditional variance of $U_i$. Loosely speaking, conditions
\eqref{eqcv2}, \eqref{eqcovvar}, and \eqref{eqirrep2} can all
be
interpreted as requiring, for any $j \notin A_0$, the existence of a
subset $S$ not containing $j$ (and disjoint from $A_0$) such that
the variables $U_i$, $i \in S$, are not very correlated. This subset
has to be large in size compared to $\log{p}$, by \eqref{eqsgrow2}.
An implicit consequence of \eqref{eqcv2}--\eqref{eqsgrow2}, as
argued in the remark following Theorem~\ref{thmgenxfirst},
is that $n/\log{p} \rightarrow\infty$.\vadjust{\goodbreak}
\end{rem}

%
\begin{rem}\label{rem14} Some readers will likely recognize condition
\eqref{eqirrep} as that of \textit{mutual incoherence} or \textit{strong
irrepresentability}, commonly used in the lasso literature on exact
support recovery [see, e.g., \citet{sharp},
\citet{irrepcond}]. This condition, in addition to a lower bound
on the magnitudes of the true coefficients, is sufficient for the
lasso solution to recover the true active set $A^*$ with probability
tending to one, at a carefully chosen value of $\lambda$.
It is important to point out that we do not place any requirements on
the magnitudes of the true nonzero coefficients; instead,
we assume directly that the lasso converges (with probability
approaching one) to some fixed model defined by $A_0, s_{A_0}$ at the
$(k_0)$th step in the path. Here, $A_0$ is large enough that it
contains the true support, $A_0 \supseteq A^*$, and the signs
$s_{A_0}$ are arbitrary---they may or may not match the signs of
the true coefficients over $A_0$. In a setting in which
the nonzero coefficients in $\beta^*$ are well separated from zero,
a condition quite similar to the irrepresentable condition can be used
to show that the lasso converges to the model with support $A_0=A^*$
and\vspace*{1pt} signs $s_{A_0}=\sign(\beta_{A_0}^*)$, at step $k_0=|A_0|$ of the
path. Our result extends beyond this case, and allows for situations
in which the lasso model converges to a possibly larger set of
``screened'' variables $A_0$, and fixed signs $s_{A_0}$.
\end{rem}

%
\begin{rem}\label{rem15} In fact, one can modify the above arguments to account
for the case that $A_0$ does not contain the entire set $A^*$ of truly
nonzero coefficients, but rather, only the ``strong'' coefficients.
While ``strong'' is rather vague, a more precise way of stating this
is to assume that $\beta^*$ has nonzero coefficients both large and
small in magnitude, and with $A_0$ corresponding to the set of
large coefficients, we assume that the (left-out) small coefficients
must be small enough that the mean of the process $g$ in \eqref{eqg2}
(with $A=A_0$ and $s_A=s_{A_0}$) grows much faster than $M^+$. The
details, though not the main ideas, of the arguments would change, and
the result would still be a conservative exponential limit for the
covariance statistic $T_k$ at step $k=k_0+1$. We may pursue this
extension in future work.
\end{rem}

%

\section{\texorpdfstring{Simulation of the null distribution.}{Simulation of the null distribution}}\label{secsim}

We investigate the null distribution of the covariance statistic
through simulations, starting with an orthogonal predictor matrix $X$,
and then considering more general forms of~$X$.

%
\begin{figure}

\includegraphics{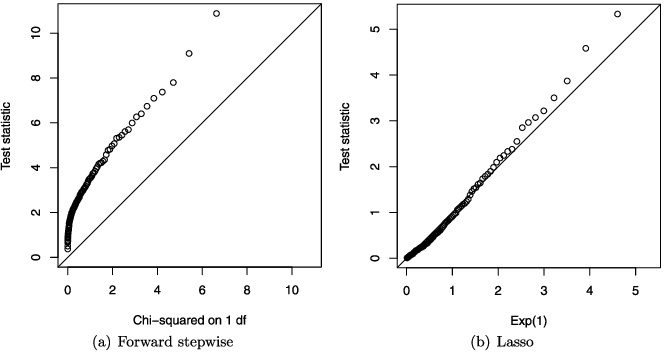}
\vspace*{-2pt}
\caption{An example with $n=100$ and $p=10$ orthogonal
predictors, and the true coefficient vector having 3 nonzero, large
components. Shown are quantile--quantile plots for the drop in $\RSS$
test applied to forward stepwise regression at the 4th step and the
covariance test for the lasso path at the 4th
step.}\label{fig2}\vspace*{-2pt}
\end{figure}

%
%
\begin{figure}[b]\vspace*{-2pt}

\includegraphics{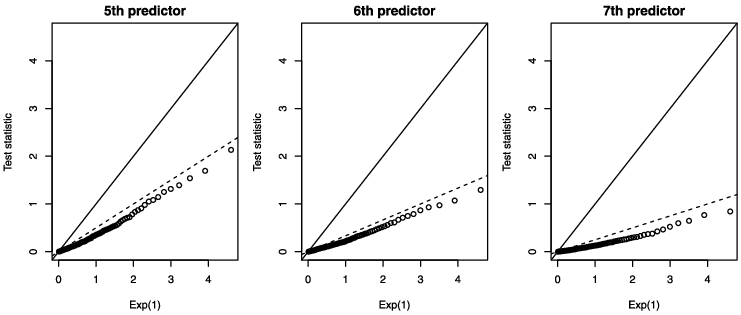}%
\vspace*{-2pt}
\caption{The same setup as in Figure~\protect\ref{fig2}, but here we
show the covariance test at the 5th, 6th
and 7th steps along the lasso path, from left to right,
respectively. The solid line has slope 1, while the broken lines have
slopes $1/2, 1/3, 1/4$, as predicted by Theorem \protect\ref
{thmorthxgen}.}\label{fig3}
\end{figure}

\subsection{\texorpdfstring{Orthogonal predictor matrix.}{Orthogonal predictor matrix}}

Similar to our example from the start of Section~\ref{seccovtest}, we
generated $n=100$ observations with $p=10$ orthogonal predictors.
The true coefficient vector $\beta^*$ contained 3 nonzero components
equal to 6, and the rest zero. The error variance was
$\sigma^2=1$, so that the truly active predictors had strong effects
and always entered the model first, with both forward stepwise and the
lasso. Figure~\ref{fig2} shows the results for testing the 4th
(truly inactive) predictor to enter, averaged over 500 simulations;
the left panel shows the chi-squared test (drop in $\RSS$) applied
at the 4th step in forward stepwise\vadjust{\goodbreak} regression, and the right panel
shows the covariance test applied at the 4th step of the lasso path.
We see that the $\Exp(1)$ distribution provides a good finite-sample
approximation for the distribution of the covariance statistic, while
$\chi^2_1$ is a poor approximation for the drop in $\RSS$.\looseness=-1

%
%
\begin{table}
\tabcolsep=0pt
\caption{Simulation results for the first
predictor to enter for a global null true model. We vary the number
of predictors $p$, correlation parameter $\rho$ and structure of the
predictor correlation matrix. Shown are the mean, variance and tail
probability $\mathbb{P}(T_1 > q_{0.95})$ of the covariance statistic $T_1$,
where $q_{0.95}$ is the $95\%$ quantile of the $\Exp(1)$ distribution,
computed over 500 simulated data sets for each setup. Standard errors
are given by ``se.''
(The panel in the bottom left corner is missing because the equal data
correlation setup is not defined for $p>n$.)}\label{tabsim}
\begin{tabular*}{\tablewidth}{@{\extracolsep{\fill}}@{}lccc ccc ccc ccc@{}}
\hline
 &\multicolumn{3}{c}{\textbf{Equal data corr}} & \multicolumn {3}{c}{\textbf{Equal pop'n corr}}
 &\multicolumn{3}{c}{$\bolds{\mathit{AR}(1)}$} & \multicolumn{3}{c@{}}{\textbf{Block diagonal}}
 \\[-6pt]
 &\multicolumn{3}{c}{\hrulefill} & \multicolumn {3}{c}{\hrulefill}
 &\multicolumn{3}{c}{\hrulefill} & \multicolumn{3}{c@{}}{\hrulefill}
 \\
$\bolds{\rho}$ & \textbf{Mean} & \textbf{Var} & \textbf{Tail pr} & \textbf{Mean} & \textbf{Var} & \textbf{Tail pr} & \textbf{Mean} & \textbf{Var}
& \textbf{Tail pr} & \textbf{Mean} & \textbf{Var} & \textbf{Tail pr} \\
\hline
\multicolumn{13}{@{}c@{}}{$n=100$, $p=10$}\\
0 & 0.966 & 1.157 & 0.062 & 1.120 & 1.951 & 0.090 & 1.017 & 1.484 & 0.070 & 1.058 & 1.548 & 0.060 \\
0.2 & 0.972 & 1.178 & 0.066 & 1.119 & 1.844 & 0.086 & 1.034 & 1.497 &0.074 & 1.069 & 1.614 & 0.078 \\
0.4 & 0.963 & 1.219 & 0.060 & 1.115 & 1.724 & 0.092 & 1.045 & 1.469 &0.060 & 1.077 & 1.701 & 0.076 \\
0.6 & 0.960 & 1.265 & 0.070 & 1.095 & 1.648 & 0.086 & 1.048 & 1.485 &0.066 & 1.074 & 1.719 & 0.086 \\
0.8 & 0.958 & 1.367 & 0.060 & 1.062 & 1.624 & 0.092 & 1.034 & 1.471 &0.062 & 1.062 & 1.687 & 0.072 \\
se & 0.007 & 0.015 & 0.001 & 0.010 & 0.049 & 0.001 & 0.013 & 0.043 &0.001 & 0.010 & 0.047 & 0.001
\\[3pt]
\multicolumn{13}{@{}c@{}}{$n=100$, $p=50$}\\
0 & 0.929 & 1.058 & 0.048 & 1.078 & 1.721 & 0.074 & 1.039 & 1.415 &0.070 & 0.999 & 1.578 & 0.048 \\
0.2 & 0.920 & 1.032 & 0.038 & 1.090 & 1.476 & 0.074 & 0.998 & 1.391 &0.054 & 1.064 & 2.062 & 0.052 \\
0.4 & 0.928 & 1.033 & 0.040 & 1.079 & 1.382 & 0.068 & 0.985 & 1.373 &0.060 & 1.076 & 2.168 & 0.062 \\
0.6 & 0.950 & 1.058 & 0.050 & 1.057 & 1.312 & 0.060 & 0.978 & 1.425 &0.054 & 1.060 & 2.138 & 0.060 \\
0.8 & 0.982 & 1.157 & 0.056 & 1.035 & 1.346 & 0.056 & 0.973 & 1.439 &0.060 & 1.046 & 2.066 & 0.068 \\
se & 0.010 & 0.030 & 0.001 & 0.011 & 0.037 & 0.001 & 0.009 & 0.041 &0.001 & 0.011 & 0.103 & 0.001
\\[3pt]
\multicolumn{13}{@{}c@{}}{$n=100$, $p=200$}\\
0 & & & & 1.004 & 1.017 & 0.054 & 1.029 & 1.240 & 0.062 & 0.930 & 1.166& 0.042 \\
0.2 & & & & 0.996 & 1.164 & 0.052 & 1.000 & 1.182 & 0.062 & 0.927 &1.185 & 0.046 \\
0.4 & & & & 1.003 & 1.262 & 0.058 & 0.984 & 1.016 & 0.058 & 0.935 &1.193 & 0.048 \\
0.6 & & & & 1.007 & 1.327 & 0.062 & 0.954 & 1.000 & 0.050 & 0.915 &1.231 & 0.044 \\
0.8 & & & & 0.989 & 1.264 & 0.066 & 0.961 & 1.135 & 0.060 & 0.914 &1.258 & 0.056 \\
se & & & & 0.008 & 0.039 & 0.001 & 0.009 & 0.028 & 0.001 & 0.007 &0.032 & 0.001 \\
\hline
\end{tabular*}
\end{table}

Figure~\ref{fig3} shows the results for testing the 5th, 6th and 7th
predictors to enter the lasso model. An $\Exp(1)$-based test will now
be conservative: at a nominal $5\%$ level, the actual type I errors
are about $1\%$, $0.2\%$ and $0.0\%$, respectively. The solid line has
slope 1, and the broken lines have slopes $1/2,1/3,1/4$,
as predicted by Theorem~\ref{thmorthxgen}.

\subsection{\texorpdfstring{General predictor matrix.}{General predictor matrix}}

In Table~\ref{tabsim}, we simulated null data (i.e., \mbox{$\beta^*=0$}), and
examined the distribution of the covariance test statistic $T_1$ for
the first predictor to enter. We varied the numbers of predictors $p$,
correlation parameter $\rho$, and structure of the predictor correlation
matrix. In the first two correlation setups, the correlation between
each pair of predictors was $\rho$, in the data and population,
respectively. In the $\AR(1)$ setup, the correlation between predictors
$j$~and~$j'$ is $\rho^{|j-j'|}$. Finally, in the block diagonal setup,
the correlation matrix has two equal-sized blocks, with population
correlation $\rho$ in each block. We computed the mean, variance and
tail probability of the covariance statistic $T_1$ over $500$
simulated data sets for each setup. We see that the $\Exp(1)$
distribution is a reasonably good approximation throughout.\looseness=1

In Table~\ref{tabsim2}, the setup was the same as in Table~\ref
{tabsim}, except that we set the first~$k$ coefficients of the
true coefficient vector equal to 4, and the rest zero, for $k=1,2,3$.
The dimensions were also fixed at $n=100$ and $p=50$.
We computed the mean, variance, and tail probability of the covariance
statistic $T_{k+1}$ for entering the next (truly inactive) $(k+1)$st
predictor, discarding those simulations in which a truly inactive
predictor was selected in the first $k$ steps. (This occurred $1.7\%$,
$4.0\%$ and $7.0\%$ of the time, resp.) Again, we see that
the $\Exp(1)$ approximation is reasonably accurate throughout.

%
%
\begin{table}
\tabcolsep=0pt
\caption{Simulation results for the
$(k+1)$st predictor to enter for a model with $k$ truly nonzero
coefficients, across $k=1,2,3$. The rest of the setup is the same as
in Table~\protect\ref{tabsim} except that the dimensions were fixed at
$n=100$ and $p=50$.
The values are conditional on the event that the $k$ truly active
variables enter in the first $k$ steps}\label{tabsim2}
\begin{tabular*}{\tablewidth}{@{\extracolsep{\fill}}@{}lccc ccc ccc ccc@{}}
\hline
&\multicolumn{3}{c}{\textbf{Equal data corr}} & \multicolumn{3}{c}{\textbf{Equal pop'n corr}} & \multicolumn{3}{c}{$\bolds{\mathit{AR}(1)}$}
& \multicolumn{3}{c@{}}{\textbf{Block diagonal}}
\\[-6pt]
 &\multicolumn{3}{c}{\hrulefill} & \multicolumn {3}{c}{\hrulefill}
 &\multicolumn{3}{c}{\hrulefill} & \multicolumn{3}{c@{}}{\hrulefill}
 \\
$\bolds{\rho}$ & \textbf{Mean} & \textbf{Var} & \textbf{Tail pr} & \textbf{Mean} & \textbf{Var} & \textbf{Tail pr} & \textbf{Mean} & \textbf{Var}
& \textbf{Tail pr} & \textbf{Mean} & \textbf{Var} & \textbf{Tail pr} \\
\hline
\multicolumn{13}{@{}c@{}}{$k=1$ and 2nd predictor to enter}\\
0 & 0.933 & 1.091 & 0.048 & 1.105 & 1.628 & 0.078 & 1.023 & 1.146 &0.064 & 1.039 & 1.579 & 0.060 \\
0.2 & 0.940 & 1.051 & 0.046 & 1.039 & 1.554 & 0.082 & 1.017 & 1.175 &0.060 & 1.062 & 2.015 & 0.062 \\
0.4 & 0.952 & 1.126 & 0.056 & 1.016 & 1.548 & 0.084 & 0.984 & 1.230 &0.056 & 1.042 & 2.137 & 0.066 \\
0.6 & 0.938 & 1.129 & 0.064 & 0.997 & 1.518 & 0.079 & 0.964 & 1.247 &0.056 & 1.018 & 1.798 & 0.068 \\
0.8 & 0.818 & 0.945 & 0.039 & 0.815 & 0.958 & 0.044 & 0.914 & 1.172 &0.062 & 0.822 & 0.966 & 0.037 \\
se & 0.010 & 0.024 & 0.002 & 0.011 & 0.036 & 0.002 & 0.010 & 0.030 &0.002 & 0.015 & 0.087 & 0.002
\\[3pt]
\multicolumn{13}{@{}c@{}}{$k=2$ and 3rd predictor to enter}\\
0 & 0.927 & 1.051 & 0.046 & 1.119 & 1.724 & 0.094 & 0.996 & 1.108 &0.072 & 1.072 & 1.800 & 0.064 \\
0.2 & 0.928 & 1.088 & 0.044 & 1.070 & 1.590 & 0.080 & 0.996 & 1.113 &0.050 & 1.043 & 2.029 & 0.060 \\
0.4 & 0.918 & 1.160 & 0.050 & 1.042 & 1.532 & 0.085 & 1.008 & 1.198 &0.058 & 1.024 & 2.125 & 0.066 \\
0.6 & 0.897 & 1.104 & 0.048 & 0.994 & 1.371 & 0.077 & 1.012 & 1.324 &0.058 & 0.945 & 1.568 & 0.054 \\
0.8 & 0.719 & 0.633 & 0.020 & 0.781 & 0.929 & 0.042 & 1.031 & 1.324 &0.068 & 0.771 & 0.823 & 0.038 \\
se & 0.011 & 0.034 & 0.002 & 0.014 & 0.049 & 0.003 & 0.009 & 0.022 &0.002 & 0.013 & 0.073 & 0.002
\\[3pt]
\multicolumn{13}{@{}c@{}}{$k=3$ and 4th predictor to enter}\\
0 & 0.925 & 1.021 & 0.046 & 1.080 & 1.571 & 0.086 & 1.044 & 1.225 &0.070 & 1.003 & 1.604 & 0.060 \\
0.2 & 0.926 & 1.159 & 0.050 & 1.031 & 1.463 & 0.069 & 1.025 & 1.189 &0.056 & 1.010 & 1.991 & 0.060 \\
0.4 & 0.922 & 1.215 & 0.048 & 0.987 & 1.351 & 0.069 & 0.980 & 1.185 &0.050 & 0.918 & 1.576 & 0.053 \\
0.6 & 0.905 & 1.158 & 0.048 & 0.888 & 1.159 & 0.053 & 0.947 & 1.189 &0.042 & 0.837 & 1.139 & 0.052 \\
0.8 & 0.648 & 0.503 & 0.008 & 0.673 & 0.699 & 0.026 & 0.940 & 1.244 &0.062 & 0.647 & 0.593 & 0.015 \\
se & 0.014 & 0.037 & 0.002 & 0.016 & 0.044 & 0.003 & 0.014 & 0.031 &0.003 & 0.016 & 0.073 & 0.002 \\
\hline
\end{tabular*}
\end{table}

\begin{figure}

\includegraphics{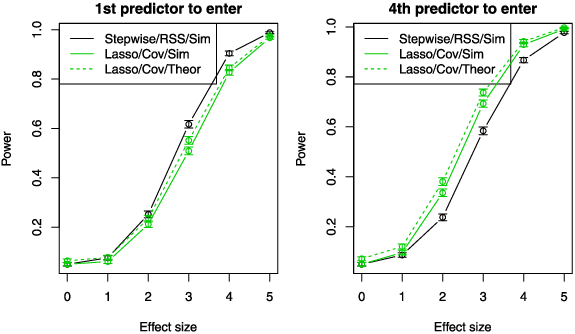}

\caption{Estimated power curves for significance tests using
forward stepwise regression and the drop in RSS statistic, as well as
the lasso and the covariance statistic. The results are averaged
over 1000 simulations with $n=100$ and $p=10$ predictors drawn
i.i.d. from $N(0,1)$ and $\sigma^2=1$.
On the left, there is one truly nonzero regression coefficient, and we
varied its magnitude (the effect size parameter on the $x$-axis). We
examined the first step of the forward stepwise and lasso
procedures. On the right, in addition to a nonzero coefficient
with varying effect size (on the $x$-axis), there are 3 additional large
coefficients in the true model. We examined the 4th step in forward
stepwise and the lasso, after the 3 strong variables have been
entered. For the power curves in both panels, we use simulation-based
cutpoints for forward stepwise to control the type I error at the 5\%
level; for the lasso we do the same, but also display the results for
the theoretically-based [$\Exp(1)$] cutpoint.
Note that in practice, simulation-based cutpoints would not
typically be available.}\label{figpow}
\end{figure}

In Figure~\ref{figpow}, we estimate the power curves for significance
testing via the drop in RSS test for forward stepwise regression, and
the covariance test for the lasso. In the former, we use
simulation-derived cutpoints, and in the latter we use the
theoretically-based $\Exp(1)$ cutpoints, to control the type I error
at the 5\% level. We find that the tests have similar power, though
the cutpoints for forward stepwise would not be typically available in
practice. For more details, see the figure caption.

%
%

\section{\texorpdfstring{The case of unknown $\sigma^2$.}{The case of unknown sigma2}}
\label{secsigma}

Up until now, we have assumed that the error variance $\sigma^2$ is
known; in practice it will typically be unknown. In this case,
provided that $n>p$, we can easily estimate it and proceed by
analogy to standard linear model
theory. In particular,\vspace*{2pt} we can estimate $\sigma^2$ by the mean squared
residual error $\hat{\sigma}^2 =
\|y-X\hbeta^\mathrm{LS}\|_2^2/(n-p)$, with $\hbeta^\mathrm{LS}$ being
the regression coefficients from $y$ on $X$ (i.e., the full model).
Plugging this estimate into the covariance statistic in
\eqref{eqcovtest} yields a new statistic $F_k$ that has an
asymptotic $F$-distribution under the null:
%
%
\begin{equation}
\label{eqcvestsig} F_k = \frac{ \langle y, X\hbeta(\lambda_{k+1})
\rangle
- \langle y,X_A \tbeta_A(\lambda_{k+1}) \rangle} {
\hat{\sigma}^2} \stackrel{d} {
\rightarrow}F_{2,n-p}.
\end{equation}
This follows because $F_k = T_k / (\hat{\sigma}^2/\sigma^2)$, the
numerator $T_k$ being asymptotically $\Exp(1) = \chi_2^2/2$, the
denominator $\hat{\sigma}^2/\sigma^2$ being asymptotically
$\chi_{n-p}^2/\break (n-p)$, and we claim that the two are independent. Why?
Note that the lasso solution path is unchanged if we replace $y$ by
$P_Xy$, so the lasso fitted values in $T_k$ are functions of $P_Xy$;
meanwhile, $\hat{\sigma}^2$ is a function of $(I-P_X)y$. The
quantities $P_Xy$ and $(I-P_X)y$ are uncorrelated, and hence
independent (recalling normality of $y$), so $T_k$ and
$\hat{\sigma}^2$ are functions of independent quantities and,
therefore, independent.

As an example, consider one of the setups from Table~\ref{tabsim},
with $n=100$, $p=80$ and predictor correlation of the $\AR(1)$ form
$\rho^{|j-j'|}$. The true
model is null, and we test the first predictor to enter along
the lasso path. (We choose $n,p$ of roughly equal sizes here to expose
the differences between the $\sigma^2$ known and unknown cases.)
Table~\ref{tabsim3} shows the results of 1000 simulations from each
of the $\rho=0$ and $\rho=0.8$ scenarios. We see that with $\sigma^2$
estimated, the $F_{2,n-p}$ distribution provides a more accurate
finite-sample approximation than does $\Exp(1)$.\vspace*{1pt}


%
%
\begin{table}
\tabcolsep=0pt
\caption{Comparison\vspace*{-1pt} of $\Exp(1)$, $F_{2,{N-p}}$,
and the observed (empirical) null distribution of the covariance
statistic, when $\sigma^2$ has been estimated. We examined 1000
simulated data sets with $n=100$, $p=80$ and the correlation
between predictors $j$ and $j'$ equal to $\rho^{|j-j'|}$. We are
testing the first step of the lasso path, and the true model is the
global null. Results are shown for $\rho=0.0$ and $0.8$. The third
column shows the tail probability $\mathbb{P}(T_1 > q_{0.95})$
computed over the 1000 simulations, where $q_{0.95}$ is the $95\%$ quantile from the
appropriate distribution [either $\Exp(1)$ or $F_{2,n-p}$]}\label{tabsim3}
\begin{tabular*}{\tablewidth}{@{\extracolsep{\fill}}@{}lcccc@{}}
\hline
& \textbf{Mean} & \textbf{Variance} & \textbf{95\% quantile} & \textbf{Tail prob}\\
\hline
\multicolumn{5}{@{}c@{}}{$\rho=0$}\\
Observed & 1.17 & 2.10 & 3.75 & \\
$\Exp(1)$ & 1.00 & 1.00& 2.99& 0.082\\
$F_{2,n-p}$ & 1.11 & 1.54 & 3.49 & 0.054
\\[3pt]
\multicolumn{5}{@{}c@{}}{$\rho=0.8$}\\
Observed & 1.14& 1.70 & 3.77 & \\
$\Exp(1)$ & 1.00 & 1.00& 2.99& 0.097\\
$F_{2,n-p}$ & 1.11& 1.54& 3.49& 0.064\\
\hline
\end{tabular*}
\end{table}

When $p \geq n$, estimation of $\sigma^2$ is not nearly as
straightforward; one idea is to estimate $\sigma^2$ from the
least squares fit on the support of the model selected by
cross-validation. One would then hope that the resulting statistic,
with this plug-in estimate of $\sigma^2$, is close in distribution to
$F_{2,n-r}$ under the null, where $r$ is the size of the model chosen
by cross-validation. This is by analogy to the low-dimensional $n >
p$ case in \eqref{eqcvestsig}, but is not supported by rigorous
theory. Simulations (withheld for brevity) show that this
approximation is not too far off, but that the variance of the
observed statistic is sometimes inflated compared that of an
$F_{2,n-r}$ distribution (this unaccounted variability is likely due
to the model selection process via cross-validation). Other authors
have argued that using cross-validation to estimate $\sigma^2$ when
$p \gg n$ is not
necessarily a good approach, as it can be anti-conservative; see,
for example, \citet{varestbycv}, \citet{scaledlasso}
for alternative techniques. In future work, we will
address the important issue of estimating $\sigma^2$ in the context of
the covariance statistic, when $p \geq n$.

\section{\texorpdfstring{Real data examples.}{Real data examples}}\label{secrealdata}

We demonstrate the use of covariance test with some real data
examples. As mentioned previously, in any serious application
of significance testing over many variables (many steps of the lasso
path), we would need to consider the issue of multiple comparisons,
which we do not here. This is a topic for future work.

\subsection{\texorpdfstring{Wine data.}{Wine data}}

Table~\ref{tabwine} shows the results for the wine quality data
taken from the UCI database. There are $p=11$ predictors, and
$n=1599$ observations, which we split
randomly into approximately equal-sized training and test sets. The
outcome is a wine quality rating, on a scale between 0 and 10.
The table shows the training set $p$-values from forward stepwise
regression (with the chi-squared test) and the lasso (with the
covariance test). Forward stepwise enters 6
predictors at the 0.05 level, while the lasso enters only 3.

%
%
\begin{table}[t]
\tabcolsep=0pt
\caption{Wine data: forward stepwise and lasso $p$-values.
The values are rounded to 3 decimal places. For the lasso, we only
show $p$-values for the steps in which a predictor entered the model
and stayed in the model for the remainder of the path (i.e., if a
predictor entered the model at a step but then later left, we do not
show this step---we only show the step corresponding to its last entry
point)}\label{tabwine}
\begin{tabular*}{\tablewidth}{@{\extracolsep{\fill}}@{}lcd{3.3}cd{2.0}cd{2.3}c@{}}
\hline
\multicolumn{4}{@{}c}{\textbf{Forward stepwise}} & \multicolumn{4}{c@{}}{\textbf{Lasso}} \\[-6pt]
\multicolumn{4}{@{}c}{\hrulefill} & \multicolumn{4}{c@{}}{\hrulefill} \\
\textbf{Step} & \textbf{Predictor} & \multicolumn{1}{c}{\textbf{RSS test}}
& \multicolumn{1}{c}{\textbf{$\bolds{p}$-value}} & \multicolumn{1}{c}{\textbf{Step}}
& \multicolumn{1}{c}{\textbf{Predictor}} & \multicolumn{1}{c}{\textbf{Cov test}} & \textbf{$\bolds{p}$-value}\\
\hline
\phantom{0}1 & Alcohol & 315.216 & 0.000 & 1 & Alcohol & 79.388 & 0.000 \\
\phantom{0}2 & Volatile\_acidity & 137.412 & 0.000 & 2 & Volatile\_acidity &77.956 & 0.000 \\
\phantom{0}3 & Sulphates & 18.571 & 0.000 & 3 & Sulphates & 10.085 & 0.000 \\
\phantom{0}4 & Chlorides & 10.607 & 0.001 & 4 & Chlorides & 1.757 & 0.173 \\
\phantom{0}5 & pH & 4.400 & 0.036 & 5 & Total\_sulfur\_dioxide & 0.622 & 0.537 \\
\phantom{0}6 & Total\_sulfur\_dioxide & 3.392 & 0.066 & 6 & pH & 2.590 & 0.076 \\
\phantom{0}7 & Residual\_sugar & 0.607 & 0.436 & 7 & Residual\_sugar & 0.318 &0.728 \\
\phantom{0}8 & Citric\_acid & 0.878 & 0.349 & 8 & Citric\_acid & 0.516 & 0.597 \\
\phantom{0}9 & Density & 0.288 & 0.592 & 9 & Density & 0.184 & 0.832 \\
10 & Fixed\_acidity & 0.116 & 0.733 & 10 & Free\_sulfur\_dioxide &0.000 & 1.000 \\
11 & Free\_sulfur\_dioxide & 0.000 & 0.997 & 11 & Fixed\_acidity &0.114 & 0.892 \\
\hline
\end{tabular*}
\end{table}
%

\begin{figure}

\includegraphics{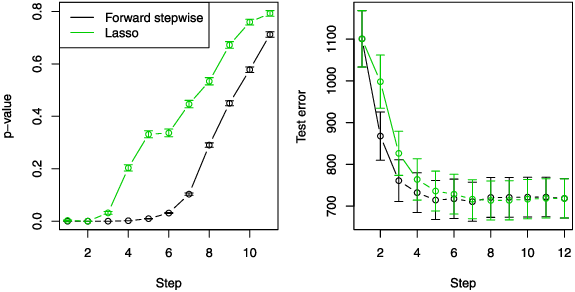}

\caption{Wine data: the data were randomly divided
500 times into roughly equal-sized training and test sets. The left
panel shows the training set $p$-values for forward stepwise regression
and the lasso. The right panel show the test set error for the
corresponding models of each size.}\label{figwine}
\end{figure}

In the left panel of Figure~\ref{figwine}, we repeated this $p$-value
computation over 500 random splits into training test sets.
The right panel shows the corresponding test set prediction error
for the models of each size. The lasso test error
decreases sharply once the 3rd predictor is added, but then somewhat
flattens out from the 4th predictor onward; this is in general
qualitative agreement with the lasso $p$-values in the left panel, the
first 3 being very small, and the 4th $p$-value being about 0.2.
This also echoes the well-known difference between hypothesis testing
and minimizing prediction error. For example, the $C_p$ statistic
stops entering variables when the $p$-value is larger than about
0.16.

%
%

%
%
\begin{figure}[b]

\includegraphics{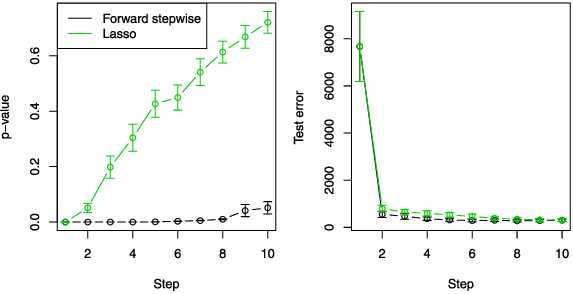}

\caption{HIV data: the data were randomly divided 50 times
into training and test sets of size 150 and 907, respectively. The
left panel shows the training set $p$-values for forward stepwise
regression and the lasso. The right panel shows the test set error
for the corresponding models of each size.}\label{fignrti}
\end{figure}

\subsection{\texorpdfstring{HIV data.}{HIV data}}

Rhee et al. (\citeyear{rhee2003}) study six nucleotide reverse transcriptase
inhibitors (NRTIs) that are used to treat HIV-1. The target of these
drugs can become resistant through mutation, and they compare a
collection of models for predicting the (log)
susceptibility of the drugs, a measure of drug resistance, based on
the location of mutations. We focused on the first drug (3TC), for
which there are $p = 217$ sites and $n=1057$
samples. To examine the behavior of the covariance test in the $p>n$
setting, we divided the data at random into training and test sets of
size 150 and 907, respectively, a total of 50
times. Figure~\ref{fignrti} shows the results, in the same format as
Figure~\ref{figwine}. We used the model chosen by cross-validation to
estimate $\sigma^2$. The covariance test for the lasso suggests that
there are only one or two important predictors (in marked contrast
to the chi-squared test for forward stepwise), and this is confirmed
by the test error plot in the right panel.

\section{\texorpdfstring{Extensions.}{Extensions}}\label{secextensions}

We discuss some extensions of the covariance statistic, beyond
significance testing for the lasso. The proposals here are
supported by simulations [in terms of having an $\Exp(1)$ null
distribution], but we do not offer any theory. This may be a direction for
future work.

\subsection{\texorpdfstring{The elastic net.}{The elastic net}}

\begin{figure}

\includegraphics{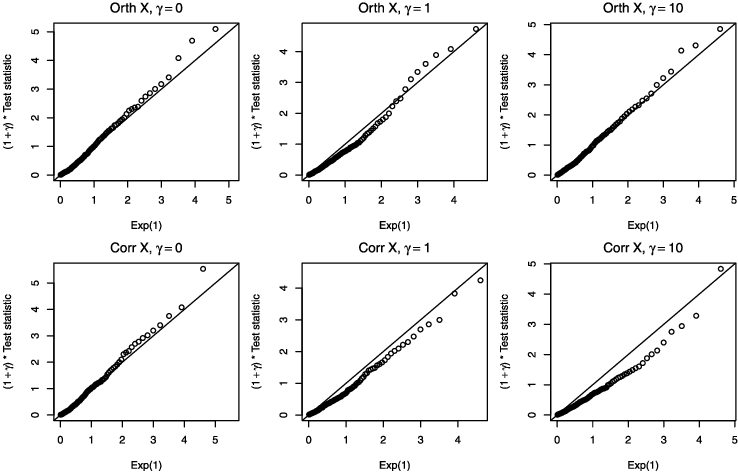}

\caption{Elastic net: an example with $n=100$ and
$p=10$, for orthogonal and correlated predictors (having pairwise
population correlation 0.5), and three different values of the ridge
penalty parameter $\gamma$.}\label{figen}
\end{figure}

The elastic net estimate [\citet{enet}] is defined as
%
%
\begin{equation}
\label{eqenet} \hbeta^\mathrm{en} = \argmin_{\beta\in\mathbb{R}^p}
\frac{1}{2}\|y-X\beta\|_2^2 + \lambda\|\beta
\|_1 + \frac{\gamma}{2} \|\beta\|_2^2,
\end{equation}
where $\gamma\geq0$ is a second tuning parameter. It is not hard to
see that this can actually be cast as a lasso estimate with
predictor\vadjust{\goodbreak}
matrix $\widetilde{X} =\bigl[{\fontsize{8.36pt}{10pt}\selectfont{
\begin{array}{c} X \\ \sqrt{\gamma} I
\end{array}}}%
\bigr] \in\mathbb{R}^{(n+p)\times p}$ and outcome $\ty= (y,0) \in
\mathbb{R}^{n+p}$. This shows that, for a fixed $\gamma$, the elastic net
solution path is piecewise linear over $\lambda$, with each knot
marking the entry (or deletion) of a variable from the active set.
We therefore define the covariance statistic in the same manner as we
did for the lasso; fixing $\gamma$, to test the predictor entering at
the $k$th step (knot $\lambda_k$) in the elastic net path, we consider
the statistic
\[
T_k 
= \bigl( \bigl\langle y, X\hbeta^\mathrm{en}(
\lambda_{k+1},\gamma) \bigr\rangle- \bigl\langle y, X_A
\tbeta_A^\mathrm{en}( \lambda_{k+1},\gamma) \bigr
\rangle \bigr) / \sigma^2,
\]
where\vspace*{1pt} as before, $\lambda_{k+1}$ is next knot in the path,
$A$ is the active set of predictors just before $\lambda_k$ and
$\tbeta^\mathrm{en}_A$ is the elastic net solution using only the
predictors $X_A$. The precise expression for the elastic net solution
in \eqref{eqenet}, for a given active set and signs, is the same
as it is for the lasso (see Section~\ref{secalternate}), but with
$(X_A^T X_A)^{-1}$ replaced by $(X_A^T X_A + \gamma I)^{-1}$. This
generally creates a complication for the theory in Sections~\ref
{secorthx} and~\ref{secgenx}. But in the orthogonal $X$ case, we
have $(X_A^T X_A + \gamma I)^{-1} = I/(1+\gamma)$ and so
\[
T_k = 1/(1+\gamma) \cdot|U_{(k)}|\bigl(|U_{(k)}|-|U_{(k+1)}|\bigr)/
\sigma^2
\]
with $U_j=X_j^T y$, $j=1,\ldots, p$. This means that for an orthogonal
$X$, under the null,
\[
(1+\gamma) \cdot T_k \stackrel{d} {\rightarrow}\Exp(1)
\]
and one is tempted to use this approximation beyond the orthogonal
setting as well.
In Figure~\ref{figen}, we evaluated the distribution of
$(1+\gamma) T_1$ (for the first predictor to enter), for orthogonal and
correlated scenarios, and for three different values of~$\gamma$.
Here, $n=100$, $p=10$ and the true model was null. It seems to be
reasonably close to $\Exp(1)$ in all cases.

%
%

\subsection{\texorpdfstring{Generalized linear models and the Cox model.}{Generalized linear models and the Cox model}}

Consider the estimate from an $\ell_1$-penalized generalized linear
model:
%
%
\begin{equation}
\label{eqlassoglm} \hbeta^\mathrm{glm} = \argmin_{\beta\in\mathbb{R}^p} -\sum
_{i=1}^n \log f(y_i;
x_i, \beta) + \lambda\|\beta\|_1,
\end{equation}
where $f(y_i; x_i, \beta)$ is an exponential family density, a
function of the predictor measurements $x_i\in\mathbb{R}^p$ and parameter
$\beta\in\mathbb{R}^p$. Note that
the usual lasso estimate in \eqref{eqlasso} is a special case
of this form when $f$ is the Gaussian density with known variance
$\sigma^2$. The natural parameter in \eqref{eqlassoglm} is $\eta_i =
x_i^T \beta$, for $i=1,\ldots, n$, related to the mean of $y_i$ via a
link function $g(\mathbb{E}[y_i | x_i]) = \eta_i$.

%
%
\begin{figure}

\includegraphics{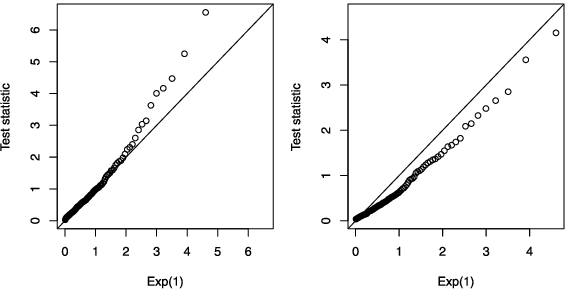}

\caption{Lasso logistic regression: an example
with $n=100$ and $p=10$ predictors, i.i.d. from $N(0,1)$. In
the left panel, all true coefficients are zero; on the right, the
first coefficient is large, and the rest are zero. Shown are
quantile--quantile plots of the covariance test statistic (at the first
and second steps, resp.), generated over 500 data sets, versus
its conjectured asymptotic distribution, $\Exp(1)$.}\label{figlogit}
\end{figure}

Having solved \eqref{eqlassoglm} with
$\lambda=0$ (i.e., this is simply maximum likelihood), producing a
vector of fitted values $\hat{\eta}=X\hbeta^\mathrm{glm}\in\mathbb
{R}^n$, we
might define degrees of freedom as\footnote{Note that in the Gaussian
case, this definition is actually $\sigma^2$ times the usual notion of
degrees of freedom; hence in the presence of a scale parameter, we
would divide the right-hand side in the definition \eqref{eqdfglm} by
this scale parameter, and we would do the same for the covariance
statistic as defined in \eqref{eqdfglm}.}
%
%
\begin{equation}
\label{eqdfglm} \df(\hat{\eta})=\sum_{i=1}^n
\Cov(y_i,\hat{\eta}_i).
\end{equation}
%
This is the implicit concept used by \citet{bradbiased} in his
definition of the ``optimism'' of the training error. The same idea
could be used to define degrees of freedom for the penalized estimate
in \eqref{eqlassoglm} for any $\lambda>0$, and this motivates the
definition of the covariance statistic, as follows. If the tuning
parameter value $\lambda=\lambda_k$ marks the entry of a new predictor
into the active set $A$, then we define the covariance statistic
%
%
\begin{equation}
\label{eqcovtestglm} T_k = \bigl\langle y, X\hbeta^\mathrm{glm}(
\lambda_{k+1}) \bigr\rangle- \bigl\langle y, X_A
\tbeta_A^\mathrm{glm}(\lambda_{k+1}) \bigr\rangle,
\end{equation}
%
where $\lambda_{k+1}$ is the next value of the tuning parameter at
which the model changes (a variable enters or leaves the active set),
and $\tbeta_A^\mathrm{glm}$ is the estimate from the penalized
generalized linear model \eqref{eqlassoglm} using only predictors in
$A$.
Unlike in the Gaussian case, the solution path in \eqref{eqlassoglm}
is not generally piecewise linear over $\lambda$, and there
is not an algorithm
to deliver the exact the values of $\lambda$ at which variables enter
the model (we still refer to these as knots in the path). However, one
can numerically approximate these knot values; for example, see
\citet{glmpath}.
By analogy to the Gaussian case, we would hope that $T_k$ has an
asymptotic $\Exp(1)$ distribution under the null. Though we have not
rigorously investigated this conjecture, simulations seem to support
it.

As example, consider the logistic regression model for binary
data. Now $\eta_i=\log(\mu_i/(1-\mu_i))$, with
$\mu_i=\mathbb{P}(y_i=1|x_i)$. Figure~\ref{figlogit} shows the simulation
results from comparing the null distribution of the covariance test
statistic in \eqref{eqcovtestglm} to $\Exp(1)$. Here, we used the
{\tt
glmpath} package in R [\citet{glmpath}] to compute an approximate
solution path and locations of knots. The null distribution of the
test statistic looks fairly close to $\Exp(1)$.

\begin{figure}

\includegraphics{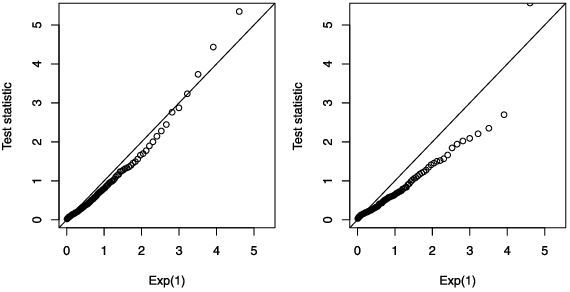}

\caption{Lasso Cox model estimate: the basic
setup is the same as in Figure~\protect\ref{figlogit} ($n$, $p$, the
distribution of the predictors $X$, the true coefficient
vector---on the left, entirely zero, and on the right, one large
coefficient). Shown are quantile--quantile plots of the covariance
test statistic (at the first and second steps, resp.),
generated over 500 data sets, versus the $\Exp(1)$
distribution.}\label{figcox}
\end{figure}

For general likelihood-based regression problems,
let $\eta=X\beta$ and $\ell(\eta)$ denote the log likelihood.
We can view maximum likelihood estimation as an iteratively weighted
least squares procedure using the outcome variable
%
%
\begin{equation}
z(\eta)= \eta+I_\eta^{-1} S_\eta,
\end{equation}
where $S_\eta=\nabla\ell(\eta)$, and
$I_\eta=\nabla^2 \ell(\eta)$. This applies, for example, to the class
of generalized linear models and Cox's proportional hazards model. For
the general $\ell_1$-penalized estimator
%
%
\begin{equation}
\label{eqlassolik} \hbeta^\mathrm{lik} = \argmin_{\beta\in\mathbb{R}^p} -\ell(X
\beta) + \lambda\|\beta\|_1,
\end{equation}
we can analogously define the covariance test statistic at a knot
$\lambda_k$, marking the entry of a predictor into the active set $A$,
as
%
%
\begin{equation}
\label{eqcovtestlik} T_k = \bigl( \bigl\langle I_0^{-1/2}
S_0, X\hbeta^\mathrm{lik}(\lambda_{k+1}) \bigr
\rangle- \bigl\langle I_0^{-1/2} S_0,
X_A\tbeta^\mathrm{lik}_A(\lambda_{k+1})
\bigr\rangle \bigr)/2
\end{equation}
with $\lambda_{k+1}$ being the next knot in the path (at which a
variable is added or deleted from the active set), and
$\tbeta_A^\mathrm{lik}$ the solution of the general penalized
likelihood problem \eqref{eqlassolik} with predictor matrix $X_A$.
For the binomial model, the statistic \eqref{eqcovtestlik} \mbox{reduces} to
expression \eqref{eqcovtestglm}. In Figure~\ref{figcox}, we computed
this statistic for Cox's proportional hazards model, using a similar
setup to that in Figure~\ref{figlogit}. The $\Exp(1)$ approximation
for its null distribution looks reasonably accurate.

%
%

\section{\texorpdfstring{Discussion.}{Discussion}}
\label{secdiscussion}

We proposed a simple \textit{covariance statistic} for testing the
significance of predictor variables as they enter the
active set, along the lasso solution path. We showed that the\vadjust{\goodbreak}
distribution of this statistic is asymptotically $\Exp(1)$, under the
null hypothesis that all truly active predictors are contained in the
current active set. [See Theorems~\ref{thmorthxgen},
\ref{thmgenxfirst} and
\ref{thmgenxgen}; the conditions required for this convergence result
vary depending on the step $k$ along the path that we are considering,
and the covariance structure of the predictor matrix $X$; the
$\Exp(1)$ limiting distribution is in some cases a conservative upper
bound under the null.] Such a result accounts for the adaptive nature
of the lasso procedure, which is not true for the usual chi-squared
test (or $F$-test) applied to, for example, forward stepwise regression.


We feel that our work has shed light not only on the lasso path (as
given by LARS), but also, at a high level, on forward stepwise
regression. Both the lasso and forward stepwise start by
entering the predictor variable most
correlated with the outcome (thinking of standardized predictors), but
the two differ in what they do next. Forward stepwise is greedy, and
once it enters this first variable, it proceeds to fit the first
coefficient fully, ignoring the effects of other predictors. The
lasso, on the
other hand, increases (or decreases) the coefficient of the first
variable only as long as its correlation with the residual is larger
than that of the inactive predictors.
Subsequent steps follow similarly.
Intuitively, it seems that
forward stepwise regression inflates coefficients unfairly, while the
lasso takes more appropriately sized steps. This intuition is
confirmed in one sense by looking at degrees of freedom (recall
Section~\ref{secdf}). The covariance test and its simple asymptotic null
distribution reveal another way in which the step sizes used by the
lasso are ``just right.''


The problem of assessing significance in an adaptive linear model fit
by the lasso is a difficult one, and what we have presented in this paper
by no means a complete solution.
We describe some current work and ideas for future projects
below.

\begin{itemize}
\item\textit{Significance test for generic lasso models.}
A natural direction to consider is the generic lasso
testing problem: given a lasso model computed at some fixed value of
$\lambda$, how do we carry out a significance test for each predictor
in the active set? Work on this is in progress.

\item\textit{Nonasymptotic null distributions.}
A geometric characterization of the first knot in the lasso path
provides an alternative test for the global null hypothesis, \mbox{$\beta^*=0$.}
When all predictors have unit norm, $\|X_i\|_2=1$, for $i=1,\ldots, p$, this
test has the form
\[
\frac{1-\Phi(\lambda_1/\sigma)}{1-\Phi(\lambda_2/\sigma)} \sim \operatorname{Unif}(0,1).
\]
Remarkably, this above result is exact (nonasymptotic),
valid for any $n$ and $p$, requiring (essentially) only Gaussianity of
the errors, and no real assumptions about the matrix $X$.
For most reasonably behaved predictor matrices $X$, the $\Exp(1)$
approximation agrees closely with this test. Details are in
\citet{geomsignif}. Work to extend this formula to subsequent steps
along the solution path, that is, to test a hypothesis beyond the global
null, is underway.

\item\textit{Generalizations to other penalties and models.}
The manuscript of \citet{geomsignif} applies to a
regularized regression setting with a \mbox{general} seminorm penalty,
and derives explicit results for the group lasso and nuclear norm
penalties (in addition to the lasso penalty). The nuclear norm result
yields a test for principal components and matrix completion.
The recent work of \citet{glassosignif} studies the covariance
test for graphical models, based on a sparse estimate of the inverse
covariance matrix.

\item\textit{Sequential procedures with false discovery rate control.}
It is also interesting to consider how the sequence of covariance test
$p$-values can be used to construct a sequential test with good power
properties, and a guaranteed bound on its false discovery rate.
A number of such approaches are proposed in \citet{fdrlasso}.

\item\textit{Proper $p$-values for forward stepwise.} Perhaps
surprisingly, a test analogous to the covariance test can be used in
forward stepwise regression, to construct valid $p$-values for
this greedy procedure. This work is in progress.
\item Other related problems include: estimation of
$\sigma^2$ when $p \geq n$, in the context of the covariance test;
power calculations and confidence
interval estimation; theory for linear models having strong
and weak signals (large and small true coefficients);
theory for the elastic net, generalized linear models, and the Cox
model.
\end{itemize}
As is clear from the above discussion, the covariance test
work has created
much excitement and activity among our close collaborators and
students. It is our hope that the current
paper will also broadly stimulate other researchers' interest in
this area, and that at some point, the joint efforts of the
community will yield a full set of inferential tools for the lasso
and other commonly used adaptive procedures.

\begin{appendix}\label{app11}
\section*{\texorpdfstring{Appendix}{Appendix}}

\subsection{\texorpdfstring{Proof of Lemma \protect\ref{lemcovtestknot}.}{Proof of Lemma 1}}\label{appcovtestknot}

By continuity of the lasso solution path at $\lambda_k$,
\[
P_Ay - \lambda_k \bigl(X_A^T
\bigr)^+ s_A = P_{A \cup\{j\}}y - \lambda_k
\bigl(X_{A\cup\{j\}}^T \bigr)^+ s_{A \cup\{j\}}
\]
and, therefore,
%
%
\begin{equation}
\label{eqcont} (P_{A\cup\{j\}}-P_A)y = \lambda_k
\bigl( \bigl(X_{A\cup\{j\}}^T \bigr)^+ s_{A \cup\{j\}} -
\bigl(X_A^T \bigr)^+ s_A \bigr).
\end{equation}
From this, we can obtain two identities: the first is
%
%
\begin{equation}
\label{eqfirst} y^T(P_{A\cup\{j\}}-P_A)y =
\lambda_k^2 \cdot\bigl\| \bigl(X_{A\cup\{j\}}^T
\bigr)^+ s_{A \cup\{j\}} - \bigl(X_A^T \bigr)^+
s_A \bigr\|_2^2,
\end{equation}
obtained by squaring both sides in \eqref{eqcont} (more precisely, taking
the inner product of the left-hand side with itself and the right-hand side
with itself), and noting that $(P_{A\cup\{j\}}-P_A)^2=P_{A\cup\{j\}}-P_A$;
the second is
%
%
\begin{eqnarray}\label{eqsecond}
&& y^T \bigl( \bigl(X_{A\cup\{j\}}^T
\bigr)^+ s_{A \cup\{j\}} - \bigl(X_A^T \bigr)^+
s_A \bigr)
\nonumber\\[-8pt]\\[-8pt]
&&\qquad = \lambda_k \cdot\bigl\| \bigl(X_{A\cup\{j\}}^T
\bigr)^+ s_{A \cup\{j\}} - \bigl(X_A^T \bigr)^+
s_A \bigr\|_2^2,\nonumber
\end{eqnarray}
obtained by taking the inner product of both sides in \eqref{eqcont}
with $y$,
and then using~\eqref{eqfirst}. Plugging \eqref{eqfirst} and \eqref
{eqsecond}
in for the first and second terms in \eqref{eqcovtestproj}, respectively,
then gives the result in \eqref{eqcovtestknot}.



\subsection{\texorpdfstring{Proof of Lemma \protect\ref{lemeventdual}.}{Proof of Lemma 4}}\label{appeventdual}

Note that
\begin{eqnarray*}
&& g(j_1,s_1) > g(j,s)
\\
&&\qquad  \Longleftrightarrow \quad
\frac{g(j_1,s_1)- ss_1R_{j,j_1} g(j_1,s_1)} {
1-ss_1R_{j,j_1}} > \frac{g(j,s)-ss_1R_{j,j_1} g(j_1,s_1)} {
1-ss_1R_{j,j_1}}
\\
&&\qquad \Longleftrightarrow\quad  g(j_1,s_1) > h^{(j_1,s_1)}(j,s),
\end{eqnarray*}
the first step following since $1-ss_1R_{j,j_1} > 0$,
and the second step following from the definition of
$h^{(j_1,s_1)}$. The intersection of the right-hand side above, over
all $(j,s)\neq(j_1,s_1)$, is equivalent to
\[
g(j_1,s_1)>g(j_1,-s_1),\qquad
g(j_1,s_1)>M(j_1,s_1).
\]
But the former inequality is the same as $g(j_1,s_1)>0$, because $g(j_1,s_1)$
and $g(j_1,-s_1)$ have opposite signs. Further, the inequality $g(j_1,s_1)>0$
is redundant, as $M(j_1,s_1) \geq0$. This gives the result.

\subsection{\texorpdfstring{Proof of Lemma \protect\ref{lemexpratio}.}{Proof of Lemma 5}}\label{appexpratio}

By l'H\^{o}pital's rule,
\[
\lim_{m\rightarrow\infty} \frac{\widebar{\Phi} (u(t,m)
)}{\widebar{\Phi}(m)} = \lim_{m\rightarrow\infty}
\frac{\phi(u(t,m) )}{\phi(m)} \cdot\frac{\partial u(t,m)}{\partial m},
\]
where $\phi$ is the standard normal density. First, note that
\[
\frac{\partial(t,m)}{\partial m} = \frac{1}{2} + \frac{m}{2\sqrt {m^2+4t}} \rightarrow1 \qquad\mbox{as } m\rightarrow\infty.
\]
Also, a straightforward calculation shows
\[
\log{\phi \bigl(u(t,m) \bigr)}-\log{\phi(m)} 
= \frac{m^2}{2} \bigl(1-\sqrt{1+4t/m^2} \bigr) -
\frac{t}{2} \rightarrow-t \qquad\mbox{as } m\rightarrow\infty,
\]
where in the last step we used the fact that
$(1-\sqrt{1+4t/m^2})/(2/m^2)\rightarrow-t/2$, again
by l'H\^{o}pital's rule. Therefore, $\phi(u(t,m))/\phi(m)
\rightarrow e^{-t}$, which completes the proof.

\subsection{\texorpdfstring{Proof of Lemma \protect\ref{lemexpconv}.}{Proof of Lemma 6}}\label{appexpconv}

Fix $\varepsilon>0$, and choose $m_0$ large enough that
\[
\biggl|\frac{\widebar{\Phi} (u(t,m/\sigma) )}{\widebar{\Phi
}(m/\sigma)} - e^{-t} \biggr| \leq\varepsilon\qquad\mbox{for all } m
\geq m_0.
\]
Starting from \eqref{eqt1bd},
\begin{eqnarray*}
\bigl|\mathbb{P}(T_1>t)-e^{-t} \bigr| &\leq& \sum
_{j_1,s_1} \int_0^\infty\biggl|
\frac{\widebar{\Phi} (u(t,m/\sigma) )} {
\widebar{\Phi}(m/\sigma)} - e^{-t} \biggr| \widebar{\Phi}(m/\sigma)
F_{M(j_1,s_1)}(dm)
\\
&\leq& \varepsilon\sum_{j_1,s_1} \int_{m_0}^\infty
\widebar{\Phi}(m/\sigma) F_{M(j_1,s_1)}(dm) + \sum
_{j_1,s_1} \int_0^{m_0}
F_{M(j_1,s_1)}(dm)
\\
&\leq& \varepsilon\sum_{j_1,s_1} \mathbb{P}
\bigl(g(j_1,s_1) > M(j_1,s_1)
\bigr) + \sum_{j_1,s_1} \mathbb{P} \bigl(M(j_1,s_1)
\leq m_0 \bigr).
\end{eqnarray*}
Above, the term multiplying $\varepsilon$ is equal to 1, and the
second term can be made arbitrarily small (say, less than
$\varepsilon$) by taking $p$ sufficiently large.

\subsection{\texorpdfstring{Proof of Theorem \protect\ref{thmgenxfirst}.}{Proof of Theorem 2}}\label{appgenxfirst}

We will show that for any fixed $m_0>0$ and $j_1,s_1$,
%
%
\begin{equation}
\label{eqmc} \mathbb{P} \bigl(M(j_1,s_1) \leq
m_0 \bigr) \leq c^{|S|},
\end{equation}
where $S \subseteq\{1,\ldots, p\}\setminus\{j_1\}$ is as in
the theorem for $j=j_1$, with size $|S| \geq d_p$,
and $c<1$ is a constant (not depending on $j_1$). This would imply that
\[
\sum_{j_1,s_1} \mathbb{P} \bigl(M(j_1,s_1)
\leq m_0 \bigr) \leq2p \cdot c^{d_p} \rightarrow0 \qquad
\mbox{as } p\rightarrow\infty,
\]
where we used the fact that $d_p/\log{p}\rightarrow\infty$
by \eqref{eqsgrow}. The above sum tending to zero
now implies the desired convergence result by Lemma~\ref{lemexpconv},
and hence it suffices to show \eqref{eqmc}. To this end, consider
\begin{eqnarray*}
M(j_1,s_1) &=& \max_{j\neq j_1, s}
\frac{sU_j -sR_{j,j_1} U_{j_1}}{1-ss_1 R_{j,j_1}}
\\
& \geq& \max_{j\neq j_1} \frac{|U_j - R_{j,j_1} U_{j_1}|}{1+|R_{j,j_1}|}
\\
& \geq& \max_{j \in S} \frac{|U_j-R_{j,j_1} U_{j_1}|}{2},
\end{eqnarray*}
where in both inequalities above we used the fact that $|R_{j,j_1}|<1$.
We can therefore use the bound
\[
\mathbb{P} \bigl(M(j_1,s_1) \leq m_0
\bigr) \leq\mathbb{P}\bigl(|V_j| \leq m_0, j \in S\bigr),
\]
where we define $V_j=(U_j - R_{j,j_1}U_{j_1})/2$ for $j\in S$.
Let $r=|S|$, and without a loss of generality, let
$S=\{1,\ldots, r\}$. We will show that
%
%
\begin{equation}
\label{eqvc} \mathbb{P}\bigl(|V_1| \leq m_0, \ldots,
|V_r| \leq m_0\bigr) \leq c^r\vadjust{\goodbreak}
\end{equation}
for $c=\Phi(2m_0/(\sigma\delta))-\Phi(-2m_0/(\sigma\delta))<1$,
by induction; this would complete the proof, as it would imply
\eqref{eqmc}. Before presenting this argument, we note
a few important facts. First, the condition in \eqref{eqcv} is really
a statement about conditional variances:
\begin{eqnarray}
\Var \bigl(U_i | U_\ell, \ell\in S\setminus\{i\} \bigr)
= \sigma^2 \cdot \bigl[1 - R_{i,S\setminus\{i\}} (R_{S\setminus\{i\},S\setminus\{i\}})^{-1}
R_{S\setminus\{i\},i} \bigr] \geq\sigma^2\delta^2\nonumber
\\
\eqntext{\mbox{for all } i\in S,}
\end{eqnarray}
where recall that $U_j=X_j^T y$, $j=1,\ldots, p$.
Second, since $U_1,\ldots, U_r$ are jointly normal, we have
%
%
\begin{eqnarray}\label{eqcondec}
\Var \bigl(U_i | U_\ell, \ell\in
S' \bigr) \geq\Var \bigl(U_i | U_\ell, \ell
\in S\setminus\{i\} \bigr) \geq\sigma^2\delta^2
\nonumber\\[-8pt]\\[-8pt]
\eqntext{\displaystyle\mbox{for any } S' \subseteq S\setminus\{i\} \mbox{ and } i \in S,}
\end{eqnarray}
which can be verified using the conditional variance formula (i.e.,
the law of total variance). Finally, the collection $V_1,\ldots, V_r$ is
independent of $U_{j_1}$,
because these random variables are jointly normal, and $\mathbb{E}[V_j
U_{j_1}]=0$ for all $j=1,\ldots, r$.

Now we give the inductive argument for \eqref{eqvc}.
For the base case, note that $V_1 \sim N(0,\tau_1^2)$, where its variance
is
\[
\tau_1^2 = \Var(V_1) =
\Var(V_1 | U_{j_1}) = \Var(U_1)/4 \geq
\sigma^2\delta^2/4,
\]
the second equality is due to the independence of $V_1$ and
$U_{j_1}$, and the last inequality comes from the fact that
conditioning can only decrease the variance, as stated above in
\eqref{eqcondec}. Hence,
\begin{eqnarray*}
\mathbb{P}\bigl(|V_1|\leq m_0\bigr) &=& \Phi(m_0/
\tau_1)-\Phi(-m_0/\tau_1)
\\
&\leq& \Phi \bigl(2m_0/(\sigma\delta) \bigr) - \Phi \bigl(-2m_0/(
\sigma \delta) \bigr) = c.
\end{eqnarray*}
Assume as the inductive hypothesis that $\mathbb{P}(|V_1| \leq m_0,
\ldots, |V_q|\leq m_0) \leq c^q$. Then
\begin{eqnarray*}
&& \mathbb{P}\bigl(|V_1| \leq m_0, \ldots, |V_{q+1}|
\leq m_0\bigr)
\\
&&\qquad = \mathbb{P} \bigl(|V_{q+1}| \leq m_0 |
|V_1| \leq m_0, \ldots, |V_q|\leq
m_0 \bigr) \cdot c^q.
\end{eqnarray*}
We have, using the independence of $V_1,\ldots, V_{q+1}$ and $U_{j_1}$,
\begin{eqnarray*}
V_{q+1} | V_1,\ldots, V_q &\stackrel{d} {=}&
V_{q+1} | V_1,\ldots, V_q, U_{j_1}
\\
&\stackrel{d} {=} & V_{q+1} | U_1,\ldots, U_q,
U_{j_1}
\\
&\stackrel{d} {=} & N \bigl(0,\tau_{q+1}^2 \bigr),
\end{eqnarray*}
where the variance is
\[
\tau_{q+1}^2 = \Var(V_{q+1} | U_1,
\ldots, U_q, U_{j_1}) = \Var(U_{q+1} |
U_1,\ldots, U_q) / 4 \geq\sigma^2
\delta^2/4
\]
and here we again used the fact that conditioning further can only reduce
the variance, as in \eqref{eqcondec}. Therefore,
\[
\mathbb{P}\bigl(|V_{q+1}| \leq m_0 | V_1,\ldots,
V_q\bigr) \leq\Phi \bigl(2m_0/(\sigma\delta) \bigr) - \Phi
\bigl(-2m_0/(\sigma\delta) \bigr) = c
\]
and so
\[
\mathbb{P}\bigl(|V_1| \leq m_0, \ldots, |V_{q+1}|
\leq m_0\bigr) \leq c \cdot c^q = c^{q+1},
\]
completing the inductive step.

\subsection{\texorpdfstring{Proof of Lemma \protect\ref{lemeventdual2}.}{Proof of Lemma 7}}\label{appeventdual2}
Notice that
\begin{eqnarray*}
&& g(j_k,s_k) > g(j,s)
\\
&&\qquad \Longleftrightarrow\quad g(j_k,s_k)
( 1 - \Sigma_{j,j'}/\Sigma_{jj}) > g(j,s) - (
\Sigma_{j,j'}/\Sigma_{jj}) g(j_k,s_k).
\end{eqnarray*}
We now handle division by $1-\Sigma_{j,j'}/\Sigma_{jj}$ in three
cases:
\begin{itemize}
\item if $1 - \Sigma_{j,j'}/\Sigma_{jj} > 0$, then
\[
g(j_k,s_k) > g(j,s)\quad\Longleftrightarrow\quad g(j_k,s_k)
> \frac{g(j,s) - (\Sigma_{j,j'}/\Sigma_{jj}) g(j_k,s_k)} {
1 - \Sigma_{j,j'}/\Sigma_{jj}};
\]
\item if $1 -\Sigma_{j,j'}/\Sigma_{jj} < 0$, then
\[
g(j_k,s_k) > g(j,s) \quad\Longleftrightarrow\quad g(j_k,s_k)
< \frac{g(j,s) - (\Sigma_{j,j'}/\Sigma_{jj}) g(j_k,s_k)} {
1 -\Sigma_{j,j'}/\Sigma_{jj}};
\]
\item if $1 -\Sigma_{j,j'}/\Sigma_{jj} = 0$, then
\[
g(j_k,s_k) > g(j,s) \quad\Longleftrightarrow\quad 0 > g(j,s) - (
\Sigma_{j,j'}/\Sigma_{jj}) g(j_k,s_k).
\]
\end{itemize}
Using this breakdown, we see that the statement
$g(j_k,s_k)>g(j,s)$ for all $(j,s)\neq(j_k,s_k)$ is then
equivalent to
\begin{eqnarray*}
g(j_k,s_k) &>& g(j_k,-s_k),\qquad
g(j_k,s_k) > M^+(j_k,s_k),
\\
g(j_k,s_k) &<& M^-(j_k,s_k),\qquad
0 > M^0(j_k,s_k).
\end{eqnarray*}
Noting that $g(j_k,s_k)$ and $g(j_k,-s_k)$ must have opposite
signs, the above is equivalent to
\begin{eqnarray*}
g(j_k,s_k) &>& 0,\qquad
g(j_k,s_k) > M^+(j_k,s_k),
\\
g(j_k,s_k) &<& M^-(j_k,s_k),\qquad
0 > M^0(j_k,s_k),
\end{eqnarray*}
which gives the result in the lemma.

\subsection{\texorpdfstring{Proof of Lemma \protect\ref{lemexpconv2}.}{Proof of Lemma 8}}\label{appexpconv2}

Define $\sigma_k = \sigma/\sqrt{C(j_k,s_k)}$ and $u(a,b) =
(b+\sqrt{b^2+4a})/2$. Exactly as before
(dropping for simplicity the notational dependence of $g,M^+$
on $j_k,s_k$),
\[
g \bigl(g-M^+ \bigr)/\sigma_k^2 > t,\qquad
g > M^+ \quad\Longleftrightarrow\quad g /\sigma_k > u \bigl(t,M^+/\sigma_k
\bigr).\vadjust{\goodbreak}
\]
Therefore, we can rewrite \eqref{eqttksum} as
\begin{eqnarray*}
\mathbb{P}(\widetilde{T}_k > t) &=& \sum_{j_k,s_k}
\mathbb{P} \bigl( g(j_k,s_k)/\sigma_k > u
\bigl(t, M^+(j_k,s_k)/\sigma_k \bigr),
\\[-1pt]
&&\hspace*{29pt} g(j_k,s_k) < M^-(j_k,s_k),
0 > M^0(j_k,s_k) \bigr).
\end{eqnarray*}
Note that we have dropped the inequality $g(j_k,s_k)>0$ from
each term, as it is implied by the first inequality
$g(j_k,s_k)/\sigma_k>u(t,M^+(j_k,s_k)/\sigma_k) \geq0$.
We can upper bound the right-hand side above by replacing
$g(j_k,s_k) < M^-(j_k,s_k)$ with
\[
g(j_k,s_k) < M^-(j_k,s_k) + u
\bigl(t\sigma_k^2,M^+(j_k,s_k)
\bigr)-M^+(j_k,s_k),
\]
because $u(a,b) \geq b$ for all $a \geq0$ and $b$. Furthermore,
Lemma~\ref{lemg2cov2}  (Appen\-dix~\ref{appg2cov2}) shows that
indeed $\sigma_k^2 = \sigma^2/C(j_k,s_k) = \Var(g(j_k,s_k))$ for
fixed $j_k,s_k$, and hence $g(j_k,s_k)/\sigma_k$ is standard normal for
fixed $j_k,s_k$. Therefore,
%
%
\begin{eqnarray}
\label{eqttksum2} \mathbb{P}(\widetilde{T}_k > t) &\leq& \sum
_{j_k,s_k} \int \bigl[\Phi \bigl(m^-/\sigma_k + u
\bigl(t,m^+/\sigma_k \bigr) - m^+/\sigma_k \bigr) -\Phi
\bigl(u \bigl(t,m^+/\sigma_k \bigr) \bigr) \bigr]\hspace*{-20pt}
\nonumber\\[-9pt]\\[-9pt]
&&\hspace*{26pt}{}\times G_{j_k,s_k} \bigl(dm^+,dm^-,dm^0 \bigr),\nonumber
\end{eqnarray}
where
\begin{eqnarray*}
&& G_{j_k,s_k} \bigl(dm^+,dm^-,dm^0 \bigr)
\\[-1pt]
&&\qquad = 1 \bigl\{m^+ < m^-,
m^0 < 0 \bigr\} \cdot F_{M^+(j_k,s_k), M^-(j_k,s_k), M^0(j_k,s_k)} \bigl(dm^+,dm^-,dm^0
\bigr)
\end{eqnarray*}
with $F_{M^+(j_k,s_k), M^-(j_k,s_k), M^0(j_k,s_k)}$
the joint distribution of $M^+(j_k,s_k)$,\break
$M^-(j_k,s_k)$, $M^0(j_k,s_k)$, and we used the fact that $g$ is
independent of $M^+,\break M^-,M^0$ for fixed $j_k,s_k$. From
\eqref{eqttksum2},
%
%
\begin{eqnarray}
\label{eqttksum3}
&& \mathbb{P}(\widetilde{T}_k > t) - e^{-t}\nonumber
\\[-1pt]
&&\qquad \leq \sum_{j_k,s_k} \int \biggl( \frac{\Phi(m^-/\sigma_k +
u(t,m^+/\sigma_k)-m^+/\sigma_k )-\Phi(u(t,m^+/\sigma
_k) )} {
\Phi(m^-/\sigma_k)-\Phi(m^+/\sigma_k)} -
e^{-t} \biggr)\hspace*{-8pt}
\\[-1pt]
&&\hspace*{59pt}{}\times
\bigl[\Phi \bigl(m^-/\sigma_k \bigr)-\Phi \bigl(m^+/
\sigma_k \bigr) \bigr] \cdot G_{j_k,s_k} \bigl(dm^+,dm^-,dm^0
\bigr),\nonumber
\end{eqnarray}
where we here used the fact that
\begin{eqnarray*}
&& \sum_{j_k,s_k} \int \bigl[\Phi \bigl(m^-/
\sigma_k \bigr)-\Phi \bigl(m^+/\sigma_k \bigr) \bigr]
G_{j_k,s_k} \bigl(dm^+,dm^-,dm^0 \bigr)
\\[-1pt]
&&\qquad = \sum_{j_k,s_k} \mathbb{P} \bigl(
g(j_k,s_k) > M^+(j_k,s_k),
g(j_k,s_k) < M^-(j_k,s_k), 0 >
M^0(j_k,s_k) \bigr)
\\[-1pt]
&&\qquad \geq\sum_{j_k,s_k} \mathbb{P} \bigl(
g(j_k,s_k)>0, g(j_k,s_k) >
M^+(j_k,s_k),
\\[-1pt]
&&\hspace*{63pt} g(j_k,s_k) <M^-(j_k,s_k), 0 > M^0(j_k,s_k)\bigr)
\\[-1pt]
&&\qquad = 1,\vadjust{\goodbreak}
\end{eqnarray*}
the last equality following by Lemma~\ref{lemeventdual2}
(i.e., each term in the last sum is exactly the probability of
$j_k,s_k$ maximizing $g$). We show in
Lemma~\ref{lemexpratio2} (Appendix~\ref{appexpratio2}) that
\[
\lim_{m^+\rightarrow\infty} \frac{\Phi(m^-+u(t,m^+)-m^+ )-\Phi
(u(t,m^+) )} {
\Phi(m^-)-\Phi(m^+)} \leq e^{-t},
\]
provided that $m^- > m^+$.
Hence, fix $\varepsilon>0$, and choose $m_0$ sufficiently large, so that
for each $k$,
\begin{eqnarray}
\frac{\Phi(m^-/\sigma_k + u(t,m^+/\sigma_k) - m^+/\sigma
_k )-
\Phi(u(t,m^+/\sigma_k) )}{\Phi(m^-/\sigma_k)-\Phi
(m^+/\sigma_k)} - e^{-t} \leq\varepsilon\nonumber
\\
\eqntext{\displaystyle\mbox{for all } m^-/\sigma_k > m^+/\sigma_k \geq
m_0.}
\end{eqnarray}
%
Working from \eqref{eqttksum3},
\begin{eqnarray*}
&& \mathbb{P}(\widetilde{T}_k > t) - e^{-t}
\\
&&\qquad \leq\varepsilon
\sum_{j_k,s_k} \int_{m^+/\sigma_k \geq m_0} \bigl[\Phi
\bigl(m^-/\sigma_k \bigr)-\Phi \bigl(m^+/\sigma_k \bigr)
\bigr] G_{j_k,s_k} \bigl(dm^+,dm^-,dm^0 \bigr)
\\
&&\quad\qquad{} + \sum_{j_k,s_k} \int_{m^+/\sigma_k \leq m_0}
G_{j_k,s_k} \bigl(dm^+,dm^-,dm^0 \bigr).
\end{eqnarray*}
Note that the first term on the right-hand side above is
$\leq\varepsilon$, and the second term is\vspace*{2pt} bounded by
$\sum_{j_k,s_k} \mathbb{P}(M^+(j_k,s_k) \leq m_0\sigma_k)$,
which by assumption can be made arbitrarily small (smaller than, say,
$\varepsilon$) by taking $p$ large enough.

\subsection{\texorpdfstring{Proof of Lemma \protect\ref{lemexpconv3}.}{Proof of Lemma 9}}\label{appexpconv3}

For now, we reintroduce the notational dependence of the process
$g$ on $A,s_A$, as this will be important.
We show in Lemma~\ref{lemg3} (Appendix~\ref{appg3}) that for any
fixed $j_k,s_k,j,s$,
\[
\frac{ g^{(A,s_A)}(j,s) - (\Sigma_{j_k,j}/\Sigma_{j_k,j_k})
g^{(A,s_A)}(j_k,s_k)} {1 - \Sigma_{j_k,j}/\Sigma_{j_k,j_k}} = g^{(A\cup
\{j_k\},s_{A\cup\{j_k\}})}(j,s),
\]
where $\Sigma_{j_k,j} = \mathbb{E}[g^{(A,s_A)}(j_k,s_k),g^{(A,s_A)}(j,s)]$,
as given in \eqref{eqg2cov}, and as usual, $s_{A\cup\{j_k\}}$ denotes
the concatenation of $s_A$ and $s_k$. According to its definition in
\eqref{eqmp}, therefore,
\[
M^+(j_k,s_k) = \max_{(j,s) \in S^+(j_k,s_k)}
g^{(A\cup\{j_k\},s_{A\cup\{j_k\}})}(j,s)
\]
and hence on the event $E(j_k,s_k)$, since we have
$g^{(A,s_A)}(j_k,s_k) > M^+(j_k,s_k)$,
\begin{eqnarray*}
&& M^+(j_k,s_k)
\\
&&\quad = \max_{(j,s) \in S^+(j_k,s_k)}
g^{(A\cup\{j_k\},s_{A\cup\{j_k\}})}(j,s)
\\
&&\qquad\hspace*{54pt} {}\times 1 \bigl\{ g^{(A\cup\{
j_k\},s_{A\cup\{j_k\}})}(j,s) < g^{(A,s_A)}(j_k,s_k)
\bigr\}
\\
&&\quad \leq \max_{j \notin A \cup\{j_k\}, s} g^{(A\cup\{j_k\},s_{A\cup\{j_k\}
})}(j,s) \cdot1 \bigl\{
g^{(A\cup\{j_k\},s_{A\cup\{j_k\}})}(j,s) < g^{(A,s_A)}(j_k,s_k) \bigr\}
\\
&&\quad =  M(j_k,s_k).
\end{eqnarray*}
This means that (now we return to writing $g^{(A,s_A)}$ as $g$, for
brevity)
\begin{eqnarray*}
&& \sum_{j_k,s_k} \mathbb{P} \bigl( \bigl
\{C(j_k,s_k) \cdot g(j_k,s_k)
\bigl(g(j_k,s_k)-M(j_k,s_k)
\bigr)/\sigma^2 >t \bigr\} \cap E(j_k,s_k)
\bigr)
\\
&&\quad \leq
\sum_{j_k,s_k} \mathbb{P} \bigl( \bigl
\{C(j_k,s_k) \cdot g(j_k,s_k)
\bigl(g(j_k,s_k)-M^+(j_k,s_k)
\bigr)/\sigma^2 >t \bigr\} \cap E(j_k,s_k)
\bigr)
\end{eqnarray*}
and so $\lim_{p\rightarrow\infty} \mathbb{P}(T_k > t) \leq
\lim_{p\rightarrow\infty} \mathbb{P}(\widetilde{T}_k > t) \leq
e^{-t}$, the desired
conclusion.

\subsection{\texorpdfstring{Proof of Theorem \protect\ref{thmgenxgen}.}{Proof of Theorem 3}}\label{appgenxgen}

Since we are assuming that $\mathbb{P}(B) \rightarrow1$, we know that
$\mathbb{P}(T_k>t | B) - \mathbb{P}(T_k > t) \rightarrow0$, so we
only need to
consider the marginal limiting distribution of $T_k$. We write
$A=A_0$ and $s_A=s_{A_0}$. The general idea here is similar to that
used in the proof of Theorem~\ref{thmgenxfirst}. Fixing $m_0$ and
$j_k,s_k$, we will show that
%
%
\begin{equation}
\label{eqmpc} \mathbb{P} \bigl(M^+ (j_k,s_k) \leq
m_0 \sigma_k \bigr) \leq c^{|S|},
\end{equation}
where $S \subseteq\{1,\ldots, p\} \setminus(A \cup\{j_k\})$ is as
in the statement of the theorem for $j=j_k$, with size
$|S| \geq d_p$, and $c<1$ is a constant (not depending on $j_k$).
Also, as in the proof of Lemma~\ref{lemexpconv2}, we
abbreviated $\sigma_k = \sigma/\sqrt{C(j_k,s_k)}$.
This bound would imply that
\[
\sum_{j_k,s_k} \mathbb{P} \bigl(M^+(j_k,s_k)
\leq m_0\sigma_k \bigr) \leq2p \cdot c^{d_p}
\rightarrow0 \qquad\mbox{as } p\rightarrow\infty,
\]
since $d_p/\log{p}\rightarrow0$. The above sum converging to zero is
precisely the condition required by Lemma~\ref{lemexpconv3}, which
then gives the desired (conservative) exponential limit for
$T_k$. Hence, it is suffices to show \eqref{eqmpc}. For this, we start
by recalling the definition of $M^+$ in \eqref{eqmp}:
\begin{eqnarray}
M^+(j_k,s_k) = \max_{(j,s) \in S^+(j_k,s_k)}
\frac{ g(j,s) - (\Sigma_{j_k,j}/\Sigma_{j_k,j_k}) g(j_k,s_k)} {
1 - \Sigma_{j_k,j}/\Sigma_{j_k,j_k}}\nonumber
\\
\eqntext{\displaystyle\mbox{where } S^+(j_k,s_k) = \biggl\{ (j,s)\dvtx j
\notin A \cup\{j_k\}, \frac{\Sigma_{j_k,j}}{\Sigma_{j_k,j_k}} < 1 \biggr\}.}
\end{eqnarray}
Here, we write $\Sigma_{j_k,j}=\mathbb{E}[g(j_k,s_k) g(j,s)]$;
note that $\Sigma_{j_k,j_k}=\sigma^2_k$
(as shown in Lemma~\ref{lemg2cov2}). First, we show that the
conditions of the theorem actually imply that
$S^+(j_k,s_k) \supseteq S \times\{-1,1\}^{|S|}$. This is true\vadjust{\goodbreak} because
for $j \in S$ and any $s\in\{-1,1\}$, we have $|R_{j,j_k}|/R_{j_k,j_k}
< \eta/(2-\eta)$ by \eqref{eqcovvar}, and
\begin{eqnarray*}
|R_{j,j_k}|/R_{j_k,j_k} < \eta/ (2-\eta) &\quad\Longrightarrow\quad&\biggl|
\frac{R_{j,j_k}}{R_{j_k,j_k}} \cdot\frac{s_k - X_{j_k}^T (X_A^T)^+
s_A} {
s - X_j^T (X_A^T)^+ s_A} \biggr| < 1
\\
&\quad\Longrightarrow\quad &\Sigma_{j_k,j}/\Sigma_{j_k,j_k} < 1.
\end{eqnarray*}
The first implication uses the assumption \eqref{eqirrep}, as
$|s_k-X_{j_k}^T (X_A^T)^+ s_A| \leq1 + \|(X_A)^+X_{j_k}\|_1 \leq
2-\eta$
and $|s-X_j^T (X_A^T)^+ s_A| \geq1 - \|(X_A)^+X_{j_k}\|_1 \geq\eta$,
and\vspace*{2pt} the second simply follows from the definition of
$\Sigma_{j_k,j}$ and $\Sigma_{j_k,j_k}$. Therefore,
\[
M^+(j_k,s_k) \geq\max_{j \in S, s}
\frac{ g(j,s) - (\Sigma_{j_k,j}/\Sigma_{j_k,j_k}) g(j_k,s_k)} {
1 - \Sigma_{j_k,j}/\Sigma_{j_k,j_k}}.
\]
Let\vspace*{2pt} $U_j=X_j^T(I-P_A)y$ and
$\theta_{j_k,j} = R_{j_k,j}/R_{j_k,j_k}$ for $j\in S$. By the
arguments given in the proof of Lemma~\ref{lemg3}, we can rewrite the
right-hand side above, yielding
\begin{eqnarray*}
M^+(j_k,s_k) &\geq&\max_{j\in S, s}
\frac{U_j - \theta_{j_k,j} U_{j_k}}{s - X_j^T
(X_{A\cup\{j_k\}}^T)^+ s_{A\cup\{j_k\}}}
\\
&=& \max_{j\in S, s} \frac{s(U_j - \theta_{j_k,j} U_{j_k})}{1 - sX_j^T
(X_{A\cup\{j_k\}}^T)^+ s_{A\cup\{j_k\}}}
\\
&\geq&\max_{j\in S} \frac{|U_j - \theta_{j_k,j} U_{j_k}|}{1 +
|X_j^T (X_{A\cup\{j_k\}}^T)^+ s_{A\cup\{j_k\}}|}
\\
&\geq&\max_{j\in S} \frac{|U_j - \theta_{j_k,j} U_{j_k}|}{2},
\end{eqnarray*}
where the last two inequalities above follow as
$|X_j^T(X_{A\cup\{j_k\}})^+s_{A\cup\{j_k\}}| \,{<}\, 1$~for all $j\in S$, which
itself follows from the assumption that $\|(X_{A\cup\{j_k\}})^+
X_j\|_\infty< 1$ for all $j\in S$, in \eqref{eqirrep2}. Hence,
\[
\mathbb{P} \bigl(M^+(j_k,s_k) \leq m_0
\sigma_k \bigr) \leq\mathbb{P}\bigl(|V_j| \leq m_0
\sigma_k, j\in S\bigr),
\]
where $V_j = (U_j-\theta_{j_k,j}U_{j_k})/2$. Writing without a loss of
generality $r=|S|$ and $S=\{1,\ldots, r\}$, it now remains
to show that
%
%
\begin{equation}
\label{eqvc2} \mathbb{P}\bigl(|V_1| \leq m_0
\sigma_k, \ldots, |V_r| \leq m_0
\sigma_k\bigr) \leq c^r.
\end{equation}
Similar to the arguments in the proof of Theorem~\ref{thmgenxfirst},
we will show \eqref{eqvc2} by induction, for the constant
$c=\Phi(2m_0\sqrt{C}/(\delta\eta))-
\Phi(-2m_0\sqrt{C}/(\delta\eta)) < 1$. Before this, it is helpful to
discuss three important facts. First, we note that \eqref{eqcv2} is
actually a lower bound on the ratio of conditional to unconditional
variances:
\begin{eqnarray*}
&& \Var \bigl(U_i | U_\ell, \ell\in S\setminus\{i\} \bigr)
/ \Var(U_i)
\\
&&\qquad = \bigl[R_{ii} - R_{i,S\setminus\{i\}}
(R_{S\setminus\{i\},S\setminus\{i\}})^{-1} R_{S\setminus\{i\},i} \bigr] / R_{ii}
\geq\delta^2\qquad \mbox{for all } i\in S.
\end{eqnarray*}
Second, conditioning on a smaller set of variables can only increase
the conditional variance:
\begin{eqnarray}
\Var \bigl(U_i | U_\ell, \ell\in S' \bigr)
\geq\Var \bigl(U_i | U_\ell, \ell\in S\setminus\{i\}
\bigr) \geq\delta^2 \sigma^2 R_{ii}\nonumber
\\
\eqntext{\displaystyle\mbox{for any } S' \subseteq S\setminus\{i\} \mbox{ and }i \in S,}
\end{eqnarray}
which holds as $U_1,\ldots, U_r$ are jointly normal. Third, and
lastly, the collection $V_1,\ldots, V_r$ is independent of $U_{j_k}$,
since these variables are all jointly normal, and it is easily
verified that $\mathbb{E}[V_j U_{j_k}] = 0$ for each $j=1,\ldots, r$.

We\vspace*{1pt} give the inductive argument for \eqref{eqvc2}. For the base
case, we have $V_1 \sim N(0,\tau_1^2)$, where
\[
\tau_1^2 = \Var(V_1) =
\Var(V_1|U_{j_k}) = \Var(U_1)/4 \geq
\delta^2 \sigma^2 R_{11} / 4.
\]
Above, in the second equality, we used that $V_1$ and $U_{j_k}$ are
independent, and in the last inequality, that conditioning on
fewer variables (here, none) only increases the variance. This means
that
\[
\mathbb{P}\bigl(|V_1| \leq m_0 \sigma_k\bigr) \leq
\mathbb{P} \bigl(|Z| \leq2m_0 \sigma_k / (\delta\sigma
\sqrt{R_{11}} ) \bigr) \leq\mathbb{P} \bigl(|Z| \leq2m_0
\sqrt{C} / (\delta\eta) \bigr)= c,
\]
where $Z$ is standard normal; note that in the
last inequality above, we applied the upper bound
\[
\frac{\sigma^2_k}{\sigma^2 R_{11}} = \frac{\Sigma_{j_k,j_k}}{\sigma
^2R_{11}} = \frac{R_{j_k,j_k}}{R_{11}} \cdot
\frac{1}{ [s_k-X_{j_k}^T (X_A^T)^+ s_{A} ]^2} \leq\frac{C}{\eta^2}.
\]
Now, for the inductive hypothesis, assume that
$\mathbb{P}(|V_1| \leq m_0\sigma_k, \ldots, |V_q| \leq m_0\sigma_k)
\leq c^q$.
Consider
\begin{eqnarray*}
&& \mathbb{P}\bigl(|V_1| \leq m_0\sigma_k, \ldots,
|V_{q+1}|\leq m_0\sigma_k\bigr)
\\
&&\qquad = \mathbb{P}
\bigl(|V_{q+1}| \leq m_0\sigma_k | |V_1|
\leq m_0\sigma_k, \ldots, |V_q|\leq
m_0\sigma_k \bigr) \cdot c^q.
\end{eqnarray*}
Using the independence of $V_1,\ldots, V_{q+1}$ and $U_{j_k}$,
\begin{eqnarray*}
V_{q+1} | V_1,\ldots, V_q &\stackrel{d} {=}&
V_{q+1} | V_1,\ldots, V_q, U_{j_k}
\\
&\stackrel{d} {=} &V_{q+1} | U_1,\ldots, U_q,
U_{j_k}
\\
&\stackrel{d} {=} &N \bigl(0,\tau_{q+1}^2 \bigr).
\end{eqnarray*}
The variance $\tau_{q+1}^2$ is
\begin{eqnarray*}
\tau_{q+1}^2 &=& \Var(V_{q+1} | U_1,
\ldots, U_q, U_{j_k})
\\
&=& \Var(U_{q+1} |
U_1,\ldots, U_q) / 4
\\
&\geq&\delta^2
\sigma^2 R_{q+1,q+1} / 4,
\end{eqnarray*}
where we again used the fact that conditioning on a smaller set of
variables only makes the variance larger. Finally,
\begin{eqnarray*}
\mathbb{P}\bigl(|V_{q+1}| \leq m_0\sigma_k |
V_1,\ldots, V_q\bigr) &\leq& \mathbb{P} \bigl(|Z|
\leq2m_0 \sigma_k / (\delta\sigma\sqrt{R_{q+1,q+1}}
) \bigr)
\\
&\leq& \mathbb{P} \bigl(|Z| \leq2m_0 \sqrt{C} / (\delta\eta)
\bigr)
\\
&=& c,
\end{eqnarray*}
where we used $\sigma_k^2/(\sigma^2 R_{q+1,q+1}) \leq C/\eta^2$ as
above, and so
\[
\mathbb{P}\bigl(|V_1| \leq m_0\sigma_k, \ldots,
|V_{q+1}|\leq m_0\sigma_k\bigr) \leq c \cdot
c^q = c^{q+1}.
\]
This completes the inductive proof.

\subsection{\texorpdfstring{Statement and proof of Lemma \protect\ref{lemg2cov2}.}{Statement and proof of Lemma 10}}\label{appg2cov2}

%
%
\begin{lemma}
\label{lemg2cov2}
For any fixed $A,s_A$, and any $j\notin A$, $s \in\{-1,1\}$, we have
\[
\Var \bigl(g(j,s) \bigr) = \frac{X_j^T (I-P_A) X_j^T\sigma^2}{[s -
X_j^T (X_A^T)^+ s_A]^2} = \frac{\sigma^2}{\|(X_{A\cup\{ j\}}^T)^+ s_{A
\cup\{j\}} -
(X_A^T)^+ s_A \|_2^2},
\]
where $s_{A\cup\{j\}}$ denotes the concatenation of $s_A$ and $s$.
\end{lemma}

\begin{pf}
We will show that
%
%
\begin{equation}
\label{eqinv} \frac{[s - X_j^T (X_A^T)^+ s_A]^2}{X_j^T (I-P_A) X_j^T} = \bigl\| \bigl(X_{A\cup\{ j\}}^T
\bigr)^+ s_{A \cup\{j\}} - \bigl(X_A^T \bigr)^+
s_A\bigr\|_2^2.
\end{equation}
The right-hand side above, after a straightforward calculation, is shown
to be equal to
%
%
\begin{equation}
\label{eqrhs2} s_{A\cup\{j\}}^T \bigl(X_{A\cup\{j\}}^T
X_{A\cup\{j\}} \bigr)^{-1} s_{A\cup\{
j\}} - s_A
\bigl(X_A^T X_A \bigr)^{-1}
s_A.
\end{equation}
Now let $z = (X_{A\cup\{j\}}^T X_{A\cup\{j\}})^{-1} s_{A\cup\{j\}}$.
In block form,
%
%
\begin{equation}
\label{eqblock} \lleft[ \matrix{ X_A^T
X_A & X_A^T X_j
\vspace*{5pt}\cr
X_j^T X_A & X_j^T
X_j} %
\rright] \lleft[ %
\matrix{z_1
\cr
z_2 } \rright] = \lleft[ %
\matrix{ s_A
\cr
s} %
\rright].
\end{equation}
Solving for $z_1$ in the first row yields
\[
z_1 = \bigl(X_A^T X_A
\bigr)^{-1} s_A - (X_A)^+ X_j
z_2
\]
and, therefore, \eqref{eqrhs2} is equal to
%
%
\begin{equation}
\label{eqrhs3} s_A^T z_1 +
sz_2 - s_A^T \bigl(X_A^T
X_A \bigr)^{-1} s_A = \bigl[s-s_A^T
(X_A)^+ X_j \bigr]z_2.
\end{equation}
Solving for $z_2$ in the second row of \eqref{eqblock} gives
\[
z_2 = \frac{s-s_A^T (X_A)^+ X_j}{X_j^T(I-P_A)X_j}.
\]
Plugging this value into \eqref{eqrhs3} produces the left-hand side
in \eqref{eqinv}, completing the proof.
\end{pf}

\subsection{\texorpdfstring{Statement and proof of Lemma \protect\ref{lemexpratio2}.}{Statement and proof of Lemma 11}}\label{appexpratio2}

%
%
\begin{lemma}
\label{lemexpratio2}
If $v=v(m)$ satisfies $v>m$, then for any $t \geq0$,
\[
\lim_{m\rightarrow\infty} \frac{\Phi(v+u(t,m)-m )-\Phi(u(t,m) )} {
\Phi(v)-\Phi(m)} \leq e^{-t}.
\]
\end{lemma}

\begin{pf}
First note, using a Taylor series expansion of $\sqrt{1+4t/m^3}$, that
for sufficiently large $m$,
%
%
\begin{equation}
\label{equbound} u(t,m) \geq m + \frac{t}{m} -\frac{t^2}{m^3}.
\end{equation}
Also, a simple calculation shows that $\partial(u(t,m)-m)/\partial m
\leq0$ for all $m$, so that
%
%
\begin{equation}
\label{equnoninc} u(t,w) - w \leq u(t,m) - m \qquad\mbox{for all } w \geq m.
\end{equation}
Now consider
\begin{eqnarray*}
\nonumber
\Phi \bigl(v+u(t,m)-m \bigr)-\Phi \bigl(u(t,m) \bigr) &=& \int
_{u(t,m)}^{v+u(t,m)-m} \frac{e^{-z^2/2}}{\sqrt{2\pi}} \,dz
\\
&=& \int_m^v \frac{e^{-(w+u(t,m)-m)^2/2}}{\sqrt{2\pi}} \,dw
\\
&\leq&\int_m^v \frac{e^{-u(t,m)^2/2}}{\sqrt{2\pi}} \,dw
\\
&\leq&\int_m^v \frac{e^{-(w+t/m-t^2/m^3)^2/2}}{\sqrt{2\pi}} \,dw,
\end{eqnarray*}
where the first inequality follows from \eqref{equnoninc}, and the
second from \eqref{equbound} (assuming $m$ is large enough).
Continuing from the last upper bound,
\[
\int_m^v \frac{e^{-(w+t/m-t^2/m^3)^2/2}}{\sqrt{2\pi}} \,dw =
e^{-t} \int_m^v \frac{e^{-w^2/2}}{\sqrt{2\pi}}
f(w,t) \,dw,
\]
where
\[
f(w,t) = \exp \biggl(\frac{t^2}{2w^2}+\frac{t^3}{w^4}-\frac
{t^4}{2w^6}
\biggr).
\]
Therefore, we have
%
%
\begin{eqnarray}\label{eqalmost}
&& \frac{\Phi(v+u(t,m)-m )-\Phi(u(t,m) )} {
\Phi(v)-\Phi(m)} - e^{-t}
\nonumber\\[-8pt]\\[-8pt]
&&\qquad \leq \biggl(
\frac{\int_m^v (e^{-w^2/2}/\sqrt{2\pi}) f(w,t) \,dw} {
\int_m^v (e^{-w^2/2}/\sqrt{2\pi}) \,dw} - 1 \biggr) \cdot e^{-t}.\nonumber
\end{eqnarray}
It is clear that $f(w,t) \rightarrow1$ as
$w\rightarrow\infty$. Fixing $\varepsilon$, choose $m_0$ large enough so
that for all $w \geq m_0$, we have $|f(w,t)-1| \leq\varepsilon$. Then
the term multiplying $e^{-t}$ on the right-hand side in
\eqref{eqalmost}, for $m \geq m_0$, is
\[
\frac{\int_m^v ({e^{-w^2/2}}/{\sqrt{2\pi}}) f(w,t) \,dw} {
\int_m^v (e^{-w^2/2}/\sqrt{2\pi}) \,dw} - 1 \leq\frac{\int_m^v
({e^{-w^2/2}}/{\sqrt{2\pi}})
| f(w,t) -1 | \,dw} {
\int_m^v (e^{-w^2/2}/\sqrt{2\pi}) \,dw} \leq\varepsilon,
\]
which shows that the right-hand side in \eqref{eqalmost} is
$\leq\varepsilon\cdot e^{-t} \leq\varepsilon$, and completes the
proof.
\end{pf}



\subsection{\texorpdfstring{Statement and proof of Lemma \protect\ref{lemg3}.}{Statement and proof of Lemma 12}}\label{appg3}

%
%
\begin{lemma}
\label{lemg3}
For any fixed $j_k,s_k,j,s$ (and fixed $A,s_A)$, we have
%
%
\begin{equation}
\label{eqg3} \qquad \frac{ g^{(A,s_A)}(j,s) - (\Sigma_{j_k,j}/\Sigma_{j_k,j_k})
g^{(A,s_A)}(j_k,s_k)} {1 - \Sigma_{j_k,j}/\Sigma_{j_k,j_k}} = g^{(A\cup
\{j_k\},s_{A\cup\{j_k\}})}(j,s),
\end{equation}
where $\Sigma_{j_k,j}$ denotes the covariance between
$g^{(A,s_A)}(j_k,s_k)$ and $g^{(A,s_A)}(j,s)$,
\[
\Sigma_{j_k,j} = \frac{X_{j_k}^T (I-P_A) X_j \sigma^2} {
[s_k - s_A^T (X_A)^+ X_{j_k}]
[s - s_A^T (X_A)^+ X_j]}.
\]
\end{lemma}

\begin{pf}
Simple manipulations of the left-hand side in \eqref{eqg3} yield the
expression
%
%
\begin{equation}
\label{eqlhs} \frac{X_j^T(I-P_A)y - \theta_{j_k,j} \cdot X_{j_k}^T
(I-P_A) y} {
s-s_A^T (X_A)^+ X_j - \theta_{j_k,j} \cdot[s_k-s_A^T (X_A)^+ X_{j_k}]},
\end{equation}
where $\theta_{j_k,j} = X_{j_k}^T (I-P_A) X_j / (X_{j_k}^T (I-P_A)
X_{j_k})$. Now it remains to show that~\eqref{eqlhs} is equal to
%
%
\begin{equation}
\label{eqrhs} \frac{X_j^T (I-P_{A\cup\{j_k\}}) y}{s - s_{A\cup\{j_k\}}^T
(X_{A\cup\{j_k\}})^+ X_j}.
\end{equation}
We show individually that the numerators and denominators in
\eqref{eqlhs} and \eqref{eqrhs} are equal. First the denominators:
starting with \eqref{eqlhs}, notice that
%
%
\begin{eqnarray}\label{eqdenom}
&& s-s_A^T (X_A)^+
X_j - \theta_{j_k,j} \bigl[s_k-s_A^T
(X_A)^+ X_{j_k} \bigr]
\nonumber\\[-8pt]\\[-8pt]
&&\qquad = s-s_{A\cup\{j_k\}}^T
\lleft[ %
\matrix{ (X_A)^+ (X_j -
\theta_{j_k,j} X_{j_k})
\vspace*{3pt}\cr
\theta_{j_k,j}} %
\rright].\nonumber
\end{eqnarray}
By the well-known formula for partial regression coefficients,
\[
\theta_{j_k,j} = \frac{X_{j_k}^T (I-P_A) X_j}{X_{j_k}^T (I-P_A)
X_{j_k}} = \bigl[(X_{A\cup\{j_k\}})^+
X_j \bigr]_{j_k},
\]
that is, $\theta_{j_k,j}$ is the $(j_k)$th coefficient in the
regression of
$X_j$ on $X_{A\cup\{j_k\}}$. Hence, to show that \eqref{eqdenom} is
equal to the denominator in \eqref{eqrhs}, we need to show that
$(X_A)^+(X_j - \theta_{j_k,j} X_{j_k})$ gives the coefficients in
$A$ in the regression of $X_j$ on $X_{A\cup\{j_k\}}$. This follows by
simply noting that the coefficients $(X_{A\cup\{j_k\}})^+X_j =
(\theta_{A,j},\theta_{j_k,j})$ satisfy the equation
\[
X_A \theta_{A,j} + X_{j_k}
\theta_{j_k,j} = P_{A\cup\{j_k\}} X_j
\]
and so solving for $\theta_{A,j}$,
\[
\theta_{A,j} = (X_A)^+ (P_{A\cup\{j_k\}}
X_j - \theta_{j_k,j}X_{j_k}) =
(X_A)^+(X_j-\theta_{j_k,j}X_{j_k}).
\]
Now for the numerators: again beginning with \eqref{eqlhs}, its
numerator is
%
%
\begin{equation}
\label{eqnumer} y^T(I-P_A) (X_j -
\theta_{j_k,j} X_{j_k})
\end{equation}
and by essentially the same argument as above, we have
\[
P_A(X_j-\theta_{j_k,j}X_{j_k}) =
P_{A\cup\{j_k\}} X_j,
\]
therefore, \eqref{eqnumer} matches the numerator in \eqref{eqrhs}.
\end{pf}
\end{appendix}

\section*{\texorpdfstring{Acknowledgements.}{Acknowledgements}}
We thank Jacob Bien, Trevor Hastie, Fred Huffer and Larry Wasserman
for helpful comments.


%

\printaddresses


\begin{thebibliography}{43}
\bibitem[\protect\citeauthoryear{Beck and Teboulle}{2009}]{fista}
%
\begin{barticle}[mr]
\bauthor{\bsnm{Beck},~\bfnm{Amir}\binits{A.}} \AND
\bauthor{\bsnm{Teboulle},~\bfnm{Marc}\binits{M.}}
(\byear{2009}).
\btitle{A fast iterative shrinkage-thresholding algorithm for linear
inverse problems}.
\bjournal{SIAM J. Imaging Sci.}
\bvolume{2}
\bpages{183--202}.
\bid{doi={10.1137/080716542}, issn={1936-4954}, mr={2486527}}
\end{barticle}
%
\bptok{imsref}%
\endbibitem

\bibitem[\protect\citeauthoryear{Becker, Bobin and Cand{\`e}s}{2011}]{nesta}
%
\begin{barticle}[mr]
\bauthor{\bsnm{Becker},~\bfnm{Stephen}\binits{S.}},
\bauthor{\bsnm{Bobin},~\bfnm{J{\'e}r{\^o}me}\binits{J.}} \AND
\bauthor{\bsnm{Cand{\`e}s},~\bfnm{Emmanuel~J.}\binits{E.~J.}}
(\byear{2011}).
\btitle{N{ESTA}: A fast and accurate first-order method for sparse recovery}.
\bjournal{SIAM J. Imaging Sci.}
\bvolume{4}
\bpages{1--39}.
\bid{doi={10.1137/090756855}, issn={1936-4954}, mr={2765668}}
\end{barticle}
%
\bptok{imsref}%
\endbibitem\vadjust{\goodbreak}

\bibitem[\protect\citeauthoryear{Becker, Cand{\`e}s and Grant}{2011}]{tfocs}
%
\begin{barticle}[mr]
\bauthor{\bsnm{Becker},~\bfnm{Stephen~R.}\binits{S.~R.}},
\bauthor{\bsnm{Cand{\`e}s},~\bfnm{Emmanuel~J.}\binits{E.~J.}} \AND
\bauthor{\bsnm{Grant},~\bfnm{Michael~C.}\binits{M.~C.}}
(\byear{2011}).
\btitle{Templates for convex cone problems with applications to sparse
signal recovery}.
\bjournal{Math. Program. Comput.}
\bvolume{3}
\bpages{165--218}.
\bid{doi={10.1007/s12532-011-0029-5}, issn={1867-2949}, mr={2833262}}
\end{barticle}
%
\bptok{imsref}%
\endbibitem

\bibitem[\protect\citeauthoryear{Boyd et~al.}{2011}]{admm}
%
\begin{barticle}[auto:STB|2014/02/12|14:17:21]
\bauthor{\bsnm{Boyd},~\bfnm{S.}\binits{S.}},
\bauthor{\bsnm{Parikh},~\bfnm{N.}\binits{N.}},
\bauthor{\bsnm{Chu},~\bfnm{E.}\binits{E.}},
\bauthor{\bsnm{Peleato},~\bfnm{B.}\binits{B.}} \AND
\bauthor{\bsnm{Eckstein},~\bfnm{J.}\binits{J.}}
(\byear{2011}).
\btitle{Distributed optimization and statistical learning via the
alternative direction method of multipliers}.
\bjournal{Faund. Trends Mach. Learn.}
\bvolume{3}
\bpages{1--122}.
\end{barticle}
%
\bptok{imsref}%
\endbibitem

\bibitem[\protect\citeauthoryear{B{\"u}hlmann}{2013}]{buhlsignif}
%
\begin{barticle}[mr]
\bauthor{\bsnm{B{\"u}hlmann},~\bfnm{Peter}\binits{P.}}
(\byear{2013}).
\btitle{Statistical significance in high-dimensional linear models}.
\bjournal{Bernoulli}
\bvolume{19}
\bpages{1212--1242}.
\bid{doi={10.3150/12-BEJSP11}, issn={1350-7265}, mr={3102549}}
\bptnote{check year}%
\end{barticle}
%
\bptok{imsref}%
\endbibitem

\bibitem[\protect\citeauthoryear{Cand{\`e}s and Plan}{2009}]{nearideal}
%
\begin{barticle}[mr]
\bauthor{\bsnm{Cand{\`e}s},~\bfnm{Emmanuel~J.}\binits{E.~J.}} \AND
\bauthor{\bsnm{Plan},~\bfnm{Yaniv}\binits{Y.}}
(\byear{2009}).
\btitle{Near-ideal model selection by {$\ell_ 1$} minimization}.
\bjournal{Ann. Statist.}
\bvolume{37}
\bpages{2145--2177}.
\bid{doi={10.1214/08-AOS653}, issn={0090-5364}, mr={2543688}}
\end{barticle}
%
\bptok{imsref}%
\endbibitem

\bibitem[\protect\citeauthoryear{Candes and Tao}{2006}]{nearopt}
%
\begin{barticle}[mr]
\bauthor{\bsnm{Candes},~\bfnm{Emmanuel~J.}\binits{E.~J.}} \AND
\bauthor{\bsnm{Tao},~\bfnm{Terence}\binits{T.}}
(\byear{2006}).
\btitle{Near-optimal signal recovery from random projections:
Universal encoding strategies?}
\bjournal{IEEE Trans. Inform. Theory}
\bvolume{52}
\bpages{5406--5425}.
\bid{doi={10.1109/TIT.2006.885507}, issn={0018-9448}, mr={2300700}}
\end{barticle}
%
\bptok{imsref}%
\endbibitem

\bibitem[\protect\citeauthoryear{Chen, Donoho and Saunders}{1998}]{bp}
%
\begin{barticle}[mr]
\bauthor{\bsnm{Chen},~\bfnm{Scott~Shaobing}\binits{S.~S.}},
\bauthor{\bsnm{Donoho},~\bfnm{David~L.}\binits{D.~L.}} \AND
\bauthor{\bsnm{Saunders},~\bfnm{Michael~A.}\binits{M.~A.}}
(\byear{1998}).
\btitle{Atomic decomposition by basis pursuit}.
\bjournal{SIAM J. Sci. Comput.}
\bvolume{20}
\bpages{33--61}.
\bid{doi={10.1137/S1064827596304010}, issn={1064-8275}, mr={1639094}}
\end{barticle}
%
\bptok{imsref}%
\endbibitem

\bibitem[\protect\citeauthoryear{de~Haan and Ferreira}{2006}]{evt}
%
\begin{bbook}[mr]
\bauthor{\bparticle{de} \bsnm{Haan},~\bfnm{Laurens}\binits{L.}}
\AND
\bauthor{\bsnm{Ferreira},~\bfnm{Ana}\binits{A.}}
(\byear{2006}).
\btitle{Extreme Value Theory: An Introduction}.
\bpublisher{Springer},
\blocation{New York}.
\bid{mr={2234156}}
\end{bbook}
%
\bptok{imsref}%
\endbibitem

\bibitem[\protect\citeauthoryear{Donoho}{2006}]{cs}
%
\begin{barticle}[mr]
\bauthor{\bsnm{Donoho},~\bfnm{David~L.}\binits{D.~L.}}
(\byear{2006}).
\btitle{Compressed sensing}.
\bjournal{IEEE Trans. Inform. Theory}
\bvolume{52}
\bpages{1289--1306}.
\bid{doi={10.1109/TIT.2006.871582}, issn={0018-9448}, mr={2241189}}
\end{barticle}
%
\bptok{imsref}%
\endbibitem

\bibitem[\protect\citeauthoryear{Efron}{1986}]{bradbiased}
%
\begin{barticle}[mr]
\bauthor{\bsnm{Efron},~\bfnm{Bradley}\binits{B.}}
(\byear{1986}).
\btitle{How biased is the apparent error rate of a prediction rule?}
\bjournal{J. Amer. Statist. Assoc.}
\bvolume{81}
\bpages{461--470}.
\bid{issn={0162-1459}, mr={0845884}}
\end{barticle}
%
\bptok{imsref}%
\endbibitem

\bibitem[\protect\citeauthoryear{Efron et~al.}{2004}]{lars}
%
\begin{barticle}[mr]
\bauthor{\bsnm{Efron},~\bfnm{Bradley}\binits{B.}},
\bauthor{\bsnm{Hastie},~\bfnm{Trevor}\binits{T.}},
\bauthor{\bsnm{Johnstone},~\bfnm{Iain}\binits{I.}} \AND
\bauthor{\bsnm{Tibshirani},~\bfnm{Robert}\binits{R.}}
(\byear{2004}).
\btitle{Least angle regression}.
\bjournal{Ann. Statist.}
\bvolume{32}
\bpages{407--499}.
\bid{doi={10.1214/009053604000000067}, issn={0090-5364}, mr={2060166}}
\bptnote{check related}%
\end{barticle}
%
\bptok{imsref}%
\endbibitem

\bibitem[\protect\citeauthoryear{Fan, Guo and Hao}{2012}]{varestbycv}
%
\begin{barticle}[mr]
\bauthor{\bsnm{Fan},~\bfnm{Jianqing}\binits{J.}},
\bauthor{\bsnm{Guo},~\bfnm{Shaojun}\binits{S.}} \AND
\bauthor{\bsnm{Hao},~\bfnm{Ning}\binits{N.}}
(\byear{2012}).
\btitle{Variance estimation using refitted cross-validation in
ultrahigh-dimensional regression}.
\bjournal{J. R. Stat. Soc. Ser. B Stat. Methodol.}
\bvolume{74}
\bpages{37--65}.
\bid{doi={10.1111/j.1467-9868.2011.01005.x}, issn={1369-7412}, mr={2885839}}
\end{barticle}
%
\bptok{imsref}%
\endbibitem

\bibitem[\protect\citeauthoryear{Friedman, Hastie and
Tibshirani}{2010}]{glmnet}
%
\begin{barticle}[auto:STB|2014/02/12|14:17:21]
\bauthor{\bsnm{Friedman},~\bfnm{J.}\binits{J.}},
\bauthor{\bsnm{Hastie},~\bfnm{T.}\binits{T.}} \AND
\bauthor{\bsnm{Tibshirani},~\bfnm{R.}\binits{R.}}
(\byear{2010}).
\btitle{Regularization paths for generalized linear models via
coordinate descent}.
\bjournal{J. Stat. Softw.}
\bvolume{33}
\bpages{1--22}.
\end{barticle}
%
\bptok{imsref}%
\endbibitem

\bibitem[\protect\citeauthoryear{Friedman et~al.}{2007}]{pco}
%
\begin{barticle}[mr]
\bauthor{\bsnm{Friedman},~\bfnm{Jerome}\binits{J.}},
\bauthor{\bsnm{Hastie},~\bfnm{Trevor}\binits{T.}},
\bauthor{\bsnm{H{\"o}fling},~\bfnm{Holger}\binits{H.}} \AND
\bauthor{\bsnm{Tibshirani},~\bfnm{Robert}\binits{R.}}
(\byear{2007}).
\btitle{Pathwise coordinate optimization}.
\bjournal{Ann. Appl. Stat.}
\bvolume{1}
\bpages{302--332}.
\bid{doi={10.1214/07-AOAS131}, issn={1932-6157}, mr={2415737}}
\end{barticle}
%
\bptok{imsref}%
\endbibitem

\bibitem[\protect\citeauthoryear{Fuchs}{2005}]{fuchs}
%
\begin{barticle}[mr]
\bauthor{\bsnm{Fuchs},~\bfnm{Jean~Jacques}\binits{J.~J.}}
(\byear{2005}).
\btitle{Recovery of exact sparse representations in the presence of
bounded noise}.
\bjournal{IEEE Trans. Inform. Theory}
\bvolume{51}
\bpages{3601--3608}.
\bid{doi={10.1109/TIT.2005.855614}, issn={0018-9448}, mr={2237526}}
\end{barticle}
%
\bptok{imsref}%
\endbibitem

\bibitem[\protect\citeauthoryear{Grazier~G'Sell, Taylor and Tibshirani}{2013}]{glassosignif}
%
\begin{bmisc}[auto:STB|2014/02/12|14:17:21]
\bauthor{\bsnm{Grazier G'Sell},~\bfnm{M.}\binits{M.}},
\bauthor{\bsnm{Taylor},~\bfnm{J.}\binits{J.}} \AND
\bauthor{\bsnm{Tibshirani},~\bfnm{R.}\binits{R.}}
(\byear{2013}).
\bhowpublished{Adaptive testing for the graphical lasso.
Preprint. Available at \arxivurl{arXiv:1307.4765}.}
\end{bmisc}
%
\bptok{imsref}%
\endbibitem

\bibitem[\protect\citeauthoryear{Grazier~G'Sell et~al.}{2013}]{fdrlasso}
%
\begin{bmisc}[auto:STB|2014/02/12|14:17:21]
\bauthor{\bsnm{Grazier G'Sell},~\bfnm{M.}\binits{M.}},
\bauthor{\bsnm{Wager},~\bfnm{S.}\binits{S.}},
\bauthor{\bsnm{Chouldechova},~\bfnm{A.}\binits{A.}} \AND
\bauthor{\bsnm{Tibshirani},~\bfnm{R.}\binits{R.}}
(\byear{2013}).
\bhowpublished{False discovery rate control for sequential selection
procedures, with application to the lasso.
Preprint. Available at \arxivurl{arXiv:1309.5352}.}
\end{bmisc}
%
\bptok{imsref}%
\endbibitem

\bibitem[\protect\citeauthoryear{Greenshtein and Ritov}{2004}]{persistence}
%
\begin{barticle}[mr]
\bauthor{\bsnm{Greenshtein},~\bfnm{Eitan}\binits{E.}} \AND
\bauthor{\bsnm{Ritov},~\bfnm{Ya'acov}\binits{Y.}}
(\byear{2004}).
\btitle{Persistence in high-dimensional linear predictor selection and
the virtue of overparametrization}.
\bjournal{Bernoulli}
\bvolume{10}
\bpages{971--988}.
\bid{doi={10.3150/bj/1106314846}, issn={1350-7265}, mr={2108039}}
\end{barticle}
%
\bptok{imsref}%
\endbibitem

\bibitem[\protect\citeauthoryear{Hastie, Tibshirani and Friedman}{2008}]{esl}
%
\begin{bbook}[auto]
\bauthor{\bsnm{Hastie},~\bfnm{Trevor}\binits{T.}},
\bauthor{\bsnm{Tibshirani},~\bfnm{Robert}\binits{R.}} \AND
\bauthor{\bsnm{Friedman},~\bfnm{Jerome}\binits{J.}}
(\byear{2008}).
\btitle{The Elements of Statistical Learning; Data Mining, Inference,
and Prediction},
\bedition{2nd} ed.
\bpublisher{Springer},
\blocation{New York}.
\bid{mr={2722294}}
\end{bbook}
%
\bptok{imsref}%
\endbibitem

\bibitem[\protect\citeauthoryear{Javanmard and Montanari}{2013a}]{montahypo2}
%
\begin{bmisc}[auto:STB|2014/02/12|14:17:21]
\bauthor{\bsnm{Javanmard},~\bfnm{A.}\binits{A.}} \AND
\bauthor{\bsnm{Montanari},~\bfnm{A.}\binits{A.}}
(\byear{2013a}).
\bhowpublished{Confidence intervals and hypothesis testing for
high-dimensional regression.
Preprint. Available at \arxivurl{arXiv:1306.3171}.}
\end{bmisc}
%
\bptok{imsref}%
\endbibitem

\bibitem[\protect\citeauthoryear{Javanmard and Montanari}{2013b}]{montahypo1}
%
\begin{bmisc}[auto:STB|2014/02/12|14:17:21]
\bauthor{\bsnm{Javanmard},~\bfnm{A.}\binits{A.}} \AND
\bauthor{\bsnm{Montanari},~\bfnm{A.}\binits{A.}}
(\byear{2013b}).
\bhowpublished{Hypothesis testing in high-dimensional regression under
the Gaussian random design model: Asymptotic theory.
Preprint. Available at \arxivurl{arXiv:1301.4240}.}
\end{bmisc}
%
\bptok{imsref}%
\endbibitem

\bibitem[\protect\citeauthoryear{Meinshausen and B{\"
u}hlmann}{2010}]{stabselect}
%
\begin{barticle}[mr]
\bauthor{\bsnm{Meinshausen},~\bfnm{Nicolai}\binits{N.}} \AND
\bauthor{\bsnm{B{\"u}hlmann},~\bfnm{Peter}\binits{P.}}
(\byear{2010}).
\btitle{Stability selection}.
\bjournal{J. R. Stat. Soc. Ser. B Stat. Methodol.}
\bvolume{72}
\bpages{417--473}.
\bid{doi={10.1111/j.1467-9868.2010.00740.x}, issn={1369-7412}, mr={2758523}}
\end{barticle}
%
\bptok{imsref}%
\endbibitem

\bibitem[\protect\citeauthoryear{Meinshausen, Meier and B{\"
u}hlmann}{2009}]{mein2009}
%
\begin{barticle}[mr]
\bauthor{\bsnm{Meinshausen},~\bfnm{Nicolai}\binits{N.}},
\bauthor{\bsnm{Meier},~\bfnm{Lukas}\binits{L.}} \AND
\bauthor{\bsnm{B{\"u}hlmann},~\bfnm{Peter}\binits{P.}}
(\byear{2009}).
\btitle{{$p$}-values for high-dimensional regression}.
\bjournal{J. Amer. Statist. Assoc.}
\bvolume{104}
\bpages{1671--1681}.
\bid{doi={10.1198/jasa.2009.tm08647}, issn={0162-1459}, mr={2750584}}
\end{barticle}
%
\bptok{imsref}%
\endbibitem

\bibitem[\protect\citeauthoryear{Minnier, Tian and Cai}{2011}]{minnier2011}
%
\begin{barticle}[mr]
\bauthor{\bsnm{Minnier},~\bfnm{Jessica}\binits{J.}},
\bauthor{\bsnm{Tian},~\bfnm{Lu}\binits{L.}} \AND
\bauthor{\bsnm{Cai},~\bfnm{Tianxi}\binits{T.}}
(\byear{2011}).
\btitle{A perturbation method for inference on regularized regression
estimates}.
\bjournal{J. Amer. Statist. Assoc.}
\bvolume{106}
\bpages{1371--1382}.
\bid{doi={10.1198/jasa.2011.tm10382}, issn={0162-1459}, mr={2896842}}
\end{barticle}
%
\bptok{imsref}%
\endbibitem

\bibitem[\protect\citeauthoryear{Osborne, Presnell and
Turlach}{2000a}]{homotopy1}
%
\begin{barticle}[mr]
\bauthor{\bsnm{Osborne},~\bfnm{M.~R.}\binits{M.~R.}},
\bauthor{\bsnm{Presnell},~\bfnm{Brett}\binits{B.}} \AND
\bauthor{\bsnm{Turlach},~\bfnm{B.~A.}\binits{B.~A.}}
(\byear{2000}a).
\btitle{A new approach to variable selection in least squares problems}.
\bjournal{IMA J. Numer. Anal.}
\bvolume{20}
\bpages{389--403}.
\bid{doi={10.1093/imanum/20.3.389}, issn={0272-4979}, mr={1773265}}
\end{barticle}
%
\bptok{imsref}%
\endbibitem

\bibitem[\protect\citeauthoryear{Osborne, Presnell and
Turlach}{2000b}]{homotopy2}
%
\begin{barticle}[mr]
\bauthor{\bsnm{Osborne},~\bfnm{Michael~R.}\binits{M.~R.}},
\bauthor{\bsnm{Presnell},~\bfnm{Brett}\binits{B.}} \AND
\bauthor{\bsnm{Turlach},~\bfnm{Berwin~A.}\binits{B.~A.}}
(\byear{2000}b).
\btitle{On the {LASSO} and its dual}.
\bjournal{J. Comput. Graph. Statist.}
\bvolume{9}
\bpages{319--337}.
\bid{doi={10.2307/1390657}, issn={1061-8600}, mr={1822089}}
\end{barticle}
%
\bptok{imsref}%
\endbibitem\vadjust{\goodbreak}

\bibitem[\protect\citeauthoryear{Park and Hastie}{2007}]{glmpath}
%
\begin{barticle}[mr]
\bauthor{\bsnm{Park},~\bfnm{Mee~Young}\binits{M.~Y.}} \AND
\bauthor{\bsnm{Hastie},~\bfnm{Trevor}\binits{T.}}
(\byear{2007}).
\btitle{{$L_ 1$}-regularization path algorithm for generalized linear models}.
\bjournal{J. R. Stat. Soc. Ser. B Stat. Methodol.}
\bvolume{69}
\bpages{659--677}.
\bid{doi={10.1111/j.1467-9868.2007.00607.x}, issn={1369-7412}, mr={2370074}}
\end{barticle}
%
\bptok{imsref}%
\endbibitem

\bibitem[\protect\citeauthoryear{Rhee et~al.}{2003}]{rhee2003}
%
\begin{barticle}[auto:STB|2014/02/12|14:17:21]
\bauthor{\bsnm{Rhee},~\bfnm{S.-Y.}\binits{S.-Y.}},
\bauthor{\bsnm{Gonzales},~\bfnm{M.~J.}\binits{M.~J.}},
\bauthor{\bsnm{Kantor},~\bfnm{R.}\binits{R.}},
\bauthor{\bsnm{Betts},~\bfnm{B.~J.}\binits{B.~J.}},
\bauthor{\bsnm{Ravela},~\bfnm{J.}\binits{J.}} \AND
\bauthor{\bsnm{Shafer},~\bfnm{R.~W.}\binits{R.~W.}}
(\byear{2003}).
\btitle{Human immunodeficiency virus reverse transcriptase and
protease sequence database}.
\bjournal{Nucleic Acids Res.}
\bvolume{31}
\bpages{298--303}.
\end{barticle}
%
\bptok{imsref}%
\endbibitem

\bibitem[\protect\citeauthoryear{Sun and Zhang}{2012}]{scaledlasso}
%
\begin{barticle}[mr]
\bauthor{\bsnm{Sun},~\bfnm{Tingni}\binits{T.}} \AND
\bauthor{\bsnm{Zhang},~\bfnm{Cun-Hui}\binits{C.-H.}}
(\byear{2012}).
\btitle{Scaled sparse linear regression}.
\bjournal{Biometrika}
\bvolume{99}
\bpages{879--898}.
\bid{doi={10.1093/biomet/ass043}, issn={0006-3444}, mr={2999166}}
\end{barticle}
%
\bptok{imsref}%
\endbibitem

\bibitem[\protect\citeauthoryear{Taylor, Loftus and Tibshirani}{2013}]{geomsignif}
%
\begin{bmisc}[auto:STB|2014/02/12|14:17:21]
\bauthor{\bsnm{Taylor},~\bfnm{J.}\binits{J.}},
\bauthor{\bsnm{Loftus},~\bfnm{J.}\binits{J.}} \AND
\bauthor{\bsnm{Tibshirani},~\bfnm{R.~J.}\binits{R.~J.}}
(\byear{2013}).
\bhowpublished{Tests in adaptive regression via the Kac--{R}ice formula.
Preprint. Available at \arxivurl{arXiv:1308.3020}.}
\end{bmisc}
%
\bptok{imsref}%
\endbibitem

\bibitem[\protect\citeauthoryear{Taylor, Takemura and
Adler}{2005}]{taylorvalid}
%
\begin{barticle}[mr]
\bauthor{\bsnm{Taylor},~\bfnm{Jonathan}\binits{J.}},
\bauthor{\bsnm{Takemura},~\bfnm{Akimichi}\binits{A.}} \AND
\bauthor{\bsnm{Adler},~\bfnm{Robert~J.}\binits{R.~J.}}
(\byear{2005}).
\btitle{Validity of the expected {E}uler characteristic heuristic}.
\bjournal{Ann. Probab.}
\bvolume{33}
\bpages{1362--1396}.
\bid{doi={10.1214/009117905000000099}, issn={0091-1798}, mr={2150192}}
\end{barticle}
%
\bptok{imsref}%
\endbibitem

\bibitem[\protect\citeauthoryear{Tibshirani}{1996}]{lasso}
%
\begin{barticle}[mr]
\bauthor{\bsnm{Tibshirani},~\bfnm{Robert}\binits{R.}}
(\byear{1996}).
\btitle{Regression shrinkage and selection via the lasso}.
\bjournal{J. Roy. Statist. Soc. Ser. B}
\bvolume{58}
\bpages{267--288}.
\bid{issn={0035-9246}, mr={1379242}}
\end{barticle}
%
\bptok{imsref}%
\endbibitem

\bibitem[\protect\citeauthoryear{Tibshirani}{2013}]{lassounique}
%
\begin{barticle}[mr]
\bauthor{\bsnm{Tibshirani},~\bfnm{Ryan~J.}\binits{Ryan~J.}}
(\byear{2013}).
\btitle{The lasso problem and uniqueness}.
\bjournal{Electron. J. Stat.}
\bvolume{7}
\bpages{1456--1490}.
\bid{doi={10.1214/13-EJS815}, issn={1935-7524}, mr={3066375}}
\bptnote{check year}%
\end{barticle}
%
\bptok{imsref}%
\endbibitem

\bibitem[\protect\citeauthoryear{Tibshirani and Taylor}{2012}]{lassodf2}
%
\begin{barticle}[mr]
\bauthor{\bsnm{Tibshirani},~\bfnm{Ryan~J.}\binits{R.~J.}} \AND
\bauthor{\bsnm{Taylor},~\bfnm{Jonathan}\binits{J.}}
(\byear{2012}).
\btitle{Degrees of freedom in lasso problems}.
\bjournal{Ann. Statist.}
\bvolume{40}
\bpages{1198--1232}.
\bid{doi={10.1214/12-AOS1003}, issn={0090-5364}, mr={2985948}}
\end{barticle}
%
\bptok{imsref}%
\endbibitem

\bibitem[\protect\citeauthoryear{van~de Geer and B{\"u}hlmann}{2013}]{vdgsignif}
%
\begin{bmisc}[auto]
\bauthor{\bsnm{van~de Geer},~\bfnm{Sara}\binits{S.}} \AND
\bauthor{\bsnm{B{\"u}hlmann},~\bfnm{Peter}\binits{P.}}
(\byear{2013}).
\bhowpublished{On asymptotically optimal confidence regions and tests
for high-dimensional models.
Preprint. Available at \arxivurl{arXiv:1303.0518}.}
\end{bmisc}
%
\bptok{imsref}%
\endbibitem
\bibitem[\protect\citeauthoryear{Wainwright}{2009}]{sharp}
%
\begin{barticle}[mr]
\bauthor{\bsnm{Wainwright},~\bfnm{Martin~J.}\binits{M.~J.}}
(\byear{2009}).
\btitle{Sharp thresholds for high-dimensional and noisy sparsity
recovery using {$\ell_ 1$}-constrained quadratic programming ({L}asso)}.
\bjournal{IEEE Trans. Inform. Theory}
\bvolume{55}
\bpages{2183--2202}.
\bid{doi={10.1109/TIT.2009.2016018}, issn={0018-9448}, mr={2729873}}
\end{barticle}
%
\bptok{imsref}%
\endbibitem

\bibitem[\protect\citeauthoryear{Wasserman and Roeder}{2009}]{screenclean}
%
\begin{barticle}[mr]
\bauthor{\bsnm{Wasserman},~\bfnm{Larry}\binits{L.}} \AND
\bauthor{\bsnm{Roeder},~\bfnm{Kathryn}\binits{K.}}
(\byear{2009}).
\btitle{High-dimensional variable selection}.
\bjournal{Ann. Statist.}
\bvolume{37}
\bpages{2178--2201}.
\bid{doi={10.1214/08-AOS646}, issn={0090-5364}, mr={2543689}}
\end{barticle}
%
\bptok{imsref}%
\endbibitem

\bibitem[\protect\citeauthoryear{Weissman}{1978}]{weissman}
%
\begin{barticle}[mr]
\bauthor{\bsnm{Weissman},~\bfnm{Ishay}\binits{I.}}
(\byear{1978}).
\btitle{Estimation of parameters and large quantiles based on the
{$k$} largest observations}.
\bjournal{J. Amer. Statist. Assoc.}
\bvolume{73}
\bpages{812--815}.
\bid{issn={0003-1291}, mr={0521329}}
\end{barticle}
%
\bptok{imsref}%
\endbibitem

\bibitem[\protect\citeauthoryear{Zhang and Zhang}{2014}]{zhangconf}
%
\begin{barticle}[auto:STB|2014/02/12|14:17:21]
\bauthor{\bsnm{Zhang},~\bfnm{C.-H.}\binits{C.-H.}} \AND
\bauthor{\bsnm{Zhang},~\bfnm{S.}\binits{S.}}
(\byear{2014}).
\btitle{Confidence intervals for low dimensional parameters in high dimensional linear models}.
\bjournal{J. R. Stat. Soc. Ser. B Stat. Methodol.}
\bvolume{76}
\bpages{217--242}.
\bid{mr={3153940}}
\end{barticle}
%
\bptok{imsref}%
\endbibitem

\bibitem[\protect\citeauthoryear{Zhao and Yu}{2006}]{irrepcond}
%
\begin{barticle}[mr]
\bauthor{\bsnm{Zhao},~\bfnm{Peng}\binits{P.}} \AND
\bauthor{\bsnm{Yu},~\bfnm{Bin}\binits{B.}}
(\byear{2006}).
\btitle{On model selection consistency of {L}asso}.
\bjournal{J. Mach. Learn. Res.}
\bvolume{7}
\bpages{2541--2563}.
\bid{issn={1532-4435}, mr={2274449}}
\end{barticle}
%
\bptok{imsref}%
\endbibitem

\bibitem[\protect\citeauthoryear{Zou and Hastie}{2005}]{enet}
%
\begin{barticle}[mr]
\bauthor{\bsnm{Zou},~\bfnm{Hui}\binits{H.}} \AND
\bauthor{\bsnm{Hastie},~\bfnm{Trevor}\binits{T.}}
(\byear{2005}).
\btitle{Regularization and variable selection via the elastic net}.
\bjournal{J. R. Stat. Soc. Ser. B Stat. Methodol.}
\bvolume{67}
\bpages{301--320}.
\bid{doi={10.1111/j.1467-9868.2005.00503.x}, issn={1369-7412}, mr={2137327}}
\end{barticle}
%
\bptok{imsref}%
\endbibitem

\bibitem[\protect\citeauthoryear{Zou, Hastie and Tibshirani}{2007}]{lassodf}
%
\begin{barticle}[mr]
\bauthor{\bsnm{Zou},~\bfnm{Hui}\binits{H.}},
\bauthor{\bsnm{Hastie},~\bfnm{Trevor}\binits{T.}} \AND
\bauthor{\bsnm{Tibshirani},~\bfnm{Robert}\binits{R.}}
(\byear{2007}).
\btitle{On the ``degrees of freedom'' of the lasso}.
\bjournal{Ann. Statist.}
\bvolume{35}
\bpages{2173--2192}.
\bid{doi={10.1214/009053607000000127}, issn={0090-5364}, mr={2363967}}
\end{barticle}
%
\bptok{imsref}%
\endbibitem

\end{thebibliography}
\end{document}